\documentclass[draft]{amsart} 
\usepackage{amssymb,latexsym,amsfonts,verbatim,amscd,ifthen} 
\usepackage{color}

\numberwithin{equation}{section}

\theoremstyle{plain}

\newtheorem{theorem}[equation]{Theorem}   
\newtheorem{lemma}[equation]{Lemma} 
\newtheorem{proposition}[equation]{Proposition} 
\newtheorem{corollary}[equation]{Corollary}

\theoremstyle{definition}

\newtheorem{example}[equation]{Example} 
\newtheorem{definition}[equation]{Definition} 
\newtheorem{remark}[equation]{Remark} 
\newtheorem{remarks}[equation]{Remarks}

\DeclareMathOperator{\HH}{H}

\DeclareMathOperator{\Coker}{Coker}

\DeclareMathOperator{\rank}{rank}
\DeclareMathOperator{\im}{Im} 
\DeclareMathOperator{\pd}{pd}

\DeclareMathOperator{\Ker}{Ker} 
\DeclareMathOperator{\lcm}{lcm}

\begin{document} 

\renewcommand{\:}{\! :\ } 
\newcommand{\p}{\mathfrak p} 
\newcommand{\m}{\mathfrak m}
\newcommand{\e}{\epsilon} 
\newcommand{\g}{\gamma} 
\newcommand{\lra}{\longrightarrow} 
\newcommand{\ra}{\rightarrow} 
\newcommand{\altref}[1]{{\upshape(\ref{#1})}} 
\newcommand{\bfa}{\boldsymbol{\alpha}} 
\newcommand{\bfb}{\boldsymbol{\beta}} 
\newcommand{\bfg}{\boldsymbol{\gamma}} 
\newcommand{\bfd}{\boldsymbol{\delta}} 
\newcommand{\bfM}{\mathbf M} 
\newcommand{\bfN}{\mathbf N}
\newcommand{\bfI}{\mathbf I} 
\newcommand{\bfC}{\mathbf C} 
\newcommand{\bfB}{\mathbf B} 
\newcommand{\bsfC}{\bold{\mathsf C}} 
\newcommand{\bsfT}{\bold{\mathsf T}}
\newcommand{\mc}{\mathcal} 
\newcommand{\smsm}{\smallsetminus} 
\newcommand{\ol}{\overline} 
\newcommand{\twedge}
           {\smash{\overset{\mbox{}_{\circ}}
                           {\wedge}}\thinspace} 

\newlength{\wdtha}
\newlength{\wdthb}
\newlength{\wdthc}
\newlength{\wdthd}
\newcommand{\elabel}[1]
           {\label{#1}  
            \setlength{\wdtha}{.4\marginparwidth}
            \settowidth{\wdthb}{\tt\small{#1}} 
            \addtolength{\wdthb}{\wdtha}
            \smash{
            \raisebox{.8\baselineskip}
            {\color{red}
             \hspace*{-\wdthb}\tt\small{#1}\hspace{\wdtha}}}}

\newcommand{\mlabel}[1] 
           {\label{#1} 
            \setlength{\wdtha}{\textwidth}
            \setlength{\wdthb}{\wdtha} 
            \addtolength{\wdthb}{\marginparsep} 
            \addtolength{\wdthb}{\marginparwidth}
            \setlength{\wdthc}{\marginparwidth}
            \setlength{\wdthd}{\marginparsep}
            \addtolength{\wdtha}{2\wdthc}
            \addtolength{\wdtha}{2\marginparsep} 
            \setlength{\marginparwidth}{\wdtha}
            \setlength{\marginparsep}{-\wdthb} 
            \setlength{\wdtha}{\wdthc} 
            \addtolength{\wdtha}{1.1ex}
            \marginpar{\vspace*{-0.3\baselineskip}
                       \tt\small{#1}\\[-0.4\baselineskip]\rule{\wdtha}{.5pt} }
            \setlength{\marginparwidth}{\wdthc} 
            \setlength{\marginparsep}{\wdthd}  }  
            
\renewcommand{\mlabel}{\label} 
\renewcommand{\elabel}{\label}

\title[Matroids and free resolutions]
      {Representations of matroids  
       and free resolutions \\ 
       for multigraded modules}
\author[A. Tchernev]{Alexandre B. Tchernev} 
\address{Department of Mathematics\\
         University at Albany, SUNY\\ 
         Albany, NY 12222}
\email{tchernev@math.albany.edu}
\begin{abstract} 
Let $\Bbbk$ be a field, let $R=\Bbbk[x_1,\dots, x_m]$ be 
a polynomial ring with the standard $\mathbb Z^m$-grading 
(multigrading), 
let $L$ be a Noetherian multigraded $R$-module, and let  
$E \smash{\overset \Phi \lra}\ G \rightarrow L \rightarrow 0$ 
\ be a finite free multigraded presentation of $L$ over $R$. Given 
a choice $S$ of a multihomogeneous basis of $E$, we construct 
an explicit canonical finite free multigraded resolution 
$T_\bullet(\Phi, S)$ of the $R$-module $L$. In the case 
of monomial ideals our construction recovers the Taylor 
resolution. A main ingredient of our work 
is a new linear algebra construction of independent 
interest, which produces from 
a representation $\phi$ over $\Bbbk$ of a matroid $\bfM$ 
a canonical finite complex of finite 
dimensional $\Bbbk$-vector spaces $T_\bullet(\phi)$ that 
is a resolution of $\Ker \phi$. We also show that 
the length of $T_\bullet(\phi)$ and the dimensions of 
its components are combinatorial invariants of 
the matroid $\bfM$, and are independent of the representation 
map $\phi$.    
\end{abstract} 

\maketitle 

\tableofcontents 

\section*{Introduction} 

In her thesis \cite{Ta66}, Diana Taylor constructs 
an explicit finite free resolution for an arbitrary 
monomial ideal. Even though the Taylor resolution 
is generally nonminimal, its canonical combinatorial description 
has turned this construction into one of the main tools  
in the study of the homological 
structure and combinatorics of monomial 
ideals. It is therefore a natural and important question 
whether a similarly explicit construction 
exists in the more general setting of multigraded modules.  
The difficulty of studying resolutions of mutigraded modules 
is apparent from the fact that most of the published work 
in this area deals exclusively with  
homological properties of monomial ideals, and   
until recently there were only a few 
general results available 
on the structure of free resolutions of 
multigraded modules that are not cyclic:  
Lescot \cite{Le88} shows that the multigraded Poincar\'e 
series of a multigraded module is rational, the papers 
of Bruns and Herzog \cite{BrHe95} and Iyengar \cite{Iy97} 
provide bounds for the shifts in a minimal multigraded 
resolution, and the papers of 
Santoni~\cite{Sa90},  
Charalambous~\cite{Ch91}, and R\"omer~\cite{Ro02}  
provide bounds for the Betti numbers of a multigraded 
module. In the last five years however 
there has been an increased interest in the study of 
the homological properties of multigraded modules. 
The papers of Yanagawa~\cite{Ya00, Ya01, Ya04} introduce 
and study the class of squarefree multigraded modules. 
The work of R\" omer~\cite{Ro01} introduced 
the Alexander duality functor for squarefree modules;  
at the same time  and independently 
Miller~\cite{Mi00} defined Alexander duality for 
arbitrary multigraded modules and used his breakthrough to 
exhibit the fundamental correspondence between minimal 
free and minimal injective resolutions of multigraded 
modules. These results provide a wonderful insight   
into the homological structure of multigraded modules, and 
they also make the problem of finding an explicit 
construction for free resolutions even more interesting 
because a solution would automatically provide 
via Alexander duality an explicit construction for injective 
resolutions. Along the lines of explicit constructions, 
Bayer and Sturmfels~\cite{BaSt98} provide a resolution 
for infinitely generated co-Artinian monomial modules,  
Charalambous and Deno~\cite{ChDe01} compute the module of 
second syzygies in a minimal free presentation of a Noetherian 
multigraded module, while Charalambous and Tchernev~\cite{ChTc03} 
introduce the notion of a generic Noetherian multigraded 
module and construct explicitly the entire minimal free 
resolution of such a module. 

The main goal of this paper is to address 
the still remaining question on the existence of an 
explicit free resolution for an arbitrary Noetherian multigraded 
module, and to answer it in  the affirmative. 
The principal tool for achieving such an explicit 
description 
is a new construction of a complex of vector spaces that 
lies at the crossroads of matroid theory and 
multilinear algebra,  
and has the potential for applications 
to both of these fields, as well as to 
other areas such as group representation theory. 
The properties of this new complex  
give a strong indication that there are 
two main combinatorial 
ingredients determining the structure of the 
minimal free resolution  
of a Noetherian multigraded module $L$.  
The first ingredient is, as expected 
from the case of monomial ideals, the combinatorics of the 
multidegrees of a (fixed) set of minimal 
generators of the first syzygy module of $L$. 
The second ingredient, which is new but probably not very 
surprising, are the combinatorial properties of the matroid 
asociated with this fixed 
set of minimal generators. What is surprising however, 
is that our construction makes it possible 
to obtain an explicit 
description of the interaction between these two ingredients; 
a detailed investigation of this interaction 
in a generic situation is carried out in~\cite{ChTc}. 

To describe the results of the present article 
in more detail, let $\Bbbk$ be a field, let $W$ be a finite 
dimensional $\Bbbk$-vector space, let $S$ be a finite set, 
let $U_S$ be the $\Bbbk$-vector space with basis the 
set of symbols $\{e_a\mid a\in S\}$, 
and let $\phi\: U_S \lra W$ be a vector space homomorphism. 
We construct in a canonical way a finite complex 
$T_\bullet(\phi)$ of finite dimensional $\Bbbk$-vector spaces 
that we call \emph{the T-complex of $\phi$}. 
The length of this complex  and the ranks of the vector spaces 
in it depend only on the matroid $\bfM$ on $S$ that is  
represented by the map $\phi$, and not on the actual map $\phi$.   
The T-complex $T_\bullet(\phi)$ behaves very well with respect to 
the operations contraction and restriction of matroids, and the 
first main result of this paper, 
Theorem~\ref{T:exactness-of-complex}, states that this complex 
is a resolution of the kernel $\Ker(\phi)$. 

In the process of defining 
$T_\bullet(\phi)$ we introduce for each 
subset $A\subseteq S$ a certain \emph{multiplicity} 
vector space $S_A(\phi)$, which appears to capture 
subtle combinatorial properties of the matroid $\bfM$, and 
seems to be of independent interest. For example, we show in 
Theorem~\ref{T:multiplicity-zero} that $S_A(\phi)\ne 0$ 
if and only if $A$ is a connected dependent set in $\bfM$. 
Furthermore, a consequence of Theorem~\ref{T:exactness-of-complex} 
is that the dimension of $S_A(\phi)$ is an invariant of the matroid 
$\bfM$, and is independent of the presentation $\phi$. The 
space $S_A(\phi)$ has a natural system of generators indexed 
by a certain collection $C_{\bfM}(A)$ of maximal chains of flats 
in the dual matroid $\bfM^*$, and it is a very interesting 
open question whether the corresponding representable 
matroid on the set $C_{\bfM}(A)$ is independent of the 
representation map $\phi$. 

Next, let $R=\Bbbk[x_1,\dots, x_m]$ be a polynomial ring 
with the usual $\mathbb Z^m$-grading (multigrading), let 
$L$ be a multigraded Noetherian $R$-module, let 
\[
E \overset{\Phi}\lra G \lra L \lra 0
\]
be a finite free (not necessarily minimal) multigraded 
presentation of $L$ over $R$, and let $S$ be a multihomogeneous 
basis of $E$. Consider $\Bbbk$ as an $R$-module via the 
canonical projection $R\lra \Bbbk$ that sends each variable 
$x_i$ to the identity element $1\in \Bbbk$, 
and let $\phi= \Bbbk \otimes_R \Phi$.  
By a standard procedure we use the multidegrees 
of the elements in $S$ to transform the complex of vector 
spaces   $T_\bullet(\phi)$ into a complex of multigraded free 
$R$-modules 
$T_\bullet(\Phi, S)$ that we call \emph{the T-resolution 
of the pair $(\Phi, S)$}. This complex is non-minimal in general, 
even if the presentation $\Phi$ is minimal.  
Its length, and the ranks of 
the free modules in $T_\bullet(\Phi, S)$ 
are completely determined by the rank of $L$ and by the 
matroid $\bfM(\Phi,S)$ on the set $S$ represented by the 
map $\phi$.    
When the matroid $\bfM(\Phi,S)$ is a uniform matroid, 
the complex $T_\bullet(\Phi, S)$ is canonically isomorphic to the 
Buchsbaum-Rim-Taylor complex from \cite{ChTc03}; in particular 
we recover the usual Taylor resolution of a monomial ideal. 
The second main result of this 
paper, Theorem~\ref{T:T-resolution-exactness}, is  
that the T-resolution of $(\Phi, S)$ is always a free multigraded 
resolution of the $R$-module $L$.  

It should be emphasized that different choices for the 
basis $S$ produce in general nonisomorphic 
matroids $\bfM(\Phi,S)$, 
and also nonisomoprhic T-resolutions. It is 
certainly an important open problem to analyze whether there 
exists an optimal choice of $S$, and what are the properties 
of a ``generic'' choice of $S$. 

This paper is written from the point of view of a commutative 
algebraist. Thus, while we employ many concepts and results 
from matroid theory, an attempt has been made to have a 
self-contained presentation of these facts.  In particular   
we review all relevant definitions about matroids, and  
we have included the proofs of all facts about matroids 
that appear after Section~1, even in the few cases 
when a reference was available in the literature.     
The material is organized as follows. 

In the first part of the paper, Section~1 through 
Section~5, we introduce all the objects necessary to 
define the T-complex $T_\bullet(\phi)$ and its 
augmentation $T_\bullet(\phi)^{+}$. The content is mostly 
of definitions, and should make for a relatively easy 
and smooth reading. More specifically, 
in Section~1 we review the 
basic notions from matroid theory that we will need, 
and in Section~2 we reintroduce (under a slightly 
different name) the somewhat neglected in the literature 
notion of a T-flat of a matroid 
and discuss its most basic properties; we 
also define the T-parts of a T-flat. In Section~3 
we consider a representation $\phi$ of 
a matroid $\bfM$ on a set $S$, and for each subset 
$A$ of $S$ we define its \emph{multiplicity space}. 
We also introduce the notions of T-space, multiplication 
maps, and  diagonal maps. In Section~4 we construct 
the maps $\phi_n^{IJ}$ that will be the building 
blocks for the differentials of the complexes 
$T_\bullet(\phi)$ and $T_\bullet(\phi)^{+}$.  
Finally, in Section~5 we give the definitions 
of the complex  
$T_\bullet(\phi)$ and of its augmentation 
$T_\bullet(\phi)^{+}$. 

Sections~6 and~7 constitute the second part of the 
paper. They are meant to serve as a reference 
to a reader who is not interested in the technical 
details of the proofs, but would like a good overview 
of the essential features of our constructions. 
In Section~6 we state all the main properties of the 
complexes  $T_\bullet(\phi)$ and $T_\bullet(\phi)^{+}$. 
In Section~7 we apply these results to construct from 
a finite free multigraded presentation 
$\Phi\: E \lra G$ of a Noetherian multigraded module 
$L$ and a chosen multihomogeneous basis $S$ of $E$ 
a finite free multigraded complex 
$T_\bullet(\Phi,S)$ and to prove 
in our second main result   
that this complex is a free resolution of $L$. 

In the next part of the article we work 
towards the proofs of all results stated in Section~6. 
In Section~8 we study the behavior of T-flats under 
restriction and contraction.  Section~9 is devoted 
to the structure of the T-parts of a T-flat, while 
Section~10 explores the relationship between 
connectedness properties of T-flats and the 
dimension of their multiplicity spaces.  
Section~11 provides key results on the 
behavior of multiplicity spaces under 
restriction and contraction. 

Starting with Section~12  we plunge into the  
most technical part of the paper and we prove 
that $T_\bullet(\phi)$ and $T_\bullet(\phi)^{+}$ 
are indeed complexes. The hard work continues 
in Sections~13 and~14, where we relate the 
T-complex of $\phi$ with the T-complex of a 
contraction of $\phi$. This provides all 
the facts that are needed in Section~15 to 
complete the proof of our first main result,  
namely that the T-complex of $\phi$ and 
its augmentation are resolutions of $\Ker(\phi)$ 
and of $\Coker(\phi)$, respectively.

I would like to thank Hara Charalambous, Bernd Sturmfels, 
Jerzy Weyman, and Thomas Zaslavski for very useful conversations 
on questions related to the material in this paper. Their comments 
contributed to my better understanding of the structures 
involved, and helped greatly improve the presentation in this 
article.  

We conclude this section with some remarks on notation. 
Throughout this paper $\Bbbk$ is a field, vector spaces 
are over $\Bbbk$, and unadorned 
tensor operations are over $\Bbbk$. 
For a finite dimensional 
vector space $X$ we denote by $\twedge X$ its 
top non-zero exterior power. We will often 
use that an exact sequence \ 
$0 \lra X \overset{\subset}\lra Y \lra Z \lra 0$ \ 
of finite dimensional vector spaces 
induces a canonical isomorphism 
$\twedge X\otimes \twedge Z\cong \twedge Y$ 
given by the formula \ 
$
x_1\wedge\dots\wedge x_p\otimes 
\bar z_1\wedge\dots\wedge\bar z_q 
\longmapsto 
x_1\wedge\dots\wedge x_p\wedge 
z_1\wedge\dots\wedge z_q  
$
\ where $p=\dim X$, $q=\dim Z$, and each  
$z_i$ is an element of $Y$ mapping onto $\bar z_i$. 
Finally, if 
$C_\bullet=(C_i,\psi_i)$ is a chain complex, its 
\emph{shift} 
$C[k]_\bullet$ is the complex $(C[k]_i,\psi[k]_i)$ 
with 
$C[k]_i=C_{i+k}$ and $\psi[k]_i=(-1)^k\psi_{i+k}$.

\section{Matroids}\mlabel{S:matroids}  

For the convenience of the reader and to establish our 
notation we recall here some notions and facts 
from matroid theory that will be needed. For further details 
and proofs the reader is refered to an introductory 
textbook on matroids such as \cite{We76} or \cite{Ox92}. 
We also note that for completeness (and in many cases due 
to lack of appropriate references)  
throughout the paper we provide a proof for every result 
on matroids that we use but have not mentioned in this 
section.

\subsection{Basics} 

A \emph{matroid} $\bfM$ on a finite set $S$ 
is given by a nonempty collection $\bfI$ of subsets of $S$ that 
satisfies the following two properties: 
\begin{enumerate} 
\item If $X\subset Y$ and $Y\in \bfI$ then $X\in \bfI$; 
\item If $X,Y \in \bfI$, and $|Y|=|X|+1$ then there is $y\in Y\smsm X$ 
such that $X\cup \{y\} \in \bfI$. 
\end{enumerate} 
The matroid $\bfM$ is called the \emph{empty} matroid if 
$S=\varnothing$ and $\bfI=\{\varnothing\}$.  
The sets in $\bfI$ are refered to as the \emph{independent} sets 
of the matroid $\bfM$. A \emph{base} of $\bfM$ is a maximal 
independent set. The matroid $\bfM$ is completely determined 
by specifying the collection $\bfB$ 
of its bases. 
A set is \emph{dependent} if it is not independent. 
A \emph{circuit} of $\bfM$ is a minimal dependent subset of $S$.   
Again, $\bfM$ is completely determined 
by specifying the collection of its circuits. 
The \emph{rank} in $\bfM$ of a subset $A$ of $S$
is the number $r^{\bfM}_A = \max\{|I|\mid I\subseteq A, I\in\bfI\}$. 
It is clear that 
\[
r_A^{\bfM} \le |A|,  
\]
with equality precisely when $A$ is independent. 
When there can be no confusion, we will omit the superscript $\bfM$. 
For any $A,B\subseteq S$ we have 
\[
r_A + r_B \ge r_{A\cup B} + r_{A\cap B}  
\]
and clearly $r_A\le r_B$ when $A\subseteq B$. 
We define the \emph{level} of $A\subseteq S$ to be the integer  
\[
\ell_A = \ell^{\bfM}_A = |A| - r^{\bfM}_A - 1.
\] 
In particular, a set $A$ is dependent if and only if $\ell_A\ge 0$. 
The \emph{closure} in $\bfM$ of a subset $A\subseteq S$ is the 
(unique) maximal subset $A^{\bsfC_{\bfM}}$ of $S$ that 
contains $A$ and has the same rank as $A$. When there can be no 
confusion we will simply write $A^{\bsfC}$. The set  
$A$ is \emph{closed} or a \emph{flat} in $\bfM$ if the closure 
of $A$ is $A$. 
The intersection of flats 
is a flat.  A maximal proper flat of $\bfM$ is called a 
\emph{hyperplane}. Every proper flat of $\bfM$ is an intersection of 
hyperplanes of $\bfM$. A flat $B$ is a \emph{cover} for a flat $A$ 
if $A\subsetneq B$ and there are no other flats between $A$ and $B$. 
It is clear this happens if and only if $r_B = r_A + 1$. Furthermore, if 
$B$ and $C$ are distinct covers for $A$ then their intersection is 
a flat containing $A$ hence equals $A$; thus the sets $B\smsm A$ and 
$C\smsm A$ are disjoint. Also, since every element $x\in S\smsm A$ 
belongs to a cover of $A$, the union of the covers of $A$ equals $S$. 
Therefore, if $A_1,\dots, A_k$ are all the distinct covers of $A$, 
we obtain a disjoint union 
\begin{equation}\elabel{E:T-flat-decomposition}
S \smsm A = (A_1\smsm A)\sqcup \dots \sqcup (A_k\smsm A). 
\end{equation}

\subsection{Minors, duality, connected components}

Let $Y$ be a subset of $S$. The \emph{restriction} of $\bfM$ to 
$Y$ is the matroid $\bfM|Y$ on the set $Y$ with collection of 
independent sets 
$\bfI|Y = \{I \mid I\subseteq Y, I\in \bfI\}$. 
It follows that for $A\subseteq Y$ one has 
\[
r_A^{\bfM|Y}=r_A^{\bfM}. 
\]
A subset $A\subseteq Y$ 
is a flat  of $\bfM|Y$ if and only if 
$A=A^{\bsfC_{\bfM}}\cap Y$. 
The \emph{contraction} of $\bfM$ to $Y$ is the matroid $\bfM.Y$ 
on $Y$ with collection of independent sets 
$\bfI.Y$ defined as follows. A subset $I\subseteq Y$ is
independent in $\bfM.Y$ if $I\cup J\in\bfI$ for every independent 
subset $J\subseteq S\smsm Y$. For a subset $A\subseteq Y$ we have 
\[
r_A^{\bfM.Y}= r^{\bfM}_{A\cup(S\smsm Y)} - r^{\bfM}_{S\smsm Y}.
\] 
A subset $A\subseteq Y$ is a flat of $\bfM.Y$ precisely when 
$A\cup (S\smsm Y)$ is a flat of $\bfM$. 
A \emph{minor} of $\bfM$ is a matroid obtained from $\bfM$ by a 
sequence of restrictions and contractions. 
  
The \emph{dual} matroid of $\bfM$ is the matroid $\bfM^*$ on $S$ 
with collection of bases given by 
$\bfB^*=\{ B \mid S\smsm B\in \bfB\}$. 
We have 
\begin{equation}\elabel{E:level-formula}
r^{\bfM^*}_A \quad = \quad |A| - r^{\bfM}_S + r^{\bfM}_{S\smsm A} 
             \quad = \quad r_S^{\bfM^*} - \ell_{S\smsm A}^{\bfM} - 1. 
\end{equation}
We also have $\bfM^*|T=(\bfM . T)^*$ and $\bfM^* . T=(\bfM|T)^*$. 
It is straightforward from the rank formula for $\bfM^*$ that 
a set $C$ is a circuit of $\bfM$ if and only if $S\smsm C$ is a 
hyperplane of $\bfM^*$. 

Let $\bfM_1$ be a matroid on a set $S_1$, and let $\bfM_2$ be a 
matroid on a set $S_2$. Then $\bfM_1 + \bfM_2$ is the matroid on 
the set $S_1\sqcup S_2$ where a subset $A\subseteq S_1\sqcup S_2$  
is independent in $\bfM_1 + \bfM_2$ if and only if the sets 
$A\cap S_1$ and $A\cap S_2$ are independent in $\bfM_1$ and 
$\bfM_2$, respectively. Thus if $\bfM$ is a matroid on $S$, and 
$\bfM=\bfM_1 + \bfM_2$ then $S=S_1\sqcup S_2$, and $\bfM_1=\bfM|S_1$, 
and $\bfM_2=\bfM|S_2$. It is immediate in this case that a 
subset $C\subseteq S$ is a circuit of $\bfM$ if and only if it 
is a circuit of either $\bfM_1$ or $\bfM_2$, and, 
furthermore, for a subset $A\subseteq S$ one has 
\[
r_A^{\bfM} = r_{A\cap S_1}^{\bfM_1} + r_{A\cap S_2}^{\bfM_2}. 
\]
A matroid $\bfM$ 
is called \emph{connected} if it is not a sum of two nonempty   
matroids. A subset $Y\subseteq S$ is called 
\emph{connected in $\bfM$} if the matroid $\bfM|Y$ is connected. 
In particular $S$ is connected in $\bfM$ precisely when $\bfM$ 
is connected. It is also straightforward from the definitions 
that a circuit is connected. 
The maximal connected subsets of $Y$ are called 
the \emph{connected components of $Y$ in $\bfM$}. The 
operation sum of matroids is associative and commutative.   
If $Y_1,\dots, Y_k$ are subsets of $Y$ 
such that $\bfM|Y=\bfM|Y_1 + \dots + \bfM|Y_k$,  then  
we write $Y=Y_1\oplus\dots\oplus Y_k$. It is clear 
in that case we also have for any subset $A\subseteq Y$ that 
\begin{equation}\elabel{E:rank-equality}
r_A = r_{A\cap Y_1} + \dots + r_{A\cap Y_k}. 
\end{equation} 
A convenient criterion which we will use often in the sequel 
without explicit reference is that $Y=Y_1\oplus\dots\oplus Y_k$ 
if and only if $Y=Y_1\sqcup\dots\sqcup Y_k$ and 
$r_Y=r_{Y_1}+\dots + r_{Y_k}$. See 
Proposition~\ref{T:direct-sum-criterion} for the proof.

\subsection{Representations}

Let $\bfM$ be a matroid on $S$. 
Let $W$ be a $\Bbbk$-vector space, and let $\phi\: S\lra W$ be 
a function of sets. We will abuse notation and write $\phi$ also 
for the canonically induced homomorphism of vector spaces 
\[
\phi\: U_S \lra W,
\] 
where 
$U=U_S$ is the $\Bbbk$-vector space with basis 
the set of symbols $\{e_a \mid a\in S\}$. 
We write $V=V(\phi)=\im \phi$ for the corresponding 
vector subspace of $W$. For a subset $I$ of 
$S$ we denote by $U_I$ the span in $U$ of the set  
$\{e_i\mid i\in I\}$,  
and by $V_I=V_I(\phi)$ the span in $W$ (hence in $V$) 
of the set $\{\phi(i)\mid i\in I\}$; thus $V_I$ is simply 
the image under $\phi$ of $U_I$.  
Given $a\in I$ we always denote by $e_a^*$ the element of 
the dual space $U_I^*$ such that $e_a^*(e_a)=1$ and 
$e_a^*(e_j)=0$ for $j\in I\smsm\{a\}$. 
We will often write $U_a$, $V_a$, and $S\smsm a$ 
instead of $U_{\{a\}}$, $V_{\{a\}}$, and $S\smsm\{a\}$, 
respectively. 

The map $\phi$ is called a \emph{representation} of the 
matroid $\bfM$ (over the field $\Bbbk$) if a subset 
$A\subseteq S$ is independent in $\bfM$  
precisely when $|A|=\dim_\Bbbk V_A$. 
If $\phi$ is a representation of $\bfM$ then clearly  
$r_I=\dim_{\Bbbk}V_I$. 

Let $\phi$ be a representation of $\bfM$. For a subset 
$Y\subseteq S$ let $\phi|Y$ be the restriction of $\phi$ 
to the subspace $U_Y$. Then 
\[
\phi|Y \: U_Y \lra W 
\] 
is a representation of the 
matroid $\bfM|Y$, and we obviously have for any subset 
$A$ of $Y$ that $V_A(\phi)=V_A(\phi|Y)$; in particular 
$V_Y(\phi) = V(\phi|Y)$. Note that if 
$Y = Y_1\sqcup \dots \sqcup Y_k$ then  
\[
\phi|Y=\phi|Y_1 + \dots +\phi|Y_k;  
\] 
if in addition $Y=Y_1\oplus\dots\oplus Y_k$ then 
the formula \altref{E:rank-equality} yields for any 
subset $A\subseteq Y$  
the direct sum decomposition of vector spaces    
\[
V_A = V_{A\cap Y_1}\oplus \dots \oplus V_{A\cap Y_k}. 
\] 
Also, let $\phi.Y=\pi_Y^{\phi}\circ(\phi|Y)$, where 
\[
\pi_Y^{\phi}\: W \lra \ol W = W/V_{S\smsm Y}(\phi)
\] 
is the canonical projection 
(when $\phi$ and/or $Y$ are clear from 
the context, we will write $\pi_Y$ or simply $\pi$ 
instead of $\pi_Y^{\phi}$). Then 
\[
\phi.Y \: U_Y \lra \ol W  
\] 
is a representation of the matroid $\bfM.Y$.

\subsection{Group actions} 

Let $\bfM$ be a matroid on a set $S$. 
We say that a group $H$ acting on the set $S$ also 
\emph{acts on the matroid} $\bfM$ if the collection 
of independent sets $\bfI$ is invariant under the 
action of $H$ on $S$. It is clear in this case 
that if $h\in H$ and $A$ is a circuit (respectively, 
base, flat, dependent set) of $\bfM$ then so is $h(A)$.  
Since a base of the dual matroid $\bfM^*$ is the 
complement of a base of $\bfM$, it follows that $H$ 
acts on $\bfM^*$ as well.  

Let $H$ be a group acting on $\bfM$, and let 
$\phi\: U_S \lra W$ be a representation of $\bfM$ 
over a field $\Bbbk$. If in addition 
$H$ also acts by linear transformations on the vector 
spaces $U_S$ and $W$ so that for every $h\in H$ and every 
$a\in S$ one has 
\[
h(e_a)=\chi(h,a)e_{h(a)}
\] 
for some (necessarily nonzero) scalar $\chi(h,a)\in\Bbbk$, 
and so that the map $\phi$ is $H$-equivariant, then   
we say that $\phi$ is an \emph{$H$-equivariant  
representation of the matroid $\bfM$}.  
In that case for every $h\in H$ and every 
$A\subseteq S$ we have   
canonical isomorphisms $h\: U_A \lra U_{h(A)}$ and 
$h\: V_A \lra V_{h(A)}$, as well as  
$h_*=(h^{-1})^*\: V_A^* \lra V_D^*$.

\begin{example}\mlabel{E:uniform-matroid}
Let $\bfM$ be the matroid on $S$ with $\bfI$ consisting 
of all subsets of $S$ of size less than or equal to $r$. 
$\bfM$ is called the \emph{uniform matroid of rank $r$ on $S$}. 
It is clear that a subset $A\subseteq S$ 
is a base of $\bfM$ precisely when $|A|=r$, 
and is a circuit precisely 
when $|A|=r+1$. Also, the formulas 
\[
r_A    = \min\{r, |A|\} \quad\text{and}\quad  
\ell_A = \max\{-1,|A|-r-1\}.
\]
are straightforward.   
\end{example}

\begin{example}\mlabel{E:matrix-matroid} 
Let $S=\{1,2,3,4\}$, let $W=\mathbb Q^2$, and let 
$\g_1$ and $\g_2$ be the standard basis vectors of $W$. 
Define $\phi\: S \lra W$ by letting $\phi(i)$ be the 
$i$th column of 
\[
M= 
\begin{pmatrix} 
1 & 1 & 1 & 1 \\
1 & 1 & 2 & 3 
\end{pmatrix}. 
\]
Thus the corresponding map $\phi\: U_S \lra W$ is  
given by 
$M$. Let $\bfM$ be the matroid represented by $\phi$.  
Then the circuits of $\bfM$ are 
$\{1,2\}$, $\{1,3,4\}$, and $\{2,3,4\}$. 
\end{example}

\section{T-flats}

For the rest 
of the paper $\bfM$ denotes a matroid on a finite set 
$S$, and all sets considered are subsets of $S$. 

The main goal of this section is to recall (under a slightly 
different name) the somewhat neglected in the literature 
notion of a 
\emph{T-flat} of $\bfM$. Its importance for us lies 
in the fact that the collection of T-flats provides 
the most convenient canonical 
indexing set for the components of the complex  
$T_\bullet(\Phi,S)$ from Section~7, which generalizes 
to the setting of multigraded modules the 
Taylor resolution of a monomial ideal. 
The study of the collection of T-flats of a matroid goes back 
to the papers of Tutte~\cite{Tu58} and~\cite{Tu71} who refers to  
T-flats of level $n\ge 0$ as ``flats of $\bfM$'' or as  
``$n$-cells of $\bfM$''. From our point of view the flats 
terminology is more appropriate. Unfortunately, in 
modern language the term ``flat of $\bfM$'' has taken on 
a different and well established meaning. Thus we decided 
on ``T-flat'' as a smallest deviation from ``flat''    
with the ``T'' intended 
as a reference to both Tutte, and Taylor.

\begin{definition}\mlabel{D:T-flat} 
(a) A set $A$ is called a \emph{T-flat} of $\bfM$ 
if and only if $S\smsm A$ is a proper flat of the dual 
matroid $\bfM^*$. 

(b) When $n\ge 0$ we write 
$\mathcal T_n = \mathcal T_n(\bfM)$ for 
the collection of T-flats of level $n$ in $\bfM$. 
We also write $\mc T_{-1}=\mc T_{-1}(\bfM)$ for the 
collection of $1$-element subsets of $S$.   
\end{definition} 

The following observations are essentially due to 
Tutte~\cite{Tu58}. 

\begin{remarks}\mlabel{T:T-remarks}  
{\rm(\cite{Tu58})}. 
Let $A$ be a T-flat in $\bfM$ of level $n$. 

(a) Then $S\smsm A$ is a proper flat in 
$\bfM^*$ of rank $r^{\bfM^*}_{S\smsm A}= r^{\bfM^*}_S - n - 1$. 
In particular $n\ge 0$, and when $n=0$ the flat $S\smsm A$ is a 
hyperplane in $\bfM^*$ hence $A$ is a circuit of $\bfM$. 
Thus the T-flats in $\bfM$ of smallest level all have level $0$ 
and are nothing but the circuits of $\bfM$.  

(b) Let $n\ge 1$, and let $B\subsetneq A$ 
be a T-flat. Then $S\smsm A \subsetneq S\smsm B$, hence 
a straightforward application of \altref{E:level-formula} yields 
$\ell^{\bfM}_B\le n-1$ with equality precisely when $S\smsm B$ is a 
cover for $S\smsm A$ in $\bfM^*$, that is, when $B$ is a maximal 
T-flat properly contained in $A$.   
Therefore the T-flats of $\bfM$ that are maximal with respect 
to the property of being properly contained in $A$ 
are precisely the T-flats of level $n-1$ contained in $A$. 

(c) It is immediate from parts (a) and (b) that a T-flat 
is minimal if and only if it is a circuit, and that every 
chain of maximal length of T-flats contained inside $A$ has 
the form 
\[
I^{(0)} \subsetneq \dots \subsetneq I^{(n)}=A 
\]
where for each $i$ the set $I^{(i)}$ is a T-flat of level $i$.  

(d) Let $n\ge 1$ and let 
$A_1,\dots,A_k$ be all the T-flats of level $n-1$ 
contained in $A$. Applying part (b) and 
\altref{E:T-flat-decomposition} yields  
a natural partition of $A$ as a disjoint union 
\begin{equation}\elabel{E:T-partition}
A = I_1 \sqcup \dots \sqcup I_k, 
\end{equation}
where $I_i=A\smsm A_i$ for each $i$. 

(e) Since $S\smsm A$ is a proper flat in $\bfM^*$ 
precisely when it is an intersection of hyperplanes of $\bfM^*$, 
it is immediate from part (a) that $A$ is a T-flat precisely when 
it is a union of circuits of $\bfM$. In particular, unions of 
T-flats are T-flats. 

(f) Let $A=\{a\}$ be a singleton, i.e.~$A\in\mc T_{-1}(\bfM)$. 
Note that $\ell_A\le 0$ with equality if and only if $A$ is 
dependent, in which case $A$ is a circuit. 
(An element $a\in S$ such that $\{a\}$ is a dependent is 
called a \emph{loop} of $\bfM$.) Thus $A$ is an element 
of both $\mc T_0(\bfM)$ and $\mc T_{-1}(\bfM)$ if and only 
if $A=\{a\}$ for some loop $a\in S$. It is clear that this 
is the only case when two sets $\mc T_i(\bfM)$ and 
$\mc T_j(\bfM)$ with $-1\le i< j$ have a nontrivial 
intersection.  
\end{remarks}

\begin{definition}\mlabel{D:T-parts} 
Let $A$ be a T-flat in $\bfM$ of level $n$, 
and let $A_1,\dots, A_k$ be all elements of 
$\mc T_{n-1}(\bfM)$ contained in $A$. 

(a) We refer to the sets $I_i=A\smsm A_i$ as the 
\emph{T-parts} of the T-flat $A$. 
  
(b) When $n\ge 1$ the natural 
partition \altref{E:T-partition} is called the 
\emph{T-partition} of $A$.  
\end{definition}

\begin{remark}\mlabel{T:H-action-on-T-flats-and-parts} 
Let $H$ be a group acting on $\bfM$, and let $h\in H$. It is 
clear that if $I$ is an element of $\mc T_n(\bfM)$ (respectively, 
a T-part of a T-flat $A$), then $h(I)$ is an element 
of $\mc T_n(\bfM)$ (repsectively, a T-part of the T-flat 
$h(A)$).  In particular, $h$ sends the T-partition of a T-flat $A$ 
into the T-partition of the T-flat $h(A)$.  
\end{remark}

The following proposition provides an alternative description 
of the T-flats of the matroid $\bfM$.  

\begin{proposition}\mlabel{T:inductive-T-flats} 
Let $A$ be a set of level $n\ge 1$. 

Then $A$ is a T-flat if and only  
if it has a decomposition $A=I_1\sqcup\dots\sqcup I_k$ such that 
$A_j=A\smsm I_j$ is a T-flat of level $n-1$ for each $j$. In that 
case $A_1,\dots,A_k$ are all the T-flats of \ $\bfM$ of 
level $n-1$ that are contained in $A$.    
\end{proposition}

\begin{proof} 
The ``only if'' part of the proposition is immediate from 
Remark~\ref{T:T-remarks}(b). To show the ``if'' part, assume 
$A=I_1\sqcup\dots\sqcup I_k$ with $A_j=A\smsm I_j$ a T-flat of level 
$n-1$ for each $j$. Since $n\ge 1$, 
we must have $k\ge 2$. Then $A=A_1\cup A_2$, and since each $A_j$ 
is a T-flat the desired conclusion follows 
by Remark~\ref{T:T-remarks}(c).   
\end{proof}

\begin{example}\mlabel{E:uniform-matroid-0}
Let $\bfM$ be the uniform matroid of rank $r$ on the set 
$S$ as in Example~\ref{E:uniform-matroid}. 
It is clear that  a subset $A\subseteq S$ is a T-flat 
of level $n\ge 0$ in $\bfM$ if and only if $|A|= n+r+1$.  
Therefore when 
$n\ge 1$ and $A$ is a T-flat of level $n$  the T-parts 
of $A$ are all the $1$-element subsets of~$A$. 
\end{example}

\begin{example}\mlabel{E:matrix-matroid-0} 
Let $\bfM$ be the matroid from Example~\ref{E:matrix-matroid}. 
The T-flats of $\bfM$ are: 
\[
\begin{aligned} 
\text{level $0$} &: \quad \{1,2\}, \{1,3,4\}, \{2,3,4\} \\ 
\text{level $1$} &: \quad \{1,2,3,4\}. 
\end{aligned} 
\]
In particular, $\{1,2,3,4\}=\{3,4\}\sqcup\{2\}\sqcup\{1\}$ is 
the T-partition of the only T-flat of positive level in $\bfM$. 
\end{example}

\section{Multiplicity spaces, T-spaces, diagonal maps}

In this section 
$\phi\: U_S \lra W$ is a representation  
over a field $\Bbbk$ of a matroid $\bfM$ on a finite set $S$. 

Our main goal 
is to introduce in Definition~\ref{D:multiplicity-space}
the \emph{multiplicity space} 
$S_A(\phi)$ of a set $A$. This is a new object, which 
appears to encode subtle combinatorial properties of 
the matroid $\bfM|A$. A detailed investigation
of the properties of these multiplicity spaces 
will be carried out in later sections. For us their  
significance lies in the fact that the multiplicity 
space of a set $A$ of level $n\ge 0$ measures the 
contribution of $A$ to the component of homological 
degree $n$ in the complex $T_\bullet(\phi)$.  
We also introduce in Definition~\ref{D:diagonal-map}
the \emph{diagonal maps} between mutiplicity spaces. 
These maps are the essential ingredients out of 
which the differentials of the complex $T_\bullet(\phi)$ 
are built.  

We begin by introducing the following objects associated 
with a T-flat $I$.    

\begin{definition}\mlabel{D:multiplicity-space-indexing} 
(a) Let $I$ be a T-flat of level $n$. We write  
\[
C_{\bfM}(I)=
\{ I^{(0)}\subsetneq\dots\subsetneq I^{(n)} = I \mid 
        \text{ where } I^{(i)}\in\mathcal T_i 
        \text{ for each } i\ge 0 \} 
\]
for the collection of all maximal chains of T-flats 
contained in the T-flat $I$. 
When there can be no confusion we will write 
simply $C(I)$. 

(b) Let $\mathbb I$ be a chain 
$I^{(0)}\subsetneq\dots\subsetneq I^{(n)}$ 
in $C(I)$. When $n=0$ then $\mathbb I$ is the only 
element of $C(I)$, and we set 
$V(\mathbb I)=S_0W=\Bbbk$. When $n\ge 1$ we define the 
vector subspace 
$V(\mathbb I) \subseteq S_nW$ of the $n$th symmetric 
power of $W$ as 
\[
V(\mathbb I)\quad = \quad 
\bigl(
V_{I^{(1)}\smsm I^{(0)}}\cap V_{I^{(0)}}
\bigr)
\ \cdot\ \ldots\ \cdot \  
\bigl(
V_{I^{(n)}\smsm I^{(n-1)}}\cap V_{I^{(n-1)}}
\bigr),   
\]
where the product of these vector spaces is taken 
in the symmetric algebra of $W$. 
\end{definition}

\begin{remark}\mlabel{T:canonical-generators} 
Since the levels of $I^{(i)}$ and $I^{(i-1)}$ 
differ by exactly one, a straightforward computation 
shows that each space 
$V_{I^{(i)}\smsm I^{(i-1)}}\cap V_{I^{(i-1)}}$ 
is either $0$ or is $1$-dimensional. Therefore the 
space $V(\mathbb I)$ is either $0$ or is 
$1$-dimensional. This observation is refined 
further in Remark~\ref{T:rank-1-intersection-a}
\end{remark}

We are now ready to give the definition of 
multiplicity space. 

\begin{definition}\mlabel{D:multiplicity-space}
Let $I\subseteq S$ be a set of level $n$.
 
(a) If $I$ is not a T-flat then we set $S_I(\phi)=0$. 
If $I$ is a T-flat, then we set 
\[
S_I(\phi) = \sum_{\mathbb I\in C(I)} V(\mathbb I).  
\]
We call the space $S_I(\phi)$ the 
\emph{multiplicity space} 
of the dependent set $I$. 
When there can be no confusion, 
we simply write $S_I$.

(b) When $n\ge 0$ we set  
\[
T_I(\phi)= S_I(\phi)^*\otimes \twedge U_I
                      \otimes \twedge V_I(\phi)^*
\]
and call the space $T_I(\phi)$ the \emph{T-space} 
of $I$. When there can be no confusion, we simply 
write $T_I$. 

(c) Note that when $A=\{a\}$ is an element of 
$\mc T_0(\bfM)$ then we have by 
Remark~\ref{T:basic-multiplicity-properties}(a) 
that $S_A^*=\Bbbk^*=\Bbbk$, also $\twedge U_A = U_a$, 
and $V_A=0$ hence $\twedge V_A^*=\Bbbk$; 
thus for the T-space $T_A$ we have 
$T_A = \Bbbk\otimes U_a\otimes\Bbbk = U_a$. 
We extend this identification to all $1$-element 
subsets of $S$ by setting 
\[
T_A(\phi) = T_A = U_a 
\] 
for every element $A=\{a\}\in \mc T_{-1}(\bfM)$.    
\end{definition}

\begin{remarks}\mlabel{T:basic-multiplicity-properties} 
Let $I$ be a T-flat of level $n$.  

(a) When $n=0$ there is exactly one chain 
$\mathbb I$ in $C(I)$, thus $S_I=V(\mathbb I)=S_0W=\Bbbk$. 

(b) We should note that the T-space of $I$ is essentially the dual 
of $S_I$ tensored with a copy of the field $\Bbbk$; in particular 
the dimensions of $T_I$ and $S_I$ are the same.  The special  
form of $T_I$ is what turns out to be necessary  
in order to define in a canonical way the differentials in our 
complex $T_\bullet(\phi)$. 

(c) While the multiplicity space $S_I(\phi)$ clearly depends on the 
representation $\phi$, we will show that its dimension is an
invariant of the matroid $\bfM|I$, and is independent of the 
representation map $\phi$. 

(d) The collection of spaces $V(\mathbb I)$ 
determine in a natural way a representation of a certain 
matroid on the set of chains $C_{\bfM}(I)$. It is a very 
interesting open question whether this matroid is in fact 
independent of the representation map $\phi$. 
\end{remarks}

\begin{remarks}\mlabel{T:H-action-on-multiplicity-spaces} 
Let $\phi$ be an $H$-equivariant representation of $\bfM$, 
and let $h\in H$.  

(a) If $\mathbb I$ is a chain 
$I^{(0)}\subsetneq\dots\subsetneq I^{(n)}$ in $C(I)$ then 
we obtain that $h(\mathbb I)$ is a chain 
\[
h(I^{(0)})\subsetneq\dots\subsetneq h(I^{(n)})
\] 
in the set $C\bigl(h(I)\bigr)$. In particular, $h$ induces 
a bijection of the sets $C(I)$ and $C\bigl(h(I)\bigr)$. 

(b) The action of $h$ on $W$ induces canonically an 
automorphism of the symmmetric algebra of $W$. 
Furtermore, if $\mathbb I$ is a chain 
$I^{(0)}\subsetneq\dots\subsetneq I^{(n)}$ in $C(I)$, 
then $h$ induces for each $i\ge 1$ an isomorphism 
\[ 
h\: \ 
V_{I^{(i-1)}}\cap V_{I^{(i)}\smsm I^{(i-1)}} 
\quad \lra \quad  
V_{h(I^{(i-1)})}\cap V_{h(I^{(i)})\smsm h(I^{(i-1)})}. 
\]
Thus $h$ induces an isomorphism of the spaces $V(\mathbb I)$ 
and $V\bigl(h(\mathbb I)\bigr)$, hence also an isomorphism 
\[
h\: S_I(\phi) \lra S_{h(I)}(\phi).   
\]
of multiplicity spaces. 

(c) Let $I$ be a T-flat of level $n$. We have an 
isomorphism $h\: T_I(\phi) \lra T_{h(I)}(\phi)$ given by  
\[
\begin{CD} 
T_I = 
S_I^*\otimes\twedge U_I\otimes\twedge V_I^* 
@> h_*\otimes \wedge h\otimes \wedge h_* >> 
S_{h(I)}^*\otimes\twedge 
U_{h(I)}\otimes\twedge V_{h(I)}^* = 
T_{h(I)} 
\end{CD}
\] 
where we write $h_*$ for $(h^{-1})^*$. 
\end{remarks}

We conclude this section with the definition of the 
\emph{multiplication maps} and the \emph{diagonal maps}. 
These maps are the essential ingredient of the 
differentials in the complex $T_\bullet(\phi)$.

\begin{definition}\mlabel{D:diagonal-map}
Let $I$ be a T-flat of level $n\ge 1$ and let 
$J\subsetneq I$ be a T-flat of level $n-1$. 
Multiplication in the symmetric algebra 
of $W$ induces an injective  map 
\[
\nu\: (V_J\cap V_{I\smsm J})\otimes S_J \lra S_I.  
\]
which we call the \emph{multiplication map}. 
Summing over all T-flats of level $n-1$ in $I$, 
we obtain a surjective map with the same name 
\begin{equation}\elabel{E:multiplication-map}
\bigoplus_{\begin{smallmatrix} 
           J\in\mathcal T_{n-1} \\ J\subset I 
           \end{smallmatrix} } 
      (V_J\cap V_{I\smsm J})\otimes S_J  
\quad\overset{\nu}\lra\quad  S_I 
\end{equation}  
Taking duals we obtain a surjective map 
\[
\Delta\: 
S_I^* \lra  (V_J\cap V_{I\smsm J})^*\otimes S_J^*    
\]
which we refer to as the \emph{diagonal map}, and 
an injective map with the same name  
\begin{equation}\elabel{E:diagonal-map} 
S_I^* \qquad\overset{\Delta}\lra\quad   
\bigoplus_{\begin{smallmatrix} 
           J\in\mathcal T_{n-1} \\ J\subset I 
           \end{smallmatrix} }  
         (V_J\cap V_{I\smsm J})^*\otimes S_J^*.  
\end{equation}
\end{definition}

\begin{remark}\mlabel{T:H-equivariant-maps} 
Let $\phi$ be $H$-equivariant representation of 
$\bfM$, and let $h\in H$. 
Let $I$ be a T-flat of level $n\ge 1$, and let $J$ be 
a T-flat of level $n-1$ contained inside $I$. 
It is straightforward to verify that the map 
$h\: S_I \lra S_{h(I)}$ from 
Remark~\ref{T:H-action-on-multiplicity-spaces}(c) 
commutes with the multiplication maps $\nu$, thus 
setting \  
$h_*= (h^{-1})^*$ \ yields a commutative diagram     
\[
\begin{CD} 
S_I^* 
@> \Delta >> 
(V_J\cap V_{I\smsm J})^*\otimes S_J^*
\\ 
@V h_* VV 
@VV h_*\otimes h_* V
\\ 
S_{h(I)}^* 
@> \Delta >> 
\bigl(
V_{h(J)}\cap V_{h(I)\smsm h(J)}
\bigr)^*\otimes S_{h(J)}^*
\end{CD}
\]
with vertical maps that are isomorphism. 
\end{remark}

\begin{example}\mlabel{E:uniform-matroid-1} 
Let $\phi\: U_S\lra W$ be a representation of the uniform 
matroid $\bfM$ of rank $r$ on $S$, and let $V=\im\phi$. 
Then for each T-flat $I$ of level $n$ in $\bfM$ we clearly 
have $V_I=V$. Furthermore,  
\[
S_I(\phi) = S_nV.  
\]
Indeed, this is trivially true when $n=0$. When $n\ge 1$ 
by Example~\ref{E:uniform-matroid-0}  
the T-flats of level $n-1$ in $I$ are all sets of the 
form  $J=I\smsm\{a\}$ for some $a\in I$, therefore   
we have 
$
\sum_{
I\supset J\in\mathcal T_{n-1}}  
V_J\cap V_{I\smsm J} = 
\sum_{a\in I} V_{I\smsm a}\cap V_a = 
\sum_{a\in I} V_a = V_I = V.
$ 
The desired equality now follows by induction on $n$ 
in view of the surjectivity of 
\altref{E:multiplication-map}.  
\end{example}

\begin{example}\mlabel{E:matrix-matroid-1}
Let $\phi\: U_S \lra W$ be the representation of the 
matroid $\bfM$ from Examples~\ref{E:matrix-matroid} 
and~\ref{E:matrix-matroid-0}. Let $\{\g_1^*,\g_2^*\}$ be 
the basis of $W^*$ dual to the standard basis of~$W$.  
Then we have 
\[
S_{\{1,2,3,4\}} = 
V_{\{1,2\}}\cap V_{\{3,4\}} + 
V_{\{1,3,4\}}\cap V_{\{2\}} + 
V_{\{2,3,4\}}\cap V_{\{1\}} = V_{\{1,2\}},  
\]
therefore for the T-spaces of the T-flats of $\bfM$ 
we obtain 
\[ 
\begin{aligned} 
T_{\{1,2\}}
&=
\mathbb Q\otimes \wedge^2 U_{\{1,2\}}
         \otimes V_{\{1,2\}}^*  
\\ 
&= 
\mathbb Q\cdot 
\big\langle 
1\otimes e_{12}\otimes (\g_1+\g_2)^*
\big\rangle  
\\[+8pt] 
T_{\{1,3,4\}}
&= 
\mathbb Q\otimes\wedge^3 U_{\{1,3,4\}}
         \otimes \wedge^2 W^* 
\\ 
&=
\mathbb Q\cdot 
\big\langle 
1\otimes e_{134}\otimes \g_{12}^* 
\big\rangle 
\\[+8pt] 
T_{\{2,3,4\}}
&= 
\mathbb Q\otimes\wedge^3 U_{\{2,3,4\}}
         \otimes \wedge^2 W^* 
\\ 
&= 
\mathbb Q\cdot 
\big\langle 
1\otimes e_{234} \otimes \g_{12}^* 
\big\rangle     
\\[+8pt]  
T_{\{1,2,3,4\}}
&= 
S^*_{\{1,2,3,4\}}\otimes \wedge^4 U_{\{1,2,3,4\}}
                 \otimes \wedge^2 W^* 
\\ 
&= 
\mathbb Q\cdot 
\big\langle 
(\g_1+\g_2)^* \otimes e_{1234}\otimes \g_{12}^* 
\big\rangle   
\end{aligned} 
\] 
where we use the shorthand $z_{a_1\dots a_k}$ for 
$z_{a_1}\wedge\dots\wedge z_{a_k}$. 
The remaining nonzero T-spaces of $\bfM$ are, 
for trivial reasons, \ 
$T_{\{i\}}=U_{\{i\}}=\mathbb Q \cdot e_i$ \  
for~$i=1,2,3,4$.  
\end{example}

\section{The homomorphisms $\phi_n^{IJ}$}

As in the previous section, we fix a representation 
$\phi\: U_S \lra W$ over $\Bbbk$ of a matroid $\bfM$ 
on a finite set $S$.  
Our goal is to 
introduce in Definition~\ref{D:differentials-of-T}
for every T-flat $I$ of level 
$n$ and for every element $J\in \mc T_{n-1}(\bfM)$ 
contained inside $I$ a homomorphism between T-spaces 
\[
\phi_n^{IJ}\: T_I(\phi)\lra T_J(\phi).
\] 
The maps 
$\phi_n^{IJ}$ will be used in the next section as the 
building blocks of the differentials of the complex 
$T_\bullet(\phi)$. First, we establish the following 
notation. 

\begin{definition}\mlabel{D:notation}  
Let $J\subseteq I$ be subsets of $S$. 

(a) We define the space  $K_{IJ}$ by the formula 
\[
K_{IJ}=V_I/V_J.
\] 

(b) Since $V_I=V_J + V_{I\smsm J}$ we have a   
natural commutative diagram 
\[
\begin{CD}
0 
@>>> 
V_J\cap V_{I\smsm J}     
@> \subseteq >> 
V_{I\smsm J} 
@>>> 
K_{IJ} 
@>>> 
0 
\\ 
@.    
@V \subseteq VV           
@V \subseteq VV   
@V = VV  
@.  
\\ 
0 
@>>> 
V_J                  
@> \subseteq >> 
V_I          
@>>> 
K_{IJ} 
@>>> 
0
\end{CD}
\]
with exact rows; we use it to identify canonically $K_{IJ}$ 
and $V_{I\smsm J}/(V_J\cap V_{I\smsm J})$. 

(c) As mentioned in the Introduction, we have a naturally  
induced diagonal isomorphism 
\[
\bfd=\bfd_{\bold{IJ}}\: \ \ \ \twedge U_I \ \ \lra \ \ 
            \twedge U_{I\smsm J}\otimes 
            \twedge U_J  
\] 
and a canonical isomorphism 
\[
\bold a = \bold{a^M_{IJ}}\: \ \ \  
\twedge V_I^* \ \ \lra \ \ 
\twedge K^*_{IJ}\otimes \twedge V^*_J   
\]
induced by the bottom row of 
the commutative diagram from Part~(b) above.  

(d) Let $I$ be a T-flat of level $n\ge 1$, and let 
$J$ be a T-flat of level $n-1$ inside $I$. 
Recall from Remark~\ref{T:canonical-generators} that 
the space $V_J\cap V_{I\smsm J}$ is either $0$ or 
has dimension~$1$. In particular, the top row of the 
commutative diagram from Part~(b) induces a canonical 
homomorphism  
\[
\bold b = \bold{b^M_{IJ}}\: \ \ \ 
(V_J\cap V_{I\smsm J})^* \otimes \twedge K^*_{IJ} 
\ \ \lra \ \   
\twedge V^*_{I\smsm J}    
\]
which is zero when $V_J\cap V_{I\smsm J}=0$, and is an 
isomorphism otherwise. 
\end{definition} 

We are now ready to state the desired definition. 

\begin{definition}\mlabel{D:differentials-of-T}
Let $I$ be a T-flat of level $n\ge 0$, 
and let $J$ be an element of $\mc T_{n-1}(\bfM)$ 
contained inside $I$. 
We define a homomorphism  
\[
\phi_n^{IJ}\: T_I(\phi) \lra  T_J(\phi) 
\]
as follows. 

(a) When $n=0$ the T-flat $I$ is a circuit   
and $J=\{a\}$, hence $I\smsm a$ is independent with 
$V_I=V_{I\smsm a}$. Furthermore        
$S_I(\phi)^*=\Bbbk^*=\Bbbk$ and 
$T_J(\phi)= U_a$, and we define the 
homomorphism \ 
$\phi_n^{IJ}= \ \phi_0^{I,a}\: \ T_I \lra U_a$ \ 
as the composition  
\[
\begin{CD} 
S_I^*\otimes \twedge U_I \otimes \twedge V_I^* 
\\ 
@VV 1\otimes\bfd\otimes 1 V                          
\\ 
\Bbbk\otimes (U_a\otimes \twedge U_{I\smsm a} ) 
                 \otimes \twedge V_{I\smsm a}^*            
\\ 
@VV \tau V                                                
\\ 
\twedge U_{I\smsm a} \otimes\twedge V_{I\smsm a}^* 
                     \otimes U_a                          
\\ 
@VV \wedge\phi  \otimes 1 \otimes 1 V                  
\\
( \twedge V_{I\smsm a}\otimes\twedge V_{I\smsm a}^* ) 
                      \otimes U_a                         
\\ 
@VV \mu \otimes 1 V                          
\\ 
U_a      
\end{CD}
\]
where $\mu \: X\otimes X^* \lra \Bbbk$ is the canonical 
evaluation map, and $\tau$ is the isomorphism 
which simply permutes the components of the tensor product 
as indicated. 

(b) When $n\ge 1$ then $J$ is a T-flat of level $n-1$ and 
we define the map $\phi_n^{IJ}$ as the composition 
\[
\begin{CD} 
T_I(\phi) 
@>\Delta\otimes\bfd\otimes\bold a >> 
\mc Q_{IJ}(\phi) 
@>\tau >> 
\mc R_{IJ}(\phi) 
@>\bold b\otimes\wedge\phi\otimes 1 >> 
\mc S_{IJ}(\phi) 
@>\mu\otimes 1  >>  
T_J(\phi)   
\end{CD}
\]
where 
\[
\begin{aligned}
\mc Q_{IJ} &= (V_J\cap V_{I\smsm J})^* 
     \otimes S_J^*  
     \otimes \bigl( \twedge U_{I\smsm J}   
     \otimes        \twedge U_J \bigr)  
     \otimes \bigl( \twedge K_{IJ}^*           
     \otimes        \twedge V_J^* \bigr),  
\\ 
\mc R_{IJ} &= \bigl[ (V_J\cap V_{I\smsm J})^* 
    \otimes \twedge K_{IJ}^* \bigr]           
    \otimes \twedge U_{I\smsm J}   
    \otimes \bigl( S_J^*\otimes \twedge U_J 
    \otimes \twedge V_J^* \bigr),           
\\ 
\mc S_{IJ} &= \twedge V_{I\smsm J}^*   
    \otimes   \twedge V_{I\smsm J}          
    \otimes T_J,     
\end{aligned}
\]
the map 
$\Delta\: S_I^* \lra (V_J\cap V_{I\smsm J})^*\otimes S_J^*$ 
is the diagonal map from Definition~\ref{D:diagonal-map}, 
and the map 
$\wedge\phi\: \twedge U_{I\smsm J} \lra \twedge V_{I\smsm J}$ 
is zero if $I\smsm J$ is not independent, and is the 
canonical isomorphism induced by $\phi$ otherwise. 
\end{definition}

\begin{remarks}\mlabel{T:kernel-of-diagonal-map} 
(a) If $I=\{a\}$ is a T-flat of level $n=0$, then the 
only element of $\mc T_{n-1}=\mc T_{-1}$ contained inside $I$ 
is $J=\{a\}=I$. In that case $T_I = U_a = T_J$, and the 
formula for \ $\phi_n^{IJ}=\phi_0^{a,a}\: U_a \lra U_a$ \ 
yields simply the identity map of $U_a$.   

(b) We note that when $n\ge 1$ the map $\phi_n^{IJ}$ is 
essentially the same as the diagonal map $\Delta$. In fact,   
all other maps appearing in the definition are isomorphisms 
except possibly the map $\wedge\phi$ which is zero if 
$I\smsm J$ is not independent. In that case however we have  
$V_J\cap V_{I\smsm J}=0$ by 
Remark~\ref{T:rank-1-intersection-a}, 
hence also both $\Delta=0$ and $\phi_n^{IJ}=0$. Therefore 
we always have 
$
\Ker(\phi_n^{IJ})=\Ker(\Delta)\otimes\twedge U_I 
                              \otimes\twedge V_I^*.   
$
\end{remarks}

\begin{remarks}\mlabel{T:H-equivariant-phi-sub-n-IJ} 
Let $\phi$ be an $H$-equivariant representation of the 
matroid $\bfM$, and consider an element $h\in H$. 

(a) It is straightforward to verify from the definitions 
that for each T-flat $I$ of level $n\ge 1$ and each T-flat 
$J$ of level $n-1$ contained inside $I$ we have a 
commutative diagram 
\[
\begin{CD} 
T_I
@> \phi_n^{IJ} >> 
T_J 
\\ 
@V h VV 
@VV h V 
\\ 
T_{h(I)} 
@> \phi_n^{h(I),h(J)} >> 
T_{h(J)} 
\end{CD}
\]
where the vertical maps are the isomorphisms from 
Remark~\ref{T:H-action-on-multiplicity-spaces}(c). 

(b) Similarly, when $I$ is a circuit and $a\in I$ 
one verifies from the definitions that  
\[
\begin{CD} 
T_I 
@> \phi_0^{I,a} >> 
U_a  
\\ 
@V h VV 
@VV h V 
\\ 
T_{h(I)} 
@> \phi_0^{h(I),h(a)} >> 
U_{h(a)}
\end{CD}
\]
is a commutative diagram. 
\end{remarks}

\begin{example}\mlabel{E:matrix-matroid-1-point-5} 
Let $\bfM$ be the matroid represented by the map $\phi$ 
from Example~\ref{E:matrix-matroid}. 
In Example~\ref{E:matrix-matroid-1} we described all the 
nonzero T-spaces and gave a basis for each one of them 
over the field $\mathbb Q$. Here we compute the images 
of these basis elements under each of the maps $\phi_n^{IJ}$. 
Consider, say, the map 
\[
\phi_1^{\{1,2,3,4\},\{1,2\}}\: 
T_{\{1,2,3,4\}}\lra T_{\{1,2\}}. 
\] 
According to Definition~\ref{D:differentials-of-T} 
the image of the basis element of $T_{\{1,2,3,4\}}$ 
is obtained  through the following sequence of 
transformations: 
\[
\begin{aligned} 
(\g_1 +\g_2)^* \ \otimes \ & \ e_{1234}  
               \ \otimes \ \g_{12}^*   
\\ 
\overset{\Delta\otimes\bfd\otimes{\bf a}}\longmapsto 
& \ (\g_1+\g_2)^*\otimes 1
             \otimes e_{34}\otimes e_{12}  
             \otimes (\g_1^*-\g_2^*)\otimes 
                      \g_2^*\big|_{V_{\{1,2\}}}  
\\ 
\overset{\tau}\longmapsto \ \   
& \ (\g_1+\g_2)^*\otimes (\g_1^*-\g_2^*) 
             \otimes e_{34}
             \otimes 1\otimes e_{12}\otimes(\g_1+\g_2)^*
\\ 
\overset{{\bf b}\otimes\wedge\phi\otimes 1}\longmapsto 
& \ \g_2^*\wedge(\g_1^*-\g_2^*)    \otimes 
 (\g_1 + 2\g_2)\wedge(\g_1+3\g_2)  \otimes 
1\otimes e_{12}\otimes (\g_1+\g_2)^*  
\\ 
\overset{=}\longmapsto \ \  
& -\g_{12}^* \otimes 
  \g_{12}    \otimes 
1 \otimes e_{12}\otimes (\g_1+\g_2)^*  
\\ 
\overset{\mu\otimes 1}\longmapsto \ \  
& -1\otimes e_{12}\otimes(\g_1+\g_2)^*.    
\end{aligned} 
\]
Proceeding in a similar fashion one obtains the following 
list for the maps $\phi_1^{IJ}$: 
\[
\begin{aligned}   
& \phi_1^{\{1,2,3,4\},\{1,2\}} &&\: \  
(\g_1 +\g_2)^* \otimes e_{1234}\otimes \g_{12}^*
&& \longmapsto &
- & 1\otimes e_{12}\otimes(\g_1+\g_2)^*  
\\ 
& \phi_1^{\{1,2,3,4\},\{1,3,4\}} & &\: \  
(\g_1 +\g_2)^* \otimes e_{1234}\otimes \g_{12}^*
& & \longmapsto & 
- & 1\otimes e_{134}\otimes \g_{12}^*  
\\ 
& \phi_1^{\{1,2,3,4\},\{2,3,4\}} & &\: \  
(\g_1 +\g_2)^* \otimes e_{1234}\otimes \g_{12}^*
& & \longmapsto &  
  & 1\otimes e_{234}\otimes \g_{12}^*.   
\end{aligned} 
\]
Finally, we give the corresponding 
list for the maps $\phi_0^{I,a}$: 
\[
\begin{aligned}  
& \phi_0^{\{1,2\},1} & &\: \  
 1 \otimes e_{12}\otimes (\g_1+\g_2)^*
& & \longmapsto &  
  & e_1  
\\ 
& \phi_0^{\{1,2\},2} & &\: \  
 1 \otimes e_{12}\otimes (\g_1+\g_2)^*
& & \longmapsto &  
- & e_2  
\\ 
& \phi_0^{\{1,3,4\},1} & &\: \  
 1 \otimes e_{134}\otimes \g_{12}^*
& & \longmapsto &  
  & e_1  
\\ 
& \phi_0^{\{1,3,4\},3} & &\: \  
 1 \otimes e_{134}\otimes \g_{12}^*
& & \longmapsto &  
-2 & e_3  
\\ 
& \phi_0^{\{1,3,4\},4} & &\: \  
 1 \otimes e_{134}\otimes \g_{12}^*
& & \longmapsto &  
  & e_4  
\\ 
& \phi_0^{\{2,3,4\},2} & &\: \  
 1 \otimes e_{234}\otimes \g_{12}^*
& & \longmapsto & 
  & e_2  
\\ 
& \phi_0^{\{2,3,4\},3} & &\: \  
 1 \otimes e_{234}\otimes \g_{12}^*
& & \longmapsto & 
-2& e_3  
\\ 
& \phi_0^{\{2,3,4\},4} & &\: \  
 1 \otimes e_{234}\otimes \g_{12}^*
& & \longmapsto & 
  & e_4.   
\end{aligned} 
\]
\end{example}

\section{The definition of $T_\bullet(\phi)$ and 
                           $T_\bullet(\phi)^{+}$}  

As in the previous section, here 
$\phi\: U_S \lra W$ is a representation over $\Bbbk$  
of a matroid $\bfM$ on a finite set $S$.  
The main goal 
is to present in Definition~\ref{D:T-complex} the 
description  
of the chain complex $T_\bullet(\phi)$ and of its 
augmentation $T_\bullet(\phi)^+$. 
We begin with the chains in homological 
degree $n$.

\begin{definition}\mlabel{D:chains-of-T} 
For $n\ge -1$ we define the space $T_n=T_n(\phi)$ by the 
formula 
\[
T_n = \bigoplus_{I \in\mathcal T_n(\bfM)} T_I(\phi)  
\]
and note that when $n=-1$ this simply yields 
$T_{-1}=\bigoplus_{a\in S}U_a = U_S$. Next, we set 
\[
\lambda = |S| - r_S^{\bfM} - 1 
\] 
and observe that for trivial reasons $T_n=0$ when  
$n\ge \lambda + 1$. Finally, we set $T_0^+ = W$, 
while for $n\ge 1$ we set $T_n^+ = T_{n-2}$. 
\end{definition}

Next we describe the differentials of our complexes. 

\begin{definition}\mlabel{D:differentials} 
Let $n\ge 0$. We define the map 
\[
\phi_n \: T_n \lra T_{n-1} 
\] 
by specifying that its restriction 
to the component $T_I(\phi)$ of $T_n$ be given by 
\[
\phi_n\big|_{T_I} \ = 
\sum_{\begin{smallmatrix} 
      J\in\mc T_{n-1}(\bfM) \\ 
      J\subseteq I
      \end{smallmatrix}} 
(-1)^{|J|} \phi_n^{IJ}.  
\]
When $n=0$ we will refer to the map $\phi_0$ also as the 
\emph{augmentation} or \emph{splicing homomorphism}. 
We also define for $n\ge 1$ the map 
\[
\phi_n^+ \: T_n^+ \lra T_{n-1}^+ 
\] 
by setting \ $\phi_1^+=\phi$, \  and by setting \ 
$\phi_n^+=\phi_{n-2}$ \ 
for $n\ge 2$.  
\end{definition}

\begin{remark}\mlabel{T:injective-differential}  
If $n\ge 1$ then by  
Remark~\ref{T:kernel-of-diagonal-map}(b) 
and the injectivity of \altref{E:diagonal-map} 
we have 
\[
\Ker\bigl(\phi_n\big|_{T_I}\bigr)=
\Ker(\Delta)\otimes\twedge U_I\otimes\twedge V_I^* = 0, 
\]
therfore the map $\phi_n\big|_{T_I}$ is injective. 
\end{remark}

\begin{remark}\mlabel{T:H-equivariant-differentials} 
Let $\phi$ be an $H$-equivariant representation of the 
matroid $\bfM$. Putting 
Remark~\ref{T:H-equivariant-phi-sub-n-IJ} together with 
the definition of the maps $\phi_n$ and $\phi_n^+$, 
it follows that for each $n\ge 0$ 
we have a canonically induced action 
of $H$ on the spaces $T_n$ and $T_n^+$ such that  
when $n\ge 1$ the maps $\phi_n$ and $\phi_n^+$ are 
$H$-equivariant. 
\end{remark}

We are now ready to achieve the 
main goal of this section.  

\begin{definition}\mlabel{D:T-complex}  
(a) We write $T_\bullet(\phi)$ for the sequence 
\[
T_\bullet(\phi)= 
0 \lra T_\lambda \overset{\phi_\lambda}\lra 
                 T_{\lambda-1} \lra \dots  \lra 
                 T_1 \overset{\phi_1}\lra T_0 \lra 0   
\]
and call it the \emph{T-complex} of the 
representation $\phi$. 

(b) We write $T_\bullet(\phi)^+$ for the sequence 
\[
\begin{CD}
T_\bullet(\phi)^+ = 0 \rightarrow   
T_{\lambda + 2}^+ 
@> \phi_{\lambda + 2}^+ >>  
T_{\lambda + 1}^+ 
\rightarrow \dots \rightarrow   
T_2^+ 
@> \phi_2^+ >> 
T_1^+ 
@> \phi_1^+ >> 
T_0^+ \rightarrow 0   
\end{CD}
\]
and call it the \emph{augmented T-complex} of the 
representation $\phi$.  
\end{definition}

\begin{remark} 
Since the maps $\phi_2^+\: T_2^+ \lra T_1^+$ and 
$\phi_1^+\: T_1^+ \lra T_0^+$ are just the maps 
$\phi_0\: T_0 \lra U_S$ and $\phi\: U_S \lra W$ 
respectively,  
the sequence $T_\bullet(\phi)^+$ is nothing but 
the shifted sequence $T_\bullet(\phi)[-2]$, spliced 
via the splicing homomorphism $\phi_0\: T_0 \lra U$
with the map $\phi\: U \lra W$. 
\end{remark}

\begin{example}\mlabel{E:uniform-matroid-2} 
Let $\phi\: U_S\lra W$ be a representation 
over $\Bbbk$ of the uniform 
matroid $\bfM$ of rank $r$ on $S$. Then 
using Example~\ref{E:uniform-matroid-1} we 
obtain for each $n\ge 0$ canonical isomorphisms  
\[
\begin{aligned} 
T_n(\phi) \quad &= \quad  
{\bigoplus}_{|I|=n+r+1} T_I(\phi)  
\\[+5pt]  
&= \quad 
{\bigoplus}_{|I|=n+r+1} 
(S_nV)^*\otimes\twedge U_I\otimes\twedge V^*   
\\[+5pt] 
&\cong \quad 
D_nV^*\otimes\wedge^{n+r+1}U_S\otimes\wedge^r V^*.  
\end{aligned}
\]
We use these to identify $T_n(\phi)$ with the space 
$D_nV^*\otimes\wedge^{n+r+1}U_S\otimes\wedge^r V^*$. 
It is immediate from the definitions that under these  
identifications for each $n\ge 0$ the vector space 
$T_n^{+}$ becomes precisely the vector space of chains 
$B_n$ of the Buchsbaum-Rim complex $B_\bullet(\phi)$ 
from \cite{BuRi64} (in the form  
described in \cite[Section 3]{ChTc03}). Furthermore, the  
maps $\phi_1^{+}$ and $-\phi_2^{+}$ transform 
exactly to the differentials $\phi$ and $\phi_2$ 
respectively, of the complex $B_\bullet(\phi)$. Finally,  
when $n\ge 3$ the morphism $\phi_n^{+}$ of 
$T_\bullet(\phi)^{+}$  becomes 
precisely equal to $(-1)^{n-2+r}$ times the differential 
$\phi_n$ of $B_\bullet(\phi)$. 

In summary, when $\phi$ represents 
a uniform matroid  the sequence $T_\bullet(\phi)^{+}$ 
is canonically isomorphic to the Buchsbaum-Rim complex 
$B_\bullet(\phi)$, in particular it is a resolution of 
$\Coker(\phi)$. 
\end{example} 

\begin{example}\mlabel{E:matrix-matroid-2}
Let $\bfM$ be the matroid represented by the map $\phi$ 
from Example~\ref{E:matrix-matroid}. Thus 
the sequence $T_\bullet(\phi)^{+}$ has the form 
\[
\begin{CD} 
0\lra T_{\{1,2,3,4\}} 
@> \phi_3^{+} >> 
T_{\{1,2\}}\oplus T_{\{1,3,4\}}\oplus T_{\{2,3,4\}} 
@> \phi_2^{+} >> 
U_S \overset{\phi_1^{+}}\lra W \lra 0.   
\end{CD}
\]
In Example~\ref{E:matrix-matroid-1} we 
gave a basis for each of the nonzero T-spaces 
over the field $\mathbb Q$. In 
Example~\ref{E:matrix-matroid-1-point-5} we described 
where these basis elements get mapped under each of the 
homomorphisms $\phi_n^{IJ}$. Putting this information 
together with the definitions of the maps $\phi_n^{+}$ 
yields that in these bases 
$T_\bullet(\phi)^{+}$ can be written as   
\[
\begin{CD} 
0 \lra \mathbb Q 
@> \phi_3^{+}=\phi_1 > 
{\begin{pmatrix} 
-1 \\ 1 \\ -1  
\end{pmatrix}} > 
\mathbb Q^3 
@> \phi_2^{+}=\phi_0 > 
{\begin{pmatrix} 
-1 & -1 &  0 \\ 
 1 &  0 & -1 \\ 
 0 &  2 &  2 \\ 
 0 & -1 & -1 
\end{pmatrix}} > 
\mathbb Q^4 
@> \phi_1^{+}=\phi > 
{\begin{pmatrix} 
1 & 1 & 1 & 1 \\ 
1 & 1 & 2 & 3 
\end{pmatrix}} >  
\mathbb Q^2 \lra 0.  
\end{CD}   
\]
It now straightforward to verify this is a complex 
that is a resolution of $\Coker(\phi)$. 
\end{example}

\section{The main properties of $T_\bullet(\phi)$ and 
                                $T_\bullet(\phi)^+$}

As in the previous two sections, throughout this one we 
fix a representation map $\phi\: U_S \lra W$ over a 
field $\Bbbk$ of a matroid $\bfM$ on a finite set $S$. 

Our goal is to present the statements of all key 
results about the complexes $T_\bullet(\phi)$ and 
$T_\bullet(\phi)^+$. Their proofs are with a few exceptions 
technically involved and are given in later sections.  
In particular, this and the next sections can be considered as 
a summary of all the main results of this paper, 
and will serve as a useful reference to a reader who is 
not interested in the technical details of the proofs but would 
like to get a good overview of the essential features of our
construction. However, since several interesting facts  
especially on the behaviour of multiplicity spaces have not 
been mentioned here, for those interested in a more detailed 
overview we recommmend to browse also through the 
statements of the results in Sections 8 through 11.   

We begin with the following fundamental assertion, which 
justifies our use of the word ``complex'' in the definitions 
of $T_\bullet(\phi)$ and $T_\bullet(\phi)^+$. 

\begin{theorem}\mlabel{T:the-T-bullets-are-complexes}  
The sequences $T_\bullet(\phi)$ and $T_\bullet(\phi)^+$ are 
complexes of vector spaces. 
\end{theorem}

Next, we describe the behaviour under the operation 
restriction of matroids. This property is a key 
ingredient in the proof of the second main result of 
this paper, Theorem~\ref{T:T-resolution-exactness}. 

\begin{theorem}\mlabel{T:exactness-combinatorial}  
Let \ $Y$ be a subset of \ $S$. 
For $n\ge -1$ we define  
\[
T_n\big|_Y \ = \ \ T_n(\phi)\big|_Y \ = 
\bigoplus_{\begin{smallmatrix} 
           I\in\mc T_n(\bfM) \\ I\subseteq Y 
           \end{smallmatrix} } 
T_I(\phi). 
\]
Similarly, we set \ $T_0^+\big|_Y = T_0^+=W$, and for 
$n\ge 1$ we set \ $T_n^+\big|_Y = T_{n-2}\big|_Y$. 
Then: 

{\rm (a)} The sequence of vector space maps 
\[
\begin{CD} 
T_\bullet(\phi)\big|_Y = 0 \rightarrow T_\lambda\big|_Y   
@> \phi_\lambda >>   
T_{\lambda-1}\big|_Y \rightarrow \dots \rightarrow T_1\big|_Y           
@> \phi_1 >>  
T_0\big|_Y \rightarrow 0  
\end{CD} 
\]
is a subcomplex of \ $T_\bullet(\phi)$, and we have \  
$T_\bullet(\phi|Y)= T_\bullet(\phi)\big|_Y$. 

{\rm (b)} The sequence of vector spaces and 
homomorphisms \  
$T_\bullet(\phi)^+\big|_Y$ defined as 
\[
\begin{CD}
0 \rightarrow T_{\lambda + 2}^+\big|_Y 
@> \phi_{\lambda + 2}^+ >>  
T_{\lambda + 1}^+\big|_Y  
\lra \dots \lra  
T_2^+\big|_Y  
@> \phi_2^+ >> 
T_1^+\big|_Y  
@> \phi_1^+ >> 
T_0^+\big|_Y \rightarrow 0   
\end{CD}
\]
is a subcomplex of \ $T_\bullet(\phi)^+$, and we have \ 
$T_\bullet(\phi|Y)^+ = \ T_\bullet(\phi)^+\big|_Y$. 
\end{theorem} 

The content of Theorem~\ref{T:exactness-combinatorial} 
is probably best stated informally: 
to obtain the component in homological degree $n$ 
of the (augmented) T-complex 
of \ $\phi|Y$, one simply needs 
to select from the degree $n$ component 
of the corresponding complex 
of $\phi$ those T-spaces that are indexed by 
subsets of $Y$. 

We also record the behavior under 
the operation sum of matroids:  

\begin{theorem}\mlabel{T:T-bullet-of-sum} 
If \ $S_1,\dots, S_k$ are subsets of \ $S$ such that \  
$\bfM = \bfM|S_1 + \dots + \bfM|S_k$, then we have 
a canonical decomposition 
\[
T_\bullet(\phi) = 
T_\bullet(\phi|S_1) \oplus \dots \oplus 
T_\bullet(\phi|S_k)  
\]
of the T-complex of \ $\phi$ as a direct sum of 
subcomplexes. 
\end{theorem} 

Next comes the behavior under the operation contraction 
of matroids. This result together with 
Theorem~\ref{T:multiplicity-exactness} are the key  
ingredients in the proof of the first main result of 
this paper, Theorem~\ref{T:exactness-of-complex}. 

\begin{theorem}\mlabel{T:contraction-embedding} 
Let \ $Y$ be a subset of \ $S$ such that \ $S\smsm Y$ is 
independent in \ $\bfM$. 

There exists a canonical morphism of complexes \ 
$
(\pi.Y)_\bullet^\phi\: 
T_\bullet(\phi.Y) \lra T_\bullet(\phi)
$ 
\ {\rm(}described in a very explicit combinatorial 
way in Section~13{\rm)} which is injective and 
an isomorphism in homology. 
\end{theorem} 

We finally arrive at the statement of  
the first main result of this paper, which asserts 
the acyclicity of our complexes. 

\begin{theorem}\mlabel{T:exactness-of-complex}
The T-complex \ 
$T_\bullet(\phi)$ is a resolution of \ $\Ker\phi$, 
and the augmented T-complex \ $T_\bullet(\phi)^+$ 
is a resolution of \ $\Coker \phi = W/V$. 
\end{theorem}

As one of the important consequences of this theorem we 
have that the numerical characteristics of our complexes 
are independent from the representation $\phi$. 

\begin{theorem}\mlabel{T:independence-from-phi} 
For any subset $A\subseteq S$ 
the dimension of the multiplicity space $S_A(\phi)$ is an 
invariant of the matroid $\bfM|A$, and does not depend 
on the representation map $\phi$. In particular, the 
length and the ranks of the components of the T-complex 
$T_\bullet(\phi)$ are invariants of the matroid $\bfM$, 
and do not depend on the representation map $\phi$.  
\end{theorem}

We conclude this section with a result which 
will be of interest to group representation theorists. 
This theorem simply brings together properties that 
we have already observed in 
Remarks~\ref{T:H-action-on-T-flats-and-parts}, 
\ref{T:H-action-on-multiplicity-spaces}, and 
\ref{T:H-equivariant-differentials}. 

\begin{theorem}\mlabel{T:group-action} 
Let $\phi$ be an $H$-equivariant representation of $\bfM$. 

Then for each $n\ge 0$ we have canonically induced 
linear actions of $H$ on the vector spaces $T_n$ 
and $T_n^+$ such that the differentials of the 
complexes $T_\bullet(\phi)$ and $T_\bullet(\phi)^+$ 
are $H$-equivariant.  \qed
\end{theorem}

While leaving a more detailed investigation for a separate 
paper, we remark here that 
a situation as in Theorem~\ref{T:group-action} arises 
for example when $H=\Sigma_n$ is the symmetric group,  
the space $W$ is an  
irreducible representation of $H$ corresponding to some 
Young diagram $\Gamma$ with $n$ boxes, the set $S$ is 
the collection of all row-standard tableaux of shape 
$\Gamma$ with entries from $1$ to $n$, and the map 
$\phi$ sends a row-standard tableau to the corresponding 
uniquely determined element of the irreducible 
representation $W$.

\section{Free resolutions of multigraded modules} 

Throughout this section 
$R=\Bbbk[x_1,\dots, x_m]$ is a polynomial ring over 
a field $\Bbbk$ with the standard $\mathbb Z^m$-grading,  
$L$ is a Noetherian $\mathbb Z^m$-graded (multigraded) 
$R$-module,  
\[
E\overset{\Phi}\lra  G\lra L \lra 0 
\]
is a finite free multigraded presentation of $L$,  
and $S$ is a multihomogeneous basis of $E$.  
The main goal is to introduce in 
Definition~\ref{D:T-resolution} a finite free complex 
of multigraded $R$-modules $T_\bullet(\Phi, S)$ which 
we call the \emph{T-resolution of the pair} $(\Phi, S)$. 
In the second main result of this paper, 
Theorem~\ref{T:T-resolution-exactness}, we show that the 
T-resolution of the pair $(\Phi, S)$ is a finite free 
resolution of the $R$-module $L$. 

Before we proceed with the statements, we introduce 
some more notation. First, we consider the field $\Bbbk$ as an 
$R$-module via the canonical projection $R\lra \Bbbk$ that 
sends each variable $x_i$ to the identity element $1\in \Bbbk$. 
The set $\mathbb Z^m$ has a partial ordering 
$\preceq$ defined for sequences $\bfa=(a_1,\dots,a_m)$ and 
$\bfb=(b_1,\dots,b_m)$ by the formula 
\[
\bfa\preceq\bfb \quad\iff\quad 
a_i\le b_i \text{ for every } i.  
\] 
With this partial order $\mathbb Z^m$ is a lattice, 
the join (or lcm) 
of $\bfa$ and $\bfb$ being their componentwise maximum 
\[
\lcm(\bfa,\bfb) = \bfa\vee\bfb = 
\bigr(\max(a_1,b_1),\dots, \max(a_m,b_m)\bigl). 
\] 
Similarly, the meet of $\bfa$ and $\bfb$ is their 
componentwise minimum. When $\bfa\in\mathbb N^m$ 
we write $x^{\bfa}$ for the monomial 
$x_1^{a_1}\cdots x_m^{a_m}\in R$. 
If $z$ is a multihomogeneous element inside  
a $\mathbb Z^m$-graded 
$R$-module, we write $\deg z$ for its multidegree. 
In particular, $\deg(x^{\bfa})=\bfa$. More generally, 
if $A=\{z_1,\dots,z_k\}$ is a collection of multihomogeneous 
elements in a multigraded $R$-module, we set 
\[
\deg A = \lcm(\deg z_1,\dots,\deg z_k). 
\]
For example, since $S$ is a collection of multihomogeneous 
elements in the multigraded $R$-module $E$, for any 
subset $I\subseteq S$ we have that $\deg I$ is the 
componentwise maximum of the multidegrees of the 
elements of $I$. 

We are now ready to introduce the key ingredients  
from which the T-resolution of $(\Phi,S)$ is built. 
First is the representation of the matroid that 
governs the linear algebraic structure of our 
resolution. 

\begin{definition}\mlabel{D:definition-of-phi}
(a) We define the $\Bbbk$-vector space $W$ as 
\[
W = \Bbbk\otimes_R G. 
\]

(b) We note that the elements $\{1\otimes_R a\mid a\in S\}$ 
form a basis of the $\Bbbk$-vector space $\Bbbk\otimes_R E$, 
and we set for each $a\in S$ 
\[
e_a = 1\otimes_R a. 
\] 
We use this to identify canonically $\Bbbk\otimes_R E$ with 
the space $U_S$. 

(c) We set $\phi= 1\otimes_R \Phi$. Thus we have a 
$\Bbbk$-vector space homomorphism 
\[
\phi\: U_S \lra W, 
\]
and we write $\bfM=\bfM(\Phi,S)$ for the matroid on the set 
$S$ represented over $\Bbbk$ by $\phi$. 
\end{definition}

Next we introduce the components out of which we will 
construct the chains and the differentials for our 
resolution. 

\begin{definition}\mlabel{D:T-resolution-components}
Let $n\ge -1$, and let $I$ be an element in $\mc T_n(\bfM)$. 

(a) We define  
\[
T_I(\Phi, S)= R\otimes_\Bbbk T_I(\phi).
\] 
A canonical multigrading on this free $R$-module 
is defined by the formula  
\[
\deg (z\otimes v) = \deg z + \deg I
\] 
for any vector $v\in T_I(\phi)$ and any monomial $z\in R$. 

(b) When $n\ge 0$ and $J$ is an element of $\mc T_{n-1}(\bfM)$ 
contained inside $I$,  we 
define the canonical morphism of multigraded free $R$-modules  
\[
\Phi_{n+2}^{IJ}\: T_I(\Phi, S) 
                  \lra         
                  T_J(\Phi, S) 
\] 
by the formula 
\[
\Phi_{n+2}^{IJ}(z\otimes v)= 
x^{\deg I - \deg J}z\otimes \phi_n^{IJ}(v).  
\]
\end{definition}

Finally, here is the definition of the T-resolution of 
the pair $(\Phi,S)$. 

\begin{definition}\mlabel{D:T-resolution}  
(a) We define $T_0(\Phi,S)=G$, 
while for $n\ge 1$ we define the multigraded free $R$-module  
\[ 
T_n(\Phi,S)= 
\bigoplus_{I\in\mc T_{n-2}(\bfM)} T_I(\Phi, S).
\] 
Note that the module $T_1(\Phi, S)$ is simply the 
free $R$-module $E$. 

(b) For $n\ge 1$ 
we define the morphism of multigraded free $R$-modules  
\[
\Phi_n \: T_n(\Phi, S) \lra T_{n-1}(\Phi, S) 
\]
as follows. When $n\ge 2$ we define $\Phi_n$ 
by requiring that its restriction to the component 
$T_I(\Phi, S)$ of the module $T_n(\Phi, S)$ be given 
by the formula 
\[
\Phi_n\big|_{T_I(\Phi, S)} = 
\sum_{\begin{smallmatrix} 
      J\in \mc T_{n-3}(\bfM) \\ J\subseteq I
      \end{smallmatrix} } 
(-1)^{|J|}\Phi_n^{IJ}.  
\] 
When $n=1$ we simply set $\Phi_1 = \Phi$. 

(c) With $\lambda=|S|-r_S^{\bfM}+1$, we define 
the sequence $T_\bullet(\Phi, S)$ as  
\[
T_\bullet(\Phi, S) = 
0\rightarrow 
T_\lambda(\Phi,S) \overset{\Phi_\lambda}\lra 
T_{\lambda -1}(\Phi, S) \rightarrow \dots \rightarrow 
T_1(\Phi, S) \overset{\Phi_1}\lra 
T_0(\Phi, S) \rightarrow 0   
\]
and call it the \emph{T-resolution of the presentation 
$\Phi$ with respect to the basis $S$}, or for short, 
the \emph{T-resolution of $(\Phi, S)$}.  
\end{definition}

\begin{remarks}\mlabel{E:uniform-matroid-3}  
(a) Suppose that a basis $S$ can be chosen so that 
the matroid $\bfM(\Phi,S)$ is uniform (in the 
terminology of \cite{ChTc03} this is the case of  
the map $\Phi$ having uniform rank). We saw in 
Example~\ref{E:uniform-matroid-2} that in such 
a situation the complex $T_\bullet(\phi)^{+}$ is 
canonically isomorphic to the Buchsbaum-Rim 
complex. It is now straightforward to verify 
from the definitions that this isomorphism 
carries through to give an isomorphism 
of the T-resolution $T_\bullet(\Phi,S)$ with 
the Taylor complex $T_\bullet(\Phi)$ from 
\cite{ChTc03}. In particular 
the T-resolution of $(\Phi,S)$ is a free 
resolution of $L=\Coker(\Phi)$, and when $\Phi$ is 
the standard minimal presentation of $R/I$ for a monomial 
ideal $I$ and $S$ is the basis of $E$ whose elements map 
to the minimal generators of $I$, then the T-resolution 
of the pair $(\Phi,S)$ recovers the usual Taylor 
resolution of $R/I$. 

(b) It is clear from part (a) that, just as the 
Taylor resolution, the complex  
$T_\bullet(\Phi,S)$ is in general not minimal. 
Furthermore, the ranks of its components depend 
substantially on the choice of the basis $S$. It 
is a very interesting open problem to investigate 
the properties of $T_\bullet(\Phi,S)$ under a 
``generic'' choice of $S$.   
\end{remarks}

The second main result of this paper is:  

\begin{theorem}\mlabel{T:T-resolution-exactness} 
The sequence $T_\bullet(\Phi, S)$ is a complex, and is a 
finite free multigraded resolution of the $R$-module 
$L=\Coker(\Phi)$. 
\end{theorem}

\begin{proof}
As $\Phi_1=\Phi$, we have $\Coker\Phi_1 = L$. Thus to 
prove the theorem it suffices to prove that the 
sequence $T=T_\bullet(\Phi, S)$ is a complex, and is 
acyclic, i.e.~it has zero homology in positive 
homological degree.   
Since all the maps in the sequence $T=T_\bullet(\Phi, S)$ 
are morphisms of multigraded modules, 
$T$ splits as a direct sum of sequences  
of multigraded vector spaces   
$T=\bigoplus_{\bfa \in \mathbb Z^n}T_{\bfa}$, where each 
$T_{\bfa}$ is a sequence of vector spaces, with all their 
vectors multihomogeneous of the same multidegree $\bfa$. 
Thus it suffices to show that each 
sequence $T_{\bfa}$ is a complex of vector spaces and  
has zero homology in positive (homological) degree. 
Let $I_{\bfa}$ be the subset of $S$ consisting of all 
elements of $S$ with multidegree $\preceq\bfa$. Since for 
an  element $A\in \mc T_n(\bfM)$ the $R$-module 
$T_A(\Phi, S)$ contributes to $T_{\bfa}$ if and only if  
$A\subseteq I_{\bfa}$, in which case it contributes 
precisely $x^{\bfa -\deg A}\otimes T_A(\phi)$, it is 
immediate that $T_{\bfa}$ can be canonically identified 
with a subcomplex of the restricted complex 
$T^+=T_\bullet(\phi|I_{\bfa})^+$, and that $T_{\bfa}$ and 
$T^+$ may possibly differ only in 
homological degree $0$. Thus $T_{\bfa}$ is a complex, and 
the desired acyclicity is now immediate from  
Theorem~\ref{T:exactness-of-complex} applied to 
the restricted complex $T^+$. 
\end{proof}

As an interesting immediate consequence we obtain a nice 
upper bound on the projective dimension of a Noetherian 
multigraded $R$-module $L$. Further applications 
of Theorem~\ref{T:T-resolution-exactness} will appear 
in~\cite{ChTc} and~\cite{Tc}.

\begin{theorem}\mlabel{T:projective-dimension} 
Let $L$ be a multigraded Noetherian $R$-module of 
rank $r$, and let $\beta_0$ and $\beta_1$ be its zeroth 
and first Betti numbers, respectively. 

Then \ $\pd_R L \le \beta_1 - \beta_0 + r + 1$. 
\end{theorem}

\begin{proof} 
Since in a minimal free presentation $\Phi$ 
of $L$ one has $\beta_0=\rank G$,
and $\beta_1 = \rank E$, the result is immediate in view 
of Theorem~\ref{T:T-resolution-exactness} and the fact 
that the resolution $T_\bullet(\Phi,S)$ has length at most  
$\lambda = \beta_1 - (\beta_0 - r) + 1$. 
\end{proof}

We conclude this section with an example. 

\begin{example}\mlabel{E:matrix-matroid-3} 
Let $R=\mathbb Q[x,y,z]$, and let 
$G$ be the free multigraded $R$-module on a basis 
$\{g_1,g_2\}$ with 
\[
\begin{aligned} 
\deg(g_1) &= (1,1,0), \\  
\deg(g_2) &= (0,0,1).
\end{aligned}
\] 
Let $E$ be the free multigraded $R$-module on a basis 
the set $S=\{1,2,3,4\}$ with 
\[
\begin{aligned} 
\deg(1) &= (3,1,1), &   
\deg(2) &= (1,3,1),   \\  
\deg(3) &= (1,1,3), &    
\deg(4) &= (1,2,2).
\end{aligned} 
\] 
Let $\Phi\: E\lra G$ be the 
morphism given in these bases by the matrix 
\[
\begin{pmatrix} 
x^2z & y^2z & z^3    & yz^2 \\ 
x^3y & xy^3 & 2xyz^2 & 3xy^2z 
\end{pmatrix}. 
\] 
Let $\g_1=1\otimes g_1$ and $\g_2=1\otimes g_2$ 
be the corresponding basis elements of 
the vector space $W=\mathbb Q\otimes_R G$. Thus the map 
\[
\phi=1\otimes_R \Phi \: U_S \lra W
\] 
is the map we considered in Examples~\ref{E:matrix-matroid}  
and~\ref{E:matrix-matroid-0}, 
and therefore the T-flats of the matroid $\bfM(\Phi,S)$ have 
multidegrees 
\[
\begin{aligned} 
\deg\{1,2\}     &= (3,3,1), & 
\deg\{1,3,4\}   &= (3,2,3),   \\ 
\deg\{2,3,4\}   &= (1,3,3), &   
\deg\{1,2,3,4\} &= (3,3,3).
\end{aligned} 
\] 
It follows from Example~\ref{E:matrix-matroid-2} 
that for the T-resolution $T_\bullet(\Phi,S)$ we have   
\[
\begin{CD} 
0 \lra T_3  
@> \Phi_3 > 
{\left(\begin{smallmatrix} 
-z^2 \\ y \\ -x^2  
\end{smallmatrix}\right)} > 
T_2  
@> \Phi_2 > 
{\left(\begin{smallmatrix} 
-y^2 & -yz^2 &  0 \\ 
 x^2 &  0    & -z \\ 
 0   & 2x^2y & 2y^2 \\ 
 0   & -x^2z & -yz 
\end{smallmatrix}\right)} > 
T_1  
@> \Phi_1=\Phi > 
{\left(\begin{smallmatrix} 
x^2z & y^2z & z^3    & yz^2 \\ 
x^3y & xy^3 & 2xyz^2 & 3xy^2z 
\end{smallmatrix}\right)} >  
T_0 \lra 0.  
\end{CD}   
\]
It can now be verified directly that this 
is a free (and in this case minimal) multigraded 
resolution of the $R$-module 
$L=\Coker(\Phi)$. 
\end{example}

\section{T-flats of minors}

In this section $\bfM$ is a matroid on a finite set $S$. 
We describe in Theorems~\ref{T:T-flats-of-minors} and 
\ref{T:T-partition-of-contraction} 
the relationship between the T-flats 
of $\bfM$ and the T-flats of the minors of $\bfM$. 
All results stated here are due to 
Tutte~\cite{Tu58}, 
and will play an essential role in the 
proofs of the main theorems of our paper. 
For completeness we have included their short proofs 
(which may differ from the arguments used 
in~\cite{Tu58}).

\begin{theorem}\mlabel{T:T-flats-of-minors}  
{\rm(\cite{Tu58})}.
Let \ $Y$ be a subset of \ $S$ and let $A$ be a subset of 
\ $Y$. Then:  
\begin{enumerate}
\item The set $A$ is a T-flat of \ $\bfM|Y$ if and only if $A$ is a 
      T-flat of \ $\bfM$. In that case the T-parts of $A$ in  
      $\bfM$ and in $\bfM|Y$ coincide. 

\item The set $A$ is a T-flat of \ $\bfM . Y$ if and only if 
      $A=I\cap Y$ for some T-flat $I$ of $\bfM$. If $A$ is a T-flat 
      of \ $\bfM.Y$ then in fact $A=B\cap Y$    
      where $B$ is the T-flat of \ $\bfM$ given by 
      $B=S\smsm (Y\smsm A)^{\bsfC_{\bfM^*}}$.  
\end{enumerate}
\end{theorem}

\begin{proof} 
Since the circuits of \ $\bfM|Y$ are precisely the circuits of \ 
$\bfM$ that are contained in $Y$, part (1) is immediate from 
Remark~\ref{T:T-remarks}(c). 

Next, $A$ is a T-flat of $\bfM.Y$ precisely when $Y\smsm A$ is a 
flat of $(\bfM.Y)^*=\bfM^*|Y$. This happens if and only if there 
is a flat $I'$ of $\bfM^*$ with  $Y\smsm A= I'\cap Y$, in which 
case also $Y\smsm A=(Y\smsm A)^{\bsfC_{\bfM^*}}\cap Y$. This in 
turn occurs    if and only if $A=I\cap Y=B\cap Y$ where 
$I=S\smsm I'$ and  $B=S\smsm (Y\smsm A)^{\bsfC_{\bfM^*}}$. 
\end{proof}

\begin{remarks}\mlabel{T:maximal-T-flat}
{\rm(\cite{Tu58})}.
Let $A$ be a T-flat of $\bfM.Y$ of level $n$, and  
$B=S\smsm (Y\smsm A)^{\bsfC_{\bfM^*}}$. 

(a) Since $(Y\smsm A)^{\bsfC_{\bfM^*}}$ is the unique 
smallest flat of $\bfM^*$ containing $Y\smsm A$, 
the set $B$ is the unique maximal T-flat 
of $\bfM$ whose intersection with $Y$ is $A$.  

(b) Since  
$r^{\bfM^*}_{S\smsm B} = r^{\bfM^*}_{Y\smsm A} = 
 r^{\bfM^*|Y}_{Y\smsm A} = r^{(\bfM.Y)^*}_{Y\smsm A} = 
 r^{(\bfM.Y)^*}_Y - n - 1 = r^{\bfM^*|Y}_Y - n - 1 = 
 r^{\bfM^*}_Y - n - 1$, 
we obtain from \altref{E:level-formula} the equality  
\[
\ell^{\bfM}_B = n + r^{\bfM^*}_S - r^{\bfM^*}_Y.
\]  
In particular $\ell_B^{\bfM}\ge n$, and 
the set $S\smsm Y$ is independent in $\bfM$ 
if and only if the set $B$ is a T-flat in $\bfM$ of 
level exactly $n$. 
\end{remarks}

\begin{theorem}\mlabel{T:T-partition-of-contraction} 
{\rm(\cite{Tu58})}. 
Let $A\subseteq Y$ be a T-flat of level $n\ge 1$ in 
the contracted matroid \ 
$\bfM.Y$. Let 
$B=S\smsm (Y\smsm A)^{\bsfC_{\bfM^*}}$, and let 
$J_1,\dots,J_k$ be those of the T-parts of \ $B$ in 
$\bfM$ that intersect \ $Y$ nontrivially. Let 
$I_i=J_i\cap Y$, let $A_i=A\smsm I_i$, and let 
$B_i=B\smsm J_i$.  

Then 
the sets $I_i$ are all the T-parts of $A$ in $\bfM.Y$ 
and $B_i=S\smsm (Y\smsm A_i)^{\bsfC_{\bfM^*}}$    
for each $i$.  
\end{theorem}

\begin{proof} 
Let $B_i=B\smsm J_i$ and $A_i=B_i\cap Y=A\smsm I_i$. 
Thus by Theorem~\ref{T:T-flats-of-minors} the set 
$A_i$ is a T-flat of $\bfM.Y$ and also  
$A_i=B'_i\cap Y$ for 
$B'_i= S\smsm(Y\smsm A_i)^{\bsfC_{\bfM^*}}$, 
in particular 
$
(Y\smsm A)^{\bsfC_{\bfM^*}} \ne 
(Y\smsm A_i)^{\bsfC_{\bfM^*}}
$. 
Note furthermore that 
$Y\smsm A_i=(S\smsm B_i)\cap Y$ and $S\smsm B_i$ is a 
cover in $\bfM^*$ of the flat $
S\smsm B=(Y\smsm A)^{\bsfC_{\bfM^*}}$. 
Since 
\[
S\smsm B = (Y\smsm A)^{\bsfC_{\bfM^*}}
          \subsetneq (Y\smsm A_i)^{\bsfC_{\bfM^*}}
          \subseteq S\smsm B_i,
\] 
it follows that 
$S\smsm B_i=(Y\smsm A_i)^{\bsfC_{\bfM^*}}$. 
Therefore each $A_i$ is a maximal T-flat of $\bfM.Y$ 
properly contained in $A$. 
Since $A=I_1\sqcup\dots\sqcup I_k$, the desired 
conclusion is 
immediate from Remark~\ref{T:T-remarks}(b). 
\end{proof} 

The following special case is of particular importance to us. 

\begin{corollary}\mlabel{T:T-flats-on-S-a} 
Let $a\in S$ be an element with $\{a\}$ independent in 
$\bfM$, and let $S_a=S\smsm \{a\}$. Let $A\subseteq S_a$, 
and let $B=S\smsm(S_a\smsm A)^{\bsfC_{\bfM^*}}$. 
Then $A$ is a T-flat of \ $\bfM.S_a$ if and only if 
for the T-flat $B$ of \ $\bfM$ we have 
$B=A$ or $B=A\cup\{a\}$. 

Furthermore, when $A$ is a T-flat of \ $\bfM.S_a$ we  
have: 
\begin{enumerate} 
\item  If $B=A$ then $A\cup\{a\}$ is not a T-flat of \ 
       $\bfM$, the T-flats of \ $\bfM$ inside $B$ coincide  
       with the T-flats of \ $\bfM.S_a$ inside $A$, and 
       $r^{\bfM}_{A\cup\{a\}}=r^{\bfM}_A + 1$.  

\item  If $B=A\cup\{a\}$ and $A$ is not a T-flat of \ 
       $\bfM$, then $\{a\}$ is not a T-part of $B$ and 
       $r^{\bfM}_{A\cup\{a\}}=r^{\bfM}_A$. 

\item  If $B=A\cup\{a\}$ and $A$ is a T-flat of \ $\bfM$, 
       then $\{a\}$ is a T-part of $B$ and 
       $r^{\bfM}_{A\cup\{a\}}=r^{\bfM}_A$.     \qed   
\end{enumerate}  
\end{corollary}

\section{The structure of the T-parts of a T-flat}

In this section $\bfM$ is a matroid on a finite set $S$. 
The main assertions on the structure of T-parts are 
Theorem~\ref{T:complement-is-independent}, 
Theorem~\ref{T:components-unite}, and 
Theorem~\ref{T:product-is-admissible}. The last of these  
theorems is essentially  due to Tutte and is a 
straightforward consequence of 
the results in \cite{Tu58}. We proceed by  
first stating our three theorems   
together with some remarks, and then we present  
their proofs including for completeness a proof 
of Theorem~\ref{T:product-is-admissible}.

\begin{theorem}\mlabel{T:complement-is-independent} 
Let $I$ be a T-flat in $\bfM$ of level $n$, and let 
$J$ be a T-part of $I$.  
\begin{enumerate} 
\item If $J$ is independent then 
      $r_J + r_{I\smsm J} = r_I + 1$. 
\item If $J$ is not independent then $J$ is a circuit,   
      $n\ge 1$, and $r_J+ r_{I\smsm J} = r_I$. 
\item If $J'$ is a proper subset of $J$ then 
      $r_{J'} + r_{I\smsm J} = r_{J'\cup(I\smsm J)}$.  
\end{enumerate} 
\end{theorem}

\begin{remark}\mlabel{T:rank-1-intersection-a} 
Let $\phi\: U_S \lra W$ be a presentation of $\bfM$ and 
let $J$ be a T-part of the T-flat $I$. 
In that setting the rank conditions 
from parts (1), (2), and (3) 
of Theorem~\ref{T:complement-is-independent} are 
equivalent to 
\begin{enumerate}  
\item $\dim_\Bbbk V_J\cap V_{I\smsm J} = 1$ when $J$ 
      is independent, 
\item $V_J\cap V_{I\smsm J} = 0$ when $J$ is a circuit, 
      and 
\item $V_{J'}\cap V_{I\smsm J}=0$ when $J'$ is a proper 
      subset of $J$,   
\end{enumerate} 
respectively. 
\end{remark} 

The following useful remarks will be needed 
several times in later sections. 
Also, they should give the reader a better 
feel for the role played by T-parts when it 
comes down to linear algebra. 

\begin{remarks}\mlabel{T:rank-1-intersection-b} 
Let $\phi\: U_S \lra W$ be a presentation of $\bfM$,  
let $I$ be a T-flat of level $n$, and  
let $I_1,\dots, I_k$ be all the T-parts  
of $I$. For each $1\le i\le k$ let 
\[
u_i=\sum_{a\in I_i}c_{ia} e_a
\] 
be an arbitrary element of the vector space $U_{I_i}$. Also, 
let 
\[
V_i = V_{I_i}\cap V_{I\smsm I_i}.      
\] 
Furthermore, for each $1\le i\le k$ let 
\[
v_i = \sum_{a\in I_i} d_{ia}e_a 
\]
be an element in $U_{I_i}$ 
such that the vector 
$w_i=\phi(v_i)$ is a basis of the $1$-dimensional space 
$V_i$ in case $I_i$ is independent, and such that  
$v_i$ is a basis for the ($1$-dimensional) 
kernel of $\phi|I_i$ in case $I_i$ is a circuit. Then: 

(a) It is clear that the vectors $v_i$ are uniquely 
determined up to a nonzero scalar multiple.

(b) If $I_i$ is a circuit then $\phi(v_i)=0$ and therefore 
$d_{ia}\ne 0$ for each $a\in I_i$. 

(c) If $I_i$ is independent then $0\ne \phi(v_i)\in V_i$, 
hence there is a $0\ne v'\in U_{I\smsm I_i}$ such that 
$\phi(v_i + v')=0$; therefore by 
Theorem~\ref{T:complement-is-independent}(3) and 
Remark~\ref{T:rank-1-intersection-a} it follows again that 
$d_{ia}\ne 0$ for each $a\in I_i$.  

(d) If $n\ge 1$ and 
$\phi(\sum_i u_i)=0$, then in view of the T-partition 
$I=I_1\sqcup\dots\sqcup I_k$ we have 
$\phi(u_i)\in V_i$ for each $i$,     
hence every $u_i$ is a multiple of $v_i$. In particular,  
any syzygy of $\phi$ on $U_I$ is a linear combination of 
the vectors $v_i$.

(e) If $n\ge 1$ and $\phi(\sum_i u_i)=0$, then   
it is immediate from parts (b), (c), and (d) above that 
whenever $u_i\ne 0$ then also $c_{ia}\ne 0$ for 
every $a\in I_i$. 
\end{remarks}

\begin{remark}\mlabel{T:rank-1-intersection-c} 
Let $I$ be a T-flat of level $n\ge 2$ with  
T-partition $I=I_1\sqcup \dots \sqcup I_k$. Then for each 
$j$ the set $I\smsm I_j$ is a T-flat of level $n-1$ 
and can be written as 
\[
I\smsm I_j = I_1\sqcup \dots \sqcup I_{j-1}\sqcup I_{j+1}
                \sqcup \dots \sqcup I_k,  
\]
however in general this is \emph{not} the  
T-partition of $I\smsm I_j$. The next theorems are aimed 
at clarifying this issue further. 
\end{remark}

\begin{theorem}\mlabel{T:components-unite}
Let $I$ be a T-flat in $\bfM$, and let $J$ be 
a T-flat in $\bfM$ of level \ $n$ properly contained 
inside $I$. 

The T-flat $J$ is a disjoint union of 
T-parts of $I$. If in addition $n\ge 1$, then each 
T-part of $J$ is the disjoint union of T-parts  of $I$. 
\end{theorem}

\begin{corollary}\mlabel{T:contained-in-circuit}
Let $I$ be a T-flat in $\bfM$, and let $J'$ be a 
T-part of $I$. There exists a circuit $C$ inside $I$ 
that contains $J'$.  
\end{corollary}

We will also need the following observations in 
Section~11. 

\begin{remarks}\mlabel{T:span-of-T-parts} 
Let $\phi\: U_S \lra W$ be a representation of $\bfM$, 
let $I=I_1\sqcup\dots\sqcup I_k$ be the T-partition of 
a T-flat $I$, and let the vectors $v_i\in U_{I_i}$ be 
chosen as in Remark~\ref{T:rank-1-intersection-b}. 
Let $J$ be a T-flat of level $n$ 
properly contained inside $I$, 
let $J'$ be an independent T-part of $J$, and let 
$v\in U_{J'}$ be such that the vector $w=\phi(v)$ is 
a basis of the ($1$-dimensional) space 
$V'=V_{J'}\cap V_{J\smsm J'}$.  
By Corollary~\ref{T:contained-in-circuit} 
the set $J'$ is contained 
in a circuit $C$ inside $J$, and let 
\[
u=\sum_{a\in C} d_a e_a
\] 
be a basis vector of the $1$-dimensional kernel 
of the map $\phi|C$. Then: 

(a) By Remarks~\ref{T:rank-1-intersection-b}(a) 
and~\ref{T:rank-1-intersection-b}(b),   
the vectors $v$ and $u$ are uniquely determined up 
to a multiple by a nonzero scalar, and 
we have $d_a \ne 0$ for each $a\in C$. 

(b) Since $\phi(u)=0$, it follows from part (a) that 
the element $u'=\sum_{a\in J'}d_a e_a$ is a multiple 
of $v$. 

(c) If $n\ge 1$ then by Theorem~\ref{T:components-unite}
the set $J'$ is the disjoint union of T-parts of $I$, say 
$J'=I_{i_1}\sqcup\dots\sqcup I_{i_t}$. Therefore, 
by Remark~\ref{T:rank-1-intersection-b}(d) for each 
$1\le p\le t$ the partial sum 
$u'_p=\sum_{a\in I_{i_p}}d_a e_a$ is a multiple 
of $v_{i_p}$. 

(d) Suppose $n\ge 1$. It is immediate from 
parts (b) and (c) above that $v$ is a linear combination 
of the vectors $v_{i_p}$.  
Applying $\phi$ to that linear combination yields that $w$  
(and hence also the entire space $V'$) lies inside the 
subspace $V_1+\dots+V_k$ of the vector space $W$.  
\end{remarks}

\begin{theorem}\mlabel{T:product-is-admissible} 
{\rm(\cite{Tu58})}. 
Let $I=I_1\oplus\dots\oplus I_k$ with $k\ge 2$ and each 
set $I_i$ nonempty. 
\begin{enumerate} 
\item We have 
      $\ell_I = \ell_{I_1} + \dots + \ell_{I_k} + k - 1$. 

\item The set $I$ is a T-flat if and only if all the sets 
      $I_i$ are T-flats. In that case a subset $J$ is a 
      T-part of $I$ if and only if  \underline{either} 
      $J=I_i$ for some circuit $I_i$ \underline{or} $J$ 
      is a T-part of some $I_i$ of level $\ell_{I_i}\ge 1$.    
\end{enumerate}
\end{theorem}

The proofs of these theorems 
require some preparation. We turn our attention first 
to the proof of Theorem~\ref{T:complement-is-independent}. 
We will need the following lemma.

\begin{lemma}\mlabel{T:all-vectors-matter} 
Let $I$ be a T-flat, and let $J\subsetneq I$ be a nonempty  
independent subset. Then \ $r_{I\smsm J} +r_J \ge r_I + 1$. 
\end{lemma}

\begin{proof} 
We induce on $n$, the level of $I$. In the case $n=0$ 
the set $I$ is a circuit, hence the set $I\smsm J$ is 
independent. Therefore 
$r_{I\smsm J}+r_J=|I\smsm J| + |J| = |I| \ge r_I + 1$.  

Assume next that $n\ge 1$ and that the lemma is 
true for T-flats of level $n-1$. Since 
$r_{I\smsm J}+r_J\ge r_I$ always, we need to show that 
we cannot have equality in this formula. Assume that is 
not the case, i.e., assume we have $r_{I\smsm J}+ r_J=r_I$.  
Let $I=I_1\sqcup\dots\sqcup I_t$ be the T-partition of $I$.  
Then for some $1\le k\le t$ the independent set 
$J'=J\cap I_k$ is not empty, and let $1\le j\le t$ be an 
integer different from $k$. Since $J$ is independent and 
$J'$ is a subset of $J$ we obtain 
$r_{I\smsm J'}+ r_{J'}=r_I$ 
by Lemma~\ref{T:zero-intersection-ranks}, and therefore 
$r_{(I\smsm I_j)\smsm J'}+ r_{J'} = r_{I\smsm I_j}$ 
by Lemma~\ref{T:zero-intersection-subsets}. This 
contradicts our induction hypothesis, because 
$I'=I\smsm I_j$ is a T-flat of level $n-1$ and $J'$ 
is a nonempty independent subset of $I'$. 
\end{proof}

\begin{proof}
[Proof of Theorem~\ref{T:complement-is-independent}]  
Assume $J$ is independent. Since $I$ has level $n$ and 
$I\smsm J$ has level $n-1$, we have 
$
r_J + r_{I\smsm J} = 
|J| + |I\smsm J| - |I| + r_I + 1 = 
r_I + 1
$. 

Next, assume that $J$ is not independent, or  
equivalently that $|J|-r_J \ge 1$. Clearly $n\ge 1$, 
because otherwise $I$ is a circuit, and $J$ has to be 
idependent as a proper subset of $I$. Since 
$|I|-r_I = n+1$ and $|I\smsm J|-r_{I\smsm J}=n$, 
we obtain that $|J|- r_I + r_{I\smsm J}=1$. Thus 
\[
1\le |J|-r_J=1+ r_I - r_{I\smsm J} - r_J \le 1, 
\]
hence $|J|-r_J=1$ and $r_J + r_{I\smsm J}= r_I$. 
We need to show that $J$ is a circuit. Assume not, 
and let $J'$ be a circuit properly contained in $J$. 
Then $|J'|-r_{J'}=1=|J|-r_J$, hence 
Lemma~\ref{T:kernels-inclusion}
yields $0\le |J\smsm J'|-r_{J\smsm J'}\le 1$. If 
$|J\smsm J'|-r_{J\smsm J'}=1$, then by 
Lemma~\ref{T:kernels-intersection} we obtain 
$0=|J'\cap(J\smsm J')|-r_{J'\cap(J\smsm J')} = 1$, a 
contradiction. Therefore $|J\smsm J'|-r_{J\smsm J'}=0$, 
yielding that $J''=J\smsm J'$ is independent. 
Furthermore we have 
$r_{J'}+r_{J''}=r_{J'}+|J|-|J'|=r_J$. 
Since $r_{J'}+r_{I\smsm J}=r_{I\smsm J''}$ by 
Lemma~\ref{T:zero-intersection-subsets}, we obtain 
$r_{J''}+r_{I\smsm J''}=r_{J''}+r_{J'}+r_{I\smsm J}=
r_J+r_{I\smsm J}= r_I$, which contradicts 
Lemma~\ref{T:all-vectors-matter}.   

Finally, let $J'$ be a proper subset of $J$. Since 
$J$ is either independent or a circuit the set $J'$ 
is independent. If $J$ is a circuit then 
$r_J+r_{I\smsm J}=r_I$, hence 
$r_{J'} + r_{I\smsm J} = r_{J'\cup(I\smsm J)}$ by 
Lemma~\ref{T:zero-intersection-subsets}. If $J$ is 
independent we get 
\[
r_{J'\cup(I\smsm J)}\le r_{J'} + r_{I\smsm J} = 
r_I + 1 - r_{J\smsm J'} \le r_{I\smsm (J\smsm J')} = 
r_{J'\cup (I\smsm J)}, 
\] 
where the last inequality follows from 
Lemma~\ref{T:all-vectors-matter}. 
This completes the proof of the theorem. 
\end{proof}

Next, we begin work towards the proof of 
Theorem~\ref{T:components-unite} and its 
corollary. We will need the following two lemmas.

\begin{lemma}\mlabel{T:circuit-parts} 
Let $I$ be a T-flat,  
and let $J$ be a circuit properly 
contained in $I$ such that 
$r_{I\smsm J} + r_J=r_I$. 
Then $\ell_I\ge 1$, and $J$ is a T-part  of $I$. 
\end{lemma}

\begin{proof} 
It is clear that $\ell_I\ge 1$ since $I$ contains 
properly a circuit. Let $I=I_1 \sqcup \dots \sqcup I_t$ 
be the T-partition of $I$ and assume that $J$ is not a  
T-part of $I$. Then for some $k$ the set 
$J_k=I_k\cap J$ is a proper nonempty subset of both 
$I_k$ and $J$. But then the set $J'=J\smsm J_k$ is a 
nonempty independent subset of the level $n-1$ 
T-flat $I'=I\smsm I_k$, and  
$r_{I'\smsm J'} + r_{J'}=r_{I'}$ by 
Lemma~\ref{T:zero-intersection-subsets}. This 
contradicts Lemma~\ref{T:all-vectors-matter}. 
\end{proof}

\begin{lemma}\mlabel{T:parts-are-maximal} 
Let $I$ be a T-flat, and let $J$ be an 
independent subset of $I$ such that 
$r_{I\smsm J} + r_J = r_I + 1$. Then $J$ is contained 
in a T-part  of the T-flat $I$. 
\end{lemma}

\begin{proof} 
The statement is obvious when $I$ is a circuit, thus 
we assume that $\ell_I\ge 1$. 
Let $I=I_1\sqcup\dots\sqcup I_t$ be the T-partition 
of $I$. Suppose  $J$ is not contained in a T-part of $I$. 
Then for each $k$ the set $J_k=J\cap I_k$ is either empty, 
or a proper subset of $J$. Since $J$ is not empty, $J_k$ 
is not empty for some $k$. Then $J'=J\smsm J_k$ is a 
nonempty independent subset of the T-flat $I'=I\smsm I_k$, 
hence $r_{I'\smsm J'} + r_{J'}\ge r_{I'} + 1$ by 
Lemma~\ref{T:all-vectors-matter}, and therefore  
$r_{I'\smsm J'} + r_{J'} = r_{I'} + 1$ by 
Lemma~\ref{T:rank-one-intersection-subsets}. But then 
Lemma~\ref{T:rank-one-intersection-splitting} yields 
$r_{I\smsm J_k} + r_{J_k}= r_I$,  which contradicts 
Lemma~\ref{T:all-vectors-matter}. 
\end{proof}

\begin{proof}
[Proof of Theorem~\ref{T:components-unite}]  
It suffices to prove the theorem in the case when the 
level of $I$ is $n+1$. The first assertion is then clear, 
since $J$ is obtained by the removal of a single T-part 
of $I$. To prove the second assertion, assume $n\ge 1$ 
and let $J'$ be a T-part of $I$ different from 
$I'=I\smsm J$. It suffices to prove that $J'$ is 
contained in a T-part of $J$. If $J'$ is 
a circuit, then we have $r_{J'}+ r_{I\smsm J'}=r_I$ 
by Theorem~\ref{T:complement-is-independent}(2), hence 
from Lemma~\ref{T:zero-intersection-subsets} we 
obtain $r_{J'}+ r_{J\smsm J'}=r_J$. Therefore $J'$ is   
a T-part of the T-flat $I'$ by 
Lemma~\ref{T:circuit-parts}. If $J'$ is independent, then 
from Theorem~\ref{T:complement-is-independent}(1) we have 
$r_{J'} + r_{I\smsm J'}= r_I + 1$, hence  
$r_{J'} + r_{J\smsm J'}= r_J + 1$ by 
Lemma~\ref{T:all-vectors-matter} and 
Lemma~\ref{T:rank-one-intersection-subsets}. 
Therefore $J'$ is contained 
in a T-part of $J$ by Lemma~\ref{T:parts-are-maximal}. 
\end{proof}

\begin{proof}
[Proof of Corollary~\ref{T:contained-in-circuit}] 
We induce on the level $n$ of $I$. Let $J''$ be a T-part 
of $I$ different from $J'$. 
When $n=0$ then $I$ itself is the desired circuit. 
When $n=1$ 
the set $C=I\smsm J''$ is a circuit containing $J$.  
Assume $n\ge 2$ and that the result is true for $n-1$. 
Then the set 
$I'=I\smsm J''$ is a T-flat of level $n-1$, and contains 
$J'$. Therefore $J'$ is contained 
in a T-part of $I'$ by Theorem~\ref{T:components-unite}, 
hence in a circuit $C$ inside $I'\subset I$ by the 
induction hypothesis.
\end{proof}

Finally, we proceed with the proof of 
Theorem~\ref{T:product-is-admissible}. The next 
lemma is a needed ingredient.

\begin{lemma}\mlabel{T:splitting}
Let $I$ and $J$ be disjoint sets such that 
$r_I + r_J= r_{I\cup J}$. 
If $I'$ is a T-flat of $\bfM$ inside $I\cup J$ such that 
$I'\cap I\ne \varnothing$ then the set 
$I'\cap I$ is a T-flat of \ $\bfM$. 
\end{lemma}

\begin{proof} 
Let $J'$ be a T-part of $I'$ that contains elements of $J$. 
By induction on the size of $I'$ it is enough to show that 
$J'$ is contained in $J$. Since $I'$ is a T-flat, by 
Corollary~\ref{T:contained-in-circuit} there is 
a circuit $C$ in $I'$ that contains $J'$. 
Therefore by Lemma~\ref{T:zero-intersection-circuits}
the set $C$, hence also $J'$, is inside the set $J$. 
\end{proof}

\begin{proof}
[Proof of Theorem~\ref{T:product-is-admissible}]  
(1) We have 
$
\ell_I=|I|-r_I - 1= 
\sum_{i=1}^k |I_i|- \sum_{i=1}^k r_{I_i} - 1 = 
\sum_{i=1}^k (\ell_{I_i} + 1) - 1 = 
\ell_{I_1} + \dots + \ell_{I_k} + k - 1
$. 

(2) Assume first that $k=2$. In that case  
if $I=I_1\oplus I_2$ is a T-flat, then $I_1$ and $I_2$ 
are T-flats  by Lemma~\ref{T:splitting}. Conversely, if 
both $I_1$ and $I_2$ are T-flats then so is 
$I=I_1\oplus I_2$ by Remark~\ref{T:T-remarks}(c). 
The first assertion of part (2) now follows from the case 
$k=2$ by an elementary induction. 

Next we consider the second assertion. 
Since $k\ge 2$ and by the first assertion  
each set $I_i$ is a T-flat, we have that $n=\ell_I\ge 1$. 
Furthermore, 
if $I_i$ is a circuit then it is a T-part of $I$ by 
Lemma~\ref{T:circuit-parts}. 
Since $I$ is the disjoint union of its T-parts, 
to complete the proof it suffices 
by Remark~\ref{T:T-remarks}(b) to show that when $J$ is a 
T-part of some non-circuit $I_j$ 
then the set $I\smsm J$ is a T-flat of level 
$n-1$. Let $I'=\bigoplus_{i\ne j} I_i$. 
Then $I=I'\oplus I_j$ 
and by Lemma~\ref{T:zero-intersection-subsets} we get 
$I\smsm J=I'\oplus (I_j\smsm J)$. 
Therefore $I\smsm J$ is a T-flat by the first assertion 
and we have by part (1)  
\[
\ell_{I\smsm J}  = \ell_{I'} + \ell_{I_j\smsm J} + 1  
                 = \ell_{I'}+(\ell_{I_j} - 1) + 1 
                 = \ell_I - 1  
                 = n - 1, 
\]
which is the desired equality. 
\end{proof}

\section{Connected T-flats and multiplicity spaces}

In this section $\bfM$ is a matroid on a finite set $S$. 
The following theorems describe the structure of 
the connected T-flats of $\bfM$. 

\begin{theorem}\mlabel{T:components-are-independent}
Let \ $I$ be a connected set of level $n\ge 0$.  
\begin{enumerate} 
\item The set \ $I$ is a T-flat, and every T-part of \ $I$ 
      is independent. 

\item The T-parts of \ $I$ are 
      all the maximal elements of the partially ordered set  
\[
\{ J\mid \text{$J$ is an independent subset of \ $I$ with \ 
               $r_{I\smsm J} + r_J = r_I + 1$} \},  
\]
      where the partial ordering is by inclusion. 
\end{enumerate} 
\end{theorem}

\begin{proof} 
(1) Let $I'$ be the union of all circuits contained in $I$. 
Since $n\ge 0$ the set $I'$ is nonempty hence 
a T-flat, and is clearly the unique maximal T-flat 
contained in $I$. Furthermore, the set $I\smsm I'$ contains no 
circuits, hence is independent, and also no circuit in $I$ 
intersects $I\smsm I'$. It follows that $I=(I\smsm I')\oplus I'$, 
and by the connectedness of $I$ we get $I=I'$, therefore $I$ is 
a T-flat. Let $J$ be a T-part of $I$. 
If $J$ is not independent, then by 
Theorem~\ref{T:complement-is-independent} it is a circuit  and  
$r_J + r_{I\smsm J}=r_I$. Thus $I=J\oplus (I\smsm J)$ which 
contradicts the connectedness of $I$. 

(2) This is immediate from part (1), 
Theorem~\ref{T:complement-is-independent},  
and Lemma~\ref{T:parts-are-maximal}. 
\end{proof} 

As direct corollary of the proof above we have  

\begin{corollary}\mlabel{T:largest-T-flat-decomposition} 
Let $A$ be a subset of $S$, and let $I$ be the union of all 
circuits contained in $A$ (in particular $I$ is empty when 
$A$ is independent, and is the unique maximal T-flat contained 
inside $A$ otherwise). Let $J=A\smsm I$. 

The set $J$ is independent, and  $A=I\oplus J$. \qed   
\end{corollary}

The next result plays a key role in the characterization of 
connected sets via their multiplicity spaces. It 
is due to Tutte and is a straightforward consequence of   
\cite[p.~148, (3.3)]{Tu58}. For 
completeness, we provide a (different) proof.

\begin{theorem}\mlabel{T:connected-contains-connected}
{\rm(\cite{Tu58})}. 
Let $I$ be a connected T-flat of level $n\ge 1$. 
Then $I$ contains a connected T-flat of level $n-1$.   
\end{theorem}

\begin{proof} 
We induce on $n$. When $n=1$ the result is immediate 
because circuits are connected, so we assume that 
$n\ge 2$ and that the result holds for connected T-flats of 
smaller levels. Pick a T-part $J$ of $I$ such that $I\smsm J$ 
has a connected component of lowest possible level. Thus 
by Proposition~\ref{T:connected-components-decomposition} 
and Theorem~\ref{T:product-is-admissible} 
the T-flat $I\smsm J$ decomposes as  
\[
I\smsm J= J_1\oplus\dots\oplus J_p
\]
for some connected T-flats $J_i$, where $J_1$ is  
the connected component of lowest possible level. 
Note that if $p=1$ then we are done, thus for the rest of this 
proof we assume that $p\ge 2$. 
 
We select a subset $C_1$ of $J_1$ as follows: if 
$J_1$ is circuit we take $C_1=J_1$, otherwise we choose 
$C_1$ to be a T-part of $J_1$ such that $J_1\smsm C_1$ is 
a connected T-flat (a T-part like that exists by our induction 
hypothesis). 
Thus by Theorem~\ref{T:product-is-admissible} the set $C_1$ is 
a T-part of $I\smsm J$, while  
by Theorem~\ref{T:components-unite} we have that  
$C_1$ is a union of T-parts of $I$. 
Let $J'$ be one of these T-parts. 
Then $I\smsm J'$ is a T-flat of level $n-1$ and contains the 
level $n-2$ T-flat  
$(I\smsm J)\smsm C_1$. Therefore $C_1'=(C_1\smsm J')\cup J$ 
is a T-part of the T-flat $I\smsm J'$. 

We claim that $I\smsm J'$ is a connected T-flat. 
Suppose not, and let  
\[
I\smsm J'= I'_1\oplus\dots\oplus I'_t 
\]
with $t\ge 2$, 
each $I'_j$ a connected T-flat, and $C'_1$ 
contained in $I'_1$. Since by 
Lemma~\ref{T:zero-intersection-subsets} and 
Proposition~\ref{T:direct-sum-criterion} we have 
\[
(I'_1\smsm C'_1)\oplus I'_2\oplus\dots\oplus I'_t = 
(I\smsm J')\smsm C'_1 = (I\smsm J)\smsm C_1    = 
(J_1\smsm C_1)\oplus J_2\oplus\dots\oplus J_p,  
\]
it follows from Lemma~\ref{T:zero-intersection-circuits}, 
Theorem~\ref{T:connected-components-decomposition}, 
and the connectedness of $I_2'$ that $I_2'$ equals 
one of the connected components   
$J_1\smsm C_1$ or $J_i$ for some $i\ge 2$. 
In particular, either $J_1\smsm C_1$ or $J_i$ for some 
$i\ge 2$ is a connected component of $I\smsm J'$. 
However, $J_1\smsm C_1$ cannot be a connected component   
of $I\smsm J'$ because of the minimality property of $J$. 
Thus without loss of 
generality we may assume that $J_2$ is a connected   
component of $I\smsm J'$. Then  
$r_{J_2}+ r_{(I\smsm J')\smsm J_2}= r_{I\smsm J'}$, 
and  $r_{J_2} + r_{I\smsm J_2}\ge r_I + 1$ 
because $I$ is connected. 
However, $J_1$ is a T-flat containing the independent 
by Theorem~\ref{T:components-are-independent}
T-part $J'$ of $I$, and disjoint from $J_2$. 
Thus by Theorem~\ref{T:complement-is-independent} and 
Lemma~\ref{T:rank-one-intersection-subsets} we have 
$r_{J'} + r_{J_1\smsm J'}= r_{J_1} + 1$. Therefore 
$r_I-r_{I\smsm J'}= r_{J_1} - r_{J_1\smsm J'}$ and we have 
inclusions $J'\subseteq J_1\subseteq I\smsm J_2\subseteq I$. 
It follows that 
$r_{J_1} - r_{J_1\smsm J'} = r_{I\smsm J_2} - 
                             r_{(I\smsm J_2)\smsm J'}$ 
by Lemma~\ref{T:rank-one-intersection-tower}. This yields 
$
r_I - r_{J_1} + r_{J_1\smsm J'} - 
r_{(I\smsm J')\smsm J_2} = 
r_I - r_{I\smsm J_2}
$. 
From here we get 
$
r_I-r_{J'}+1-r_{(I\smsm J')\smsm J_2}=
r_I-r_{I\smsm J_2}
$,  
hence 
$
r_{I\smsm J'} - r_{(I\smsm J')\smsm J_2} = 
r_I - r_{I\smsm J_2}
$. 
Therefore we obtain 
\[
0=r_{J_2}+ r_{(I\smsm J')\smsm J_2}- r_{I\smsm J'} 
 =r_{J_2}+ r_{I\smsm J_2} - r_I \ge 1, 
\]
yielding the desired contradiction. 
Thus the T-flat $I\smsm J'$ is connected, hence  
the connected T-flat $I$ contains a connected 
T-flat of level $n-1$. 
\end{proof}

\begin{proposition}\mlabel{T:T-part-is-circuit} 
Let $A$ be a T-flat of level $n$ which is not connected, 
and let $I$ be a connected T-flat of level $n-1$ inside 
$A$. Then the T-part $A\smsm I$ is a circuit. 
\end{proposition}

\begin{proof} 
By Lemma~\ref{T:zero-intersection-circuits}, 
Proposition~\ref{T:connected-components-decomposition}, 
Theorem~\ref{T:product-is-admissible}, 
and the fact that $A$ is not connected, 
the connected T-flat $I$ is a connected component 
of $A$, we have $A=(A\smsm I)\oplus I$, and $A\smsm I$ 
is a T-flat of level $0$ hence a citcuit. 
\end{proof}

Finally, let $\phi\:U_S \lra W$ be a 
representation of the matroid $\bfM$ over a field $\Bbbk$. 
We have the following important characterization 
of connected T-flats in terms of their multiplicity 
spaces. 

\begin{theorem}\mlabel{T:multiplicity-zero} 
We have $S_I(\phi)\ne 0$ if and only if 
$I$ is a connected T-flat. 
\end{theorem}

\begin{proof} 
To prove the ``only if'' part of the theorem we 
assume that $I$ is not a connected T-flat and we will 
show that $S_I(\phi)=0$. This is clear from the definition 
if $I$ is not a T-flat, so we assume that $I$ is a T-flat 
but is not connected. Let $n$ be the level of $I$. 
Since the T-flat $I$ is not connected, 
we have $n\ge 1$. We prove our assertion by induction 
on $n$. When $n=1$ the non-connected T-flat $I$ must by 
Theorem~\ref{T:product-is-admissible} 
and Theorem~\ref{T:connected-components-decomposition} 
have a decomposition into connected components 
of the form $I=I_1 \oplus I_2$ with $I_1$ and $I_2$ 
circuits, in particular $I=I_1\sqcup I_2$ is the 
T-partition of $I$. Thus by 
Remark~\ref{T:rank-1-intersection-a} we have 
$V_{I_1}\cap V_{I_2}=0$, hence from 
Definition~\ref{D:multiplicity-space} we get  
$S_I = V_{I_1}\cap V_{I_2}=0$. 
Therefore we assume $n\ge 2$ and that our assertion is 
true for T-flats of level $n-1$.  In view of the 
surjectivity of the map \altref{E:multiplication-map}, 
it is enough to show that for each T-flat $J$ of 
level $n-1$ in $I$ we have either 
$V_{J}\cap V_{I\smsm J}=0$ or $S_J=0$. 
Since $I$ is not connected, 
by Theorem~\ref{T:product-is-admissible} and 
Theorem~\ref{T:connected-components-decomposition} 
the T-part $I\smsm J$ is either a circuit, or is 
a T-part of one of the connected components 
$A$ of $I$ with $\ell_A\ge 1$, and we have a non-trivial 
direct sum decomposition $I=A\oplus (I\smsm A)$. 
If $I\smsm J$ is a circuit, then by 
Theorem~\ref{T:complement-is-independent} and 
Remark~\ref{T:rank-1-intersection-a} we have 
$V_{I\smsm J}\cap V_J=0$. If $I\smsm J$ is not a circuit, 
then by Lemma~\ref{T:zero-intersection-subsets}  
we get a non-trivial direct sum decomposition 
$J=A\smsm (I\smsm J) \oplus I\smsm A$,  
yielding that the T-flat $J$ is not connected.   
Therefore by our induction hypothesis $S_J=0$, 
which completes the proof of our assertion, hence of 
the ``only if'' part of the theorem. 

To prove the ``if'' part of the theorem 
we again induce on $n$, the level of $I$. The result is 
trivial for $n=0$, thus we assume that $n\ge 1$ and 
that the assertion is true for connected T-flats of level 
$n-1$. By Theorem~\ref{T:components-are-independent}
for each T-flat $J$ of level $n-1$ in $I$ 
the space $V_J\cap V_{I\smsm J}$ is $1$-dimensional.  
In view of the injectivity of the map $\nu$, 
it is enough to show that for some T-flat $J$ of level 
$n-1$ in $I$ the multiplicity space $S_J$ is not zero. 
By Theorem~\ref{T:connected-contains-connected} the 
connected T-flat $I$  
contains a connected T-flat $J$ of level $n-1$.  
Therefore $S_J\ne 0$ by our induction hypothesis. 
\end{proof}

\section{Multiplicity spaces for minors} 

In this section $\phi\: U_S \lra W$ is a representation 
over $\Bbbk$ of a matroid $\bfM$ on a finite set $S$.  
We study the relationship between the multiplicity 
spaces of $\bfM$ and the multiplicity spaces of the 
minors of $\bfM$. We also present 
(modulo Theorem~\ref{T:the-T-bullets-are-complexes}) 
the proofs of Theorem~\ref{T:exactness-combinatorial} 
and Theorem~\ref{T:T-bullet-of-sum}. 

We begin with the (relatively simpler) 
behavior under the operation restriction. 

\begin{theorem}
\mlabel{T:multiplicity-spaces-for-restriction}
Let \ $Y$ be a subset of \ $S$, and let $I$ be a subset of 
level $n$ in  $Y$. 
\begin{enumerate} 
\item $I\in\mathcal T_n(\bfM|Y)$ if and only if \ 
      $I\in\mathcal T_n(\bfM)$. 
\item $S_I(\phi)=S_I(\phi|Y)$ and \ $T_I(\phi)=T_I(\phi|Y)$.   
\item If \ $I\in\mathcal T_n(\bfM)$ and 
      $J\subseteq I$ is an element of \ 
      $\mathcal T_{n-1}(\bfM)$ then \ 
      $\phi_n^{IJ}=(\phi|Y)_n^{IJ}$.   
\end{enumerate} 
\end{theorem}

\begin{proof} 
When $n=-1$ part (1) is immediate from the definitions. 
When $n\ge 0$ part (1) is immediate by  
Theorem~\ref{T:T-flats-of-minors}. It follows from that 
same theorem that a chain $\mathbb I$ is in $C_{\bfM}(I)$ 
if and only if it is a chain in $C_{\bfM|Y}(I)$. Since for 
any subset $A\subseteq Y$ the spaces $U_A$ and $V_A$ are 
the same for $\phi$ and for $\phi|Y$, part (2) follows 
from the definitions of multiplicity space and T-space.  
Part (3) is now also clear, since all spaces and maps 
appearing in the definitions of $\phi_n^{IJ}$ and 
$(\phi|Y)_n^{IJ}$ are the same.
\end{proof}

We are now ready to give the proofs of 
Theorem~\ref{T:exactness-combinatorial} and 
Theorem~\ref{T:T-bullet-of-sum}.

\begin{proof}
[Proof of Theorem~\ref{T:exactness-combinatorial}] 
It is immediate from 
Theorem~\ref{T:multiplicity-spaces-for-restriction}(2) 
that for each $n\ge 0$ we have 
$T_n(\phi|Y)=T_n(\phi)\big|_Y$ and 
$T_n^+(\phi|Y)=T_n^+(\phi)\big|_Y$. Similarly, from 
Theorem~\ref{T:multiplicity-spaces-for-restriction}(3) 
we obtain for each  $n\ge 1$ that 
$\phi_n\big|_{T_n|_Y}=(\phi|Y)_n$ and 
$\phi_n^+\big|_{T_n^+|_Y}=(\phi|Y)_n^+$. 
The desired conclusion 
now follows in view of 
Theorem~\ref{T:the-T-bullets-are-complexes}. 
\end{proof}

\begin{proof}
[Proof of Theorem~\ref{T:T-bullet-of-sum}] 
The inclusion $\supseteq$ is clear from 
Theorem~\ref{T:exactness-combinatorial}. If  
a T-flat $A$ of $\bfM$ is not contained in some $S_i$, 
then we have by Lemma~\ref{T:zero-intersection-subsets} 
a nontrivial decomposition 
$A=A\cap S_1 \oplus \dots \oplus A\cap S_k$, and 
therefore $T_A(\phi)=0$ by 
Theorem~\ref{T:multiplicity-zero}. 
The desired equality of complexes is now immediate.  
\end{proof}

Next, we turn our attention to the significantly more 
complicated behavior of multiplicity spaces under the 
operation contraction.

\begin{theorem}
\mlabel{T:multiplicity-spaces-for-contraction} 
Let \ $Y$ be a subset of $S$ with $S\smsm Y$ 
independent in $\bfM$, let $A\subseteq Y$ be a 
T-flat of \ $\bfM.Y$ of level $n$, let 
$B=S\smsm (Y\smsm A)^{\bsfC_{\bfM^*}}$, let 
$\ol W = W/V_{S\smsm Y}$, and let 
$\pi\: W \lra \ol W$ be the 
canonical projection map.   

Then the canonical surjection of symmetric powers 
$\pi_n\: S_n W \lra S_n\ol W$ 
induces by restriction 
a surjection of multiplicity spaces \ 
$
\pi_{Y,A}^{\phi}=\pi_{Y,A} \: 
S_B(\phi) \lra S_A(\phi.Y)
$.   
\end{theorem}

\begin{proof} 
Note that by Remark~\ref{T:maximal-T-flat}(b) the set 
$B$ is a T-flat in $\bfM$ of level $n$. 
Let $\mathbb I$ be a chain 
$I^{(0)}\subsetneq\dots\subsetneq I^{(n)}$ in 
$C_{\bfM}(B)$ and for each $k$ let $J^{(k)}=I^{(k)}\cap Y$. 
Since $S\smsm Y$ is independent, each set $J^{(k)}$ is 
nonempty. Let $C_{\bfM}^Y(B)$ be the collection of all 
chains in $C_{\bfM}(B)$ such that $J^{(k)}$ and $J^{(k-1)}$ 
are distinct for every $k\ge 1$. Note that if $\mathbb I$ is 
not in $C_{\bfM}^Y(B)$ then for some $k\ge 1$ we have 
$J^{(k)}=J^{(k-1)}$, hence  
$I^{(k)}\smsm I^{(k-1)}\subseteq S\smsm Y$; therefore 
$V(\mathbb I)\subseteq V_{S\smsm Y}(\phi)S^{n-1}W$, yielding 
that $\pi_n\bigl(V(\mathbb I)\bigr)=0$. 
On the other hand  if $\mathbb I$ is in $C_{\bfM}^Y(B)$    
then by Theorem~\ref{T:T-flats-of-minors}  and 
Remark~\ref{T:maximal-T-flat}(b) each set 
$J^{(k)}$ is a T-flat in $\bfM.Y$ of level $k$; thus 
intersection with $Y$ produces a chain $\mathbb I\cap Y$ 
in $C_{\bfM . Y}(A)$, and,  
by Theorem~\ref{T:T-partition-of-contraction},  
every chain in $C_{\bfM . Y}(A)$ can be obtained in this 
way from a chain in $C_{\bfM}^Y(B)$. 

Next, assume $\mathbb I$ is in $C_{\bfM}^Y(B)$ and 
consider the set $I_k=I^{(k)}\smsm I^{(k-1)}$. It 
is a T-part of $I^{(k)}$, and the nonempty set 
$J_k=I_k\cap Y$ is a T-part of $J^{(k)}$ in $\bfM . Y$. 
Recall that $J_k$ is independent in $\bfM . Y$ if and only 
if $J_k\cup L$ is independent in $\bfM$ for every 
independent in $\bfM$ subset $L$ of $S\smsm Y$. 
Therefore if $J_k$ is independent in $\bfM . Y$ we have  
$V_{J_k}(\phi)\cap V_{S\smsm Y}(\phi) = 0$, and also 
\[
V_{I_k}(\phi)\cap V_{S\smsm Y}(\phi) = 
V_{I_k\smsm J_k}(\phi).
\]
Furthermore, for each $I$ the space $V_I(\phi)$ maps 
under $\pi$ onto $V_{I\cap Y}(\phi . Y)$. This induces 
for each $k\ge 1$ a map  
\[
\pi^{(k)}\: 
V_{I^{(k-1)}}(\phi)\cap V_{I_k}(\phi) \lra  
V_{J^{(k-1)}}(\phi . Y)\cap V_{J_k}(\phi . Y).
\] 
We claim that $\pi^{(k)}$ is surjective. Indeed, 
if $J_k$ is not independent in $\bfM . Y$ then 
$V_{J^{(k-1)}}(\phi . Y)\cap V_{J_k}(\phi . Y) = 0$, 
hence $\pi^{(k)}$ is surjective in that case. Assume that 
$J_k$ is independent in $\bfM . Y$. Then 
$V_{J^{(k-1)}}(\phi . Y)\cap V_{J_k}(\phi . Y)$ is a 
$1$-dimensional vector space, and $I_k$ is independent 
in $\bfM$. Since the kernel of $\pi^{(k)}$ is the space 
\[
\ker\pi^{(k)}= 
V_{I^{(k-1)}}(\phi)\cap V_{I_k}(\phi)
                   \cap V_{S\smsm Y}(\phi)
\]  
and $V_{I^{(k-1)}}(\phi)\cap V_{I_k}(\phi)$ is 
$1$-dimensional, to prove surjectivity it is enough 
to show that 
$
V_{I^{(k-1)}}(\phi)\cap V_{I_k}(\phi)
                   \cap V_{S\smsm Y}(\phi)=0
$. 
Since $I'=I_k\smsm J_k$ is a strict subset of the 
independent T-part $I_k$, we have by 
Theorem~\ref{T:complement-is-independent}  that 
$
r^{\bfM}_{I'\cup I^{(k-1)}} = 
r^{\bfM}_{I'}+r^{\bfM}_{I^{(k-1)}}
$. 
Therefore  
\[
V_{I^{(k-1)}}(\phi)\cap V_{I_k}(\phi)
                   \cap V_{S\smsm Y}(\phi) = 
V_{I^{(k-1)}}(\phi)\cap V_{I'}(\phi) = 0, 
\]
completing the proof of our claim.  

It follows that for each chain $\mathbb I$ in 
$C_{\bfM}^Y(B)$ the space $V(\mathbb I)$ is mapped 
under $\pi_n$ surjectively 
onto the space $V(\mathbb I\cap Y)$. Since every 
chain in $C_{\bfM . Y}(A)$ comes from a chain in 
$C_{\bfM}^Y(B)$, and $\pi_n\bigl(V(\mathbb I)\bigr)=0$ 
when $\mathbb I$ is not in $C_{\bfM}^Y(B)$, it is 
clear from the definitions that 
$\pi_n(S_B(\phi))=S_A(\phi . Y)$.  
\end{proof}

\begin{remark}\mlabel{T:contraction-composition} 
Let $Y\subseteq Y_1$ be subsets of $S$ such that 
$S\smsm Y$ is independent, 
let $A$ be a T-flat of $\bfM.Y$, 
let $B=S\smsm (Y\smsm A)^{\bsfC_{\bfM^*}}$, 
and let $A_1= Y_1\cap B$. 
It is straightforward to verify 
that $\phi.Y = (\phi.Y_1).Y$, that 
$\pi_Y^{\phi}=\pi_Y^{\phi.Y_1}\circ \pi_{Y_1}^{\phi}$, 
and therefore we have 
\[
\pi_{Y,A}^{\phi} = 
\pi_{Y,A}^{\phi.Y_1}\circ \pi_{Y_1,A_1}^{\phi}.
\] 
In particular, to study the properties of the maps 
$\pi_{Y,*}^{\phi}$ one may contract one element at 
a time.  We examine this situation in more detail 
in the following two theorems. 
\end{remark}

\begin{theorem}\mlabel{T:multiplicity-kernels} 
Let $a\in S$ with \ $\{a\}$ independent in $\bfM$, let 
$A\subseteq S_a$ be a T-flat of $\bfM . S_a$, and let 
$B=S\smsm (S_a\smsm A)^{\bsfC_{\bfM^*}}$. Then: 
\begin{enumerate} 
\item  If the set \ $\{a\}$ is not a T-part of $B$ in 
       $\bfM$ then the surjective homomorphism   
       $\pi_{S_a,A} \: S_B(\phi) \lra S_A(\phi .S_a)$ 
       is an isomorphism. 

\item  If the set \ $\{a\}$ is a T-part of $B$ then 
       we have a canonical complex 
\[
\begin{CD}
0 \lra S_A(\phi)\otimes V_a(\phi) 
\overset{\nu}\lra S_B(\phi)
@> \pi_{S_a,A} >> 
S_A(\phi .S_a) \lra 0  
\end{CD} 
\] 
of vector spaces over the field $\Bbbk$. 
\end{enumerate} 
\end{theorem}

\begin{proof} 
(1) It suffices to show that $S_B(\phi)$ and 
$S_A(\phi.S_a)$ have the same dimension. Note that  
$(\bfM.S_a)|A = (\bfM|A\cup a).A$, and similarly 
$(\phi.S_a)|A = (\phi|A\cup a).A$. Since 
by Corollary~\ref{T:T-flats-on-S-a} we have 
$B\subseteq A\cup a$, in view of 
Theorem~\ref{T:T-flats-of-minors} and 
Theorem~\ref{T:multiplicity-spaces-for-restriction} 
we may assume that $S=A\cup a$ (hence also $S_a=A$).  

Suppose now that 
$a\notin B$. Then by Corollary~\ref{T:T-flats-on-S-a} 
we have $B=A$, and $r^{\bfM}_S=r_B^{\bfM}+r^{\bfM}_a$. 
Therefore $\bfM=\bfM|B + \bfM|a$, and it is 
straightforward to verify that in such a 
situation we have  
$W= V_B(\phi)\oplus V_a(\phi) \oplus W'$ for some 
subspace $W'$. It follows that the kernel of the map 
$\pi_n\: S_nW \lra S_n\ol W$ does not intersect 
$S_nV_B(\phi)$, hence the induced map 
$\pi_{S_a,A}\: S_B(\phi) \lra S_A(\phi.A)$ 
is injective, therefore an isomorphism by 
Theorem~\ref{T:multiplicity-spaces-for-restriction} 
as desired.     

Consider now the remaining case $B=A\cup a$. 
If $B$ is a circuit in $\bfM$ then so is $A$ in 
$\bfM.A$, and the desired equality of dimensions is 
immediate from the definitions. Thus we assume in the 
sequel that $B$ has level $\ell_B^{\bfM}\ge 1$.   
Let $B=I_1\sqcup\dots\sqcup I_k$ be the T-partition 
of $B$, and let $I=I_k$ be the T-part of $B$ that 
contains (properly) $a$. For each $i$ let 
$B_i=B\smsm I_i$ and let 
$V_i=V_{B_i}(\phi)\cap V_{I_i}(\phi)$. Thus each 
$V_i$ is either $0$ or $1$-dimensional, and it follows   
from Remark~\ref{T:span-of-T-parts}(d) that    
the space $S_B(\phi)$ is contained  
inside the subalgebra of the symmetric algebra $S(W)$ 
generated by the subspace 
$\widetilde W=V_1 + \dots + V_k\subset W$.   
Note that the nonzero vector $\e_a=\phi(a)$ is not in 
$\widetilde W$. Indeed, 
otherwise we would have $\e_a=v_1+\dots+v_k$ with  
each $v_i\in V_i$, 
hence $\e_a-c_kv_k\in V_k$; thus $\e_a\in V_k$ and 
therefore $0\ne \e_a\in V_{B_k}(\phi)\cap V_a(\phi)$  
which contradicts 
Theorem~\ref{T:complement-is-independent}(3). 
It follows that 
$(\e_aS_{n-1}W) \cap S_B(\phi) = 0$, hence the map 
$\pi_{S_a,A}$ is injective as desired.

(2) We have by Corollary~\ref{T:T-flats-on-S-a} that 
$A=B\smsm\{a\}$ is a T-flat in $\bfM$. 
Also, since $\{a\}$ is independent, we have that 
$V_A(\phi) \cap V_a(\phi)$ is $1$-dimensional, therefore 
equals $V_a(\phi)$. Finally, since the kernel of 
$\pi_{S_a,A}$ is the space 
$(\e_a S_{n-1}W) \cap S_B(\phi)$, the 
desired conclusion is immediate from the definitions. 
\end{proof}

We close this section with the statement of a key 
result that strengthens the conclusion of part (2) of 
Theorem~\ref{T:multiplicity-kernels}. We postpone the 
proof till Section~15.

\begin{theorem}\mlabel{T:multiplicity-exactness} 
Let $a\in S$ with \ $\{a\}$ independent in $\bfM$, let 
$A\subseteq S_a$ be a T-flat of $\bfM . S_a$, and let 
$B=S\smsm (S_a\smsm A)^{\bsfC_{\bfM^*}}$.  
If \ $\{a\}$ is a T-part of $B$ then 
\[
\begin{CD}
0 \lra S_A(\phi)\otimes V_a(\phi) \overset{\nu}\lra 
       S_B(\phi)
@> \pi_{S_a,A} >> 
S_A(\phi .S_a) \lra 0  
\end{CD} 
\] 
is an exact sequence of vector spaces. 
\end{theorem}

\section{$T_\bullet(\phi)$ and 
         $T_\bullet(\phi)^+$ are complexes} 

In this section $\phi\: U_S \lra W$ is a 
representation over $\Bbbk$ of a matroid $\bfM$ on a 
finite set $S$. Our  goal is to give the proof of 
Theorem~\ref{T:the-T-bullets-are-complexes}. The 
main ingredient in the proof is the following result.

\begin{theorem}\mlabel{T:composition-zero-high-degree}
Let $I$ be a T-flat of level $n\ge 1$ and let 
$A\subset I$ be an element of $\mc T_{n-2}(\bfM)$. Let 
$J_1, \dots , J_k$ be the T-flats of level $n-1$ inside  
$I$ that contain $A$. Then 
\[
\sum_{j=1}^k (-1)^{|J_j|}
\phi_{n-1}^{J_jA} \circ \phi_n^{I J_j} = 0. 
\]
\end{theorem}

The proof of this theorem relies on explicit 
computations. 
In order for us to do these computations
effectively, we will need to fix a 
linear ordering with certain properties 
on the elements of the T-flat $I$. 
The exact properties required from this ordering 
will be  specified later, but once an ordering has 
been fixed, we will use the 
following notation. 

\begin{definition} 
Let $I$ be a set of level $n$ with a fixed linear 
oredring on its elements. 

(a) If $I$ is independent we set $\widehat I=I$. 

(b) If $I$ is a T-flat, we introduce the following 
objects. Let $I^{(0)}=I$, and let $i^{(0)}$ be the 
smallest element of $I^{(0)}$. Proceeding inductively 
on $r$, given we have already defined $I^{(r)}$ and 
$i^{(r)}$ for some $0\le r < n$,  
we write $I^{(r+1)}$ for the only T-flat of level 
$n-r-1$ inside $I^{(r)}$ that does not contain 
$i^{(r)}$ and set $i^{(r+1)}$ for its smallest 
element. Note that if $C$ is a circuit inside 
$I^{(r)}$ then $C$ is a union of T-parts 
of $I^{(r)}$, therefore if $C$ does not contain 
$i^{(r)}$ it must be a circuit inside $I^{(r+1)}$. 
Finally we set 
\[
\widehat{I} = I\smsm \{i^{(0)},\dots, i^{(n)} \}. 
\]

(c) If $I$ is neither independent nor a T-flat, 
then by 
Corollary~\ref{T:largest-T-flat-decomposition}
we have a nontrivial decomposition $I=K\oplus J$ 
where $K$ is the unique maximal T-flat contained 
in $I$, and $J$ is independent. In that case we set 
$\widehat I=\widehat K\cup \widehat J$. 
\end{definition}

\begin{remark}\mlabel{T:I-hat-is-independent}
It is clear from this construction that 
$\widehat I$ does not contain any circuits, and 
since $r_I = |\widehat I|$ the set $\widehat I$ 
is a maximal independent set in $I$. 
\end{remark}

\begin{proof}
[Proof of Theorem~\ref{T:composition-zero-high-degree}]  
We assume that both $I$ and $A$ are connected, since 
otherwise the result is immediate from 
Theorem~\ref{T:multiplicity-zero}. 
Let $I_j=I\smsm J_j$, thus the sets $I_1, \dots, I_k$ 
are the T-parts of $I$ that are disjoint from $A$. 
Since $A\in\mc T_{n-2}$, 
it is clear that $k\ge 2$. When $n\ge 2$ we set 
$J=A$, and when $n=1$ we set $J=I_{k+1}$ where 
$I_{k+1}$ is the 
T-part of $I$ that contains the singleton $A=\{a\}$.   

We choose a linear order on the set $I$ such that 
the elements of the set $J$ are the biggest $|J|$ 
elements of $I$. Then it 
is clear that $\hat J\subset \hat I$. 
In addition we require that all the elements of 
$I_i=I\smsm J_i$ are smaller than the elements of 
$I_j$ for $i< j$, and, when $n=1$, 
that $a$ is the smallest element of $J$. 
Let $i_p$ be the smallest element of $I_p$ and let  
$I_p'=I_p\smsm \{i_p\}$, in particular when $n=1$ 
we have $a=i_{k+1}$.  Thus 
$i^{(0)} = i_1 = \min\{ i_p \mid 1\le p\le k \}$, 
we set $I' = I\smsm\{i^{(0)}\}$, and, when $n=1$ 
we set $J'=J\smsm a$.  
Also we have $i^{(1)}=i_2$, and clearly 
\[
\widehat I = 
I_1'\sqcup I_2'\sqcup I_3\sqcup\dots
    \sqcup I_k\sqcup \widehat J.  
\]
Similarly, for $J_i$ we have 
$\widehat J\subset\widehat{J_i}$, and 
\[
\begin{aligned} 
\widehat{J_1} &= 
I_2'\sqcup I_3\sqcup\dots\sqcup I_k
    \sqcup \widehat J,   
\\ 
\widehat{J_2} &= 
I_1'\sqcup I_3\sqcup\dots\sqcup I_k
    \sqcup \widehat J, \quad\text{and} 
\\ 
\widehat{J_i} &= 
I_1'\sqcup I_2\sqcup\dots\sqcup I_{i-1}
    \sqcup I_{i+1}\sqcup\dots\sqcup I_k
    \sqcup \widehat J  
    \quad\text{for }i\ge 3. 
\end{aligned}   
\]

Using the linear ordering on $I$  
we identify each subset of $I$ with the 
increasing sequence of 
its elements. For a sequence 
$K=(a_1,\dots,a_q)$ we set  
\[
e_K = 
e_{a_1}\wedge \dots 
       \wedge e_{a_q} \in \wedge^q U
\quad\text{and}\quad  
v_K = 
\phi(e_{a_1})\wedge\dots
             \wedge\phi(e_{a_q})\in
             \wedge^q V.
\] 
If $K$ is a subset of $I$ and 
the elements $\phi(e_{a_1}), \dots, \phi(e_{a_q})$ 
form a basis of $V_K$, we write $e^*_K$ and $v^*_K$ 
for the unique elements of $\wedge^q U_K^*$ and 
$\wedge^q V^*_K$ respectively, such that 
$e^*_K(e_K)=1$ and $v^*_K(v_K)=1$. 
In particular $v_{\widehat I}$ is a basis 
of $\wedge^{r_I} V_I$ and $v^*_{\widehat I}$ 
is the dual basis of $\wedge^{r_I} V^*_I$.
With this notation we also have 
\[
\begin{aligned}
v_{\widehat{J_1}} &= 
v_{I_2'}\wedge v_{I_3}\wedge \dots \wedge v_{I_k} 
                      \wedge v_{\widehat J}  \\ 
v_{\widehat{J_2}} &= 
v_{I_1'}\wedge v_{I_3}\wedge \dots \wedge v_{I_k} 
                      \wedge v_{\widehat J}   \\ 
v_{\widehat{J_j}} &= 
v_{I_1'}\wedge v_{I_2}\wedge\dots\wedge v_{I_{j-1}}  
        \wedge v_{I_{j+1}}\wedge\dots\wedge v_{I_k}
        \wedge v_{\widehat J} 
        \quad\text{ for }\quad 3\le j\le k 
\end{aligned}
\]
Next, for $p=1,\dots, k$ 
(or $p=1,\dots,k+1$ when $n=1$) let 
\[
u_p = \sum_{i\in I_p} x_i e_i 
\]
be the unique element of $U_{I_p}$ 
such that $x_{i_p}=1$ and 
$w_p= \phi(u_p)\in V_{I_p}\cap V_{I\smsm I_p}$. 
Since $I$ is connected for each $j$ we have that 
the T-part $I_j$ is independent, and  
\[
v_{I_j}= w_j\wedge v_{I_j'}.
\] 
Also, for each $j\ge 2$ let 
\[
v_j=u_1 + \sum_{2\le i\ne j}\alpha_{ij}u_i
\] 
be the unique element in $U_{J_j\smsm J}$ such that 
$y_j=\phi(v_j)\in V_J\cap V_{J_j\smsm J}$ 
(we know that $v_j$ has 
to be a unique linear combination of 
the $u_i$ because of 
Remarks~\ref{T:rank-1-intersection-b}, 
which also imply that 
we have $\alpha_{ij}\ne 0$ for all $i\ne j$). 
When $j=1$ we define 
analogously $v_1$, and $y_1$, 
and again we have a unique linear 
combination with non zero coefficients   
\[
v_1=u_2 + \sum_{3\le i}\alpha_{i1}u_i.
\] 
It follows that 
$v_{\widehat I}=v_{I_1'}\wedge v_{\widehat{J_1}}$ 
and 
$
v_{\widehat I} = 
(-1)^{|I_1'||I_2'|}v_{I_2'}\wedge v_{\widehat{J_2}}
$, 
while for $j\ge 3$ we have   
\[
\begin{aligned}
v_{\widehat I} &= 
v_{I_1'}\wedge v_{I_2'}\wedge 
v_{I_3}\wedge\dots\wedge v_{I_k}
       \wedge v_{\widehat J}  
\\[+7pt] 
&= 
v_{I_1'}\wedge v_{I_2'}\wedge 
v_{I_3}\wedge\dots\wedge v_{I_{j-1}}\wedge 
(w_j\wedge v_{I_j'})\wedge v_{I_{j+1}} 
    \wedge\dots \wedge v_{I_k}
    \wedge v_{\widehat J}  
\\[+5pt] 
&= 
v_{I_1'}\wedge v_{I_2'}\wedge 
v_{I_3}\wedge\dots \wedge 
\left( \textstyle{ 
       \frac{\phi(v_1) - w_2 - 
             \sum_{3\le i\ne j}\alpha_{i1}w_i}
            {\alpha_{j1}} }  
       \wedge v_{I_j'}  
\right)\wedge\dots \wedge 
v_{I_k}\wedge v_{\widehat J}  
\\[+1pt] 
&= 
v_{I_1'}\wedge v_{I_2'}\wedge 
v_{I_3}\wedge\dots\wedge v_{I_{j-1}}\wedge 
\left( \frac{-w_2}{\alpha_{j1}} \wedge v_{I_j'} 
\right)\wedge v_{I_{j+1}} 
       \wedge\dots \wedge v_{I_k}
       \wedge v_{\widehat J}  
\\[+1pt]  
&= 
\left( 
\textstyle{
\frac{(-1)^{|I_2|+\dots +|I_{j-1}|}}
     {\alpha_{j1}}} 
\right)
v_{I_1'}\wedge (w_2\wedge v_{I_2'})\wedge 
v_{I_3} \wedge\dots \wedge 
v_{I_j'}\wedge\dots \wedge 
v_{I_k} \wedge v_{\widehat J}  
\\[+3pt] 
&= 
\left( 
\textstyle{ 
\frac{(-1)^{|I_2|+\dots +|I_{j-1}|}}
     {\alpha_{j1}} }
\right)
v_{I_1'}\wedge v_{I_2} \wedge 
v_{I_3}\wedge\dots\wedge 
v_{I_{j-1}}\wedge v_{I_j'}\wedge 
v_{I_{j+1}}\wedge\dots \wedge 
v_{I_k}    \wedge v_{\widehat J}  
\\[+1pt] 
&= 
\left( 
\textstyle{
\frac{ (-1)^{|I_j|
       (\sum_{i=2}^{j-1}|I_i|)+|I_j'||I_1'|}}
     { \alpha_{j1} } } 
\right) 
v_{I_j'}   \wedge v_{I_1'}\wedge 
v_{I_2}    \wedge\dots\wedge 
v_{I_{j-1}}\wedge v_{I_{j+1}}\wedge\dots\wedge 
v_{I_k}    \wedge v_{\widehat J}.   
\end{aligned} 
\]
Therefore for $3\le j\le k$ we obtain the formula 
\begin{equation}
\elabel{E:computing-v-widehat-I-for-j-ge-3}
v_{\widehat I} = 
(-1)^{\tau_j}\frac{1}{ \alpha_{j1} }  \ 
v_{I_j'}\wedge v_{\widehat{J_j}}.    
\end{equation}
where 
$\tau_j=1+ |I_1| + |I_j|(1+\sum_{i=1}^{j-1}|I_i|)$. 

In order to prove the Theorem it is enough to show that  
\[
\sum_{j=1}^k (-1)^{|J_j|}
(\phi_n^{I,J_j})^* \circ (\phi_{n-1}^{J_j,A})^* = 0,  
\]
and for the rest of this argument we will concentrate 
on proving this equality. Also, despite their 
similarity there are enough distinctions that 
unfortunately force us 
to treat the cases $n\ge 2$ and $n=1$ separately.

\emph{Case $n\ge 2$}. In this case $A$ is a 
connected T-flat. Furthermore, if $J_j$ is not 
connected, then (since we are assuming 
that $J=A$ is connected) 
Proposition~\ref{T:T-part-is-circuit}
implies that $J_j\smsm J$ must 
be a circuit, and we must have $y_j=0$. 
When $J_j$ is  connected we get 
\begin{equation}\elabel{E:computing-v-sub-J}
v_{J_j\smsm J} = 
\begin{cases} 
y_1\wedge v_{I_2'}
   \wedge v_{I_3}\wedge\dots\wedge v_{I_k} 
   &\text{if } j=1; 
\\
y_2\wedge v_{I_1'}
   \wedge v_{I_3}\wedge\dots\wedge v_{I_k} 
   &\text{if } j=2; 
\\ 
y_j\wedge v_{I_1'}
   \wedge v_{I_2}\wedge\dots\wedge v_{I_{j-1}} 
   \wedge v_{I_{j+1}}\wedge\dots\wedge v_{I_k} 
   &\text{if } j\ge 3. 
\end{cases}
\end{equation} 

Let us take an arbitrary element  
$z$ in 
$
T_J(\phi)^* = 
S_J(\phi)\otimes \wedge^{|J|}U_J^* 
         \otimes \wedge^{r_J}V_J
$. 
Thus $z$ has the form 
$z=x\otimes e^*_J\otimes v_{\widehat J}$ 
for some $x\in S_J(\phi)$. 
Computing the image $t_j$ of $z$ 
in $T_{J_j}(\phi)^*$ we obtain 
\[
t_j = (\phi_{n-1}^{J_j,J})^*(z) = 
y_j x\otimes e^*_{J_j}
     \otimes v_{\widehat J_j}. 
\]  
Indeed, this is trivially true when $J_j$ is 
not connected since then $y_j=0$. When $J_j$ 
is connected, this follows from 
\altref{E:computing-v-sub-J} in view of the 
definition of $y_j$, the definition of 
$\phi_{n-1}^{J_jJ}$, and the fact that 
the map 
$
\twedge V_{J_j\smsm J}^*\lra 
\twedge U_{J_j\smsm J}^* 
$ 
induced by $\phi^*$ sends $v_{J_j\smsm J}^*$ 
to $e_{J_j\smsm J}^*$. Similarly, computing the 
image of $t_j$ in $T_I(\phi)^*$ for $j\ge 3$ we 
obtain 
\[
(\phi_n^{IJ_j})^*(t_j)= 
(-1)^{\delta_j}\alpha_{j1} 
w_jy_jx\otimes e_I^*\otimes v_{\widehat I} 
\]
where 
$\delta_j = \tau_j + |I_j|\sum_{i=1}^{j-1}|I_i|$, 
and therefore 
$(-1)^{\delta_j}=(-1)^{1+|I_1|+|I_j|}$. For the 
remaining possibilities $j=1$ and $j=2$ we obtain 
\[
(\phi_n^{IJ_j})^*(t_j) = 
\begin{cases} 
w_1y_1\otimes e_I^*
      \otimes v_{\widehat I} 
      &\text{ if } j=1; 
\\ 
(-1)^{1+|I_1|+|I_2|} 
w_2y_2\otimes e_I^*
      \otimes v_{\widehat I} 
      &\text{ if } j=2. 
\end{cases}
\]
In particular, this yields  
\[
\begin{aligned}
&\sum_{j=1}^k (-1)^{|J_j|}(\phi_n^{IJ_j})^*(t_j) 
\\ 
&= 
\left( 
(-1)^{|J_1|}w_1y_1x - (-1)^{|I_1|+|I|}  w_2y_2x - 
\sum_{j=3}^k (-1)^{|I_1|+|I|}\alpha_{j1}w_jy_jx 
\right) 
\otimes e_I^* \otimes v_{\widehat I} 
\\ 
&= 
(-1)^{|I|+|I_1|}
\left( 
w_1y_1x - w_2y_2x - 
\sum_{j=3}^k \alpha_{j1} w_jy_jx 
\right)
\otimes e_I^* \otimes v_{\widehat I}.   
\end{aligned}
\]
Thus to complete our proof of the case $n\ge 2$, 
it suffices to show that the expression   
\begin{equation}\elabel{E:key-equality} 
Z = 
w_1y_1 - w_2y_2 - 
\sum_{j=3}^k \alpha_{j1}w_jy_j  
\end{equation}
is equal to zero 
in the symmetric algebra of the space $W$. 
Since 
$
w_1y_1=w_2w_1 + \sum_{i=3}^k \alpha_{i1}w_iw_1
$ 
and 
$
w_2y_2=w_1w_2 + \sum_{i=3}^k \alpha_{i2}w_iw_2
$, 
while for $j\ge 3$ we have 
$
w_jy_j = 
w_1w_j + \sum_{2\le i\ne j}\alpha_{ij}w_iw_j
$, 
we obtain
\[
Z = - 
\sum_{j=3}^k 
(\alpha_{j2} + \alpha_{j1}\alpha_{2j})w_jw_2 - 
\sum_{j=3}^{k-1} 
\sum_{i=j+1}^k 
(\alpha_{j1}\alpha_{ij}+\alpha_{i1}\alpha_{ji})
w_iw_j.  
\]  
Let $U'$ be the subspace of $U_{I\smsm J}$ 
spanned by the (independent) vectors 
$u_1,\dots,u_k$. Note that by 
Remark~\ref{T:rank-1-intersection-b}(d)  
a vector in $U_{I\smsm J}$ 
gets mapped under $\phi$ inside $V_J$ only if 
that vector belongs to the subspace $U'$, in 
particular $K'=\ker\phi\cap U_{I\smsm J}$ is a 
subspace of $U'$. Let $V'=\phi(U')$. Thus 
$V_{I\smsm J}\cap V_J = V'\cap V_J$, and 
(since the level of $J$ is $n-2$) we get 
\[
\begin{aligned}
\dim V'\cap V_J  
&= 
\dim V_{I\smsm J} +\dim V_J - \dim V_I 
\\ 
&= 
r_{I\smsm J} + (|J|-(n-2)-1) - (|I| -n -1)
\\  
&= 
r_{I\smsm J} - |I\smsm J| -n +2 -1 +n + 1 
\\ 
&= 
r_{I\smsm J} - |I\smsm J| + 2 
\\ 
&= 
2 - \dim K' 
\end{aligned}
\]              
Therefore the kernel of the induced by $\phi$ 
map $U'\lra W/V_J$ is exactly $2$-dimensional. 
On the other hand, each of the vectors 
$v_1,\dots, v_k$ belong to that kernel. 
It follows that the matrix 
\[
\left(
\begin{matrix}
0          & 1 & 1          & \dots &  1         & \dots 
                                    &  1         & \dots 
                                    &  1        \\ 
1          & 0 & \alpha_{23}& \dots & \alpha_{2j}& \dots 
                                    & \alpha_{2i}& \dots
                                    &\alpha_{2k}\\ 
\alpha_{31}& \alpha_{32} & 0& \dots &            & \dots  
                                    &            & \dots
                                    &\alpha_{3k}\\ 
\vdots     & \vdots &\vdots &       & \vdots     &  
                                    & \vdots     &  
                                    & \vdots    \\ 
\alpha_{j1}& \alpha_{j2} &  & \dots & 0          & \dots 
                                    &\alpha_{ji} & \dots 
                                    &           \\ 
\vdots     &        &\vdots &       & \vdots     &        
                                    & \vdots     &       
                                    & \vdots    \\ 
\alpha_{i1}& \dots  &       & \dots & \alpha_{ij}& \dots 
                                    & 0          & \dots 
                                    &           \\ 
\vdots     &        &\vdots &       & \vdots     &        
                                    & \vdots     &        
                                    & \alpha_{k-1,k}    \\ 
\alpha_{k1}&\alpha_{k2}&    &\dots  &            & \dots 
                                    &            & \alpha_{k,k-1} 
                                    & 0         \\  
\end{matrix}
\right)
\] 
has rank at most two, hence all  minors of size $3$ are 
zero. Since the minor on rows $1,2,j$ and columns $1,2,j$ is 
precisely the coefficient of $w_jw_2$ in $Z$, and since the 
minor on rows $1,j,i$ and columns $1,j,i$ is precisely 
the coefficient of $w_iw_j$ in $Z$, the desired conclusion 
is immediate.

\emph{Case $n=1$}. In this case $J=I_{k+1}$ is the 
independent T-part of the connected T-flat $I$ of 
level $1$ that contains the singleton $A=\{a\}$. 
Furthermore, since each $J_j$ is a circuit, there 
is a unique 
$\alpha_j\in\Bbbk$ such that $y_j=\alpha_j w_{k+1}$.  
Since we have 
$v_J=w_{k+1}\wedge v_{J'}$, we obtain 
\begin{equation}\elabel{E:computing-v-sub-J-again}
\begin{aligned}
v_{J_j\smsm a} &= 
\begin{cases} 
y_1\wedge v_{I_2'}\wedge v_{I_3}
   \wedge\dots\wedge v_{I_k}  \wedge v_{J'} 
   &\text{if } j=1; 
\\
y_2\wedge v_{I_1'}\wedge v_{I_3}
   \wedge\dots\wedge v_{I_k}  \wedge v_{J'} 
   &\text{if } j=2; 
\\ 
y_j\wedge v_{I_1'}\wedge v_{I_2}
   \wedge\dots\wedge v_{I_{j-1}} 
   \wedge v_{I_{j+1}}
   \wedge\dots\wedge v_{I_k}  \wedge v_{J'} 
   &\text{if } 3\le j  
\end{cases}
\\ 
&= (-1)^{|J_j|-|J|-1}\alpha_j v_{\widehat{J_j}}
\end{aligned}
\end{equation} 

In order to complete the proof of the theorem it 
is enough to show that  
\[
\sum_{j=1}^k (-1)^{|J_j|}
[(\phi_1^{I,J_j})^* \circ (\phi_0^{J_j,a})^*]
(e^*_a) = 0.   
\]
Computing the image $t_j$ of $e_a^*$ in 
$
T_{J_j}(\phi)^* = 
\twedge U_{J_j}^* 
\otimes \twedge V_{J_j}
$ 
we obtain 
\[
t_j = 
(\phi_0^{J_j,a})^*(e_a^*) = 
- \alpha_j 
(e^*_{J_j}\otimes v_{\widehat J_j}). 
\]  
Indeed, this follows from 
\altref{E:computing-v-sub-J-again} in view of 
the definition of $\phi_0^{J_j,a}$, and the fact 
that the isomorphism induced by $\phi^*$  
\[
\twedge V_{J_j\smsm a}^*\otimes U_a^* \lra 
\twedge U_{J_j\smsm a}^*\otimes U_a^* 
\overset{\tau}\lra 
U_a^*\otimes \twedge U_{J_j\smsm a}^* \lra 
\twedge U_{J_j}^*
\] 
sends $v_{J_j\smsm a}^*\otimes e_a^*$ to 
the element $e_a^*\wedge\psi\in\twedge U_{J_j}^*$ 
where $\psi\in\wedge^{|J_j|-1}U_{J_j}^*$ is an 
element such that $\psi(e_{J_j\smsm a})=1$ 
(which implies  
\[
[e_a^*\wedge\psi](e_{J_j})=  
(-1)^{|J_j|-|J|}[e_a^*\wedge\psi]
                (e_a\wedge e_{J_j\smsm a})= 
(-1)^{|J_j|-|J|}    
\] 
and therefore 
$e_a^*\wedge\psi = (-1)^{|J_j|-|J|}e_{J_j}^*$). 

Computing the image of 
$t_j$ in $T_I(\phi)^*$ for $j\ge 3$ we obtain 
as in the Case~$n\ge 2$ that 
\[
(\phi_1^{IJ_j})^*(t_j)= 
(-1)^{\delta_j}\alpha_{j1}\alpha_j 
w_j\otimes e_I^*\otimes v_{\widehat I} 
\]
where 
$\delta_j = \tau_j + 1+|I_j|\sum_{i=1}^{j-1}|I_i|$, 
and therefore 
$(-1)^{\delta_j}=(-1)^{|I_1|+|I_j|}$. For the 
remaining possibilities $j=1$ and $j=2$ we obtain 
\[
(\phi_1^{IJ_j})^*(t_j) = 
\begin{cases} 
-\alpha_1 
w_1\otimes e_I^*\otimes v_{\widehat I} 
&\text{ if } j=1; 
\\ 
(-1)^{|I_1|+|I_2|} 
\alpha_2 w_2\otimes e_I^*\otimes v_{\widehat I} 
&\text{ if } j=2. 
\end{cases}
\]
In particular, this yields  
\[
\begin{aligned}
& \sum_{j=1}^k (-1)^{|J_j|}
(\phi_1^{IJ_j})^*(t_j) 
\\ 
& = 
\left( 
- (-1)^{|J_1|}\alpha_1w_1 + 
  (-1)^{|I_1|+|I|}\alpha_2 w_2
+ \sum_{j=3}^k (-1)^{|I_1|+|I|}
\alpha_{j1}\alpha_j w_j
\right) 
\otimes e_I^* \otimes v_{\widehat I} 
\\ 
&= 
(-1)^{|I|+|I_1|+1}
\left(
\alpha_1 w_1 - \alpha_2 w_2 - 
\sum_{j=3}^k \alpha_{j1}\alpha_j w_j
\right)
\otimes e_I^* \otimes v_{\widehat I}.   
\end{aligned}
\]
Thus to complete our proof, it suffices to show 
that the vector    
\begin{equation}\elabel{E:key-equality-again} 
Z = 
\alpha_1 w_1 - \alpha_2 w_2 - 
\sum_{j=3}^k \alpha_{j1}\alpha_j w_j  
\end{equation}
is equal to zero in the space $W$.  

Let $U'$ be the subspace of $U_I$ spanned by 
the (independent) vectors 
$u_1,\dots,u_k,u$. Note that by 
Remark~\ref{T:rank-1-intersection-b}(d)  
a vector in $U_I$ 
gets mapped to $0$ under $\phi$ only if that 
vector belongs to the subspace $U'$, in 
particular  
$K'=\ker\phi\cap U_I$ is a subspace of $U'$. 
Since $I$ has level $1$, we have 
$\dim K' = |I|-r_I = 2$. 
Furthermore, the vectors 
$v_1-\alpha_1 u, \dots, v_k-\alpha_k u$ belong 
to $K'$, and the first two are clearly 
independent hence form a basis 
for $K'$. It follows that the matrix 
\[
\left(
\begin{matrix}
0          & 1 & 1          & \dots &  1         & \dots 
                                    &  1        \\ 
1          & 0 & \alpha_{23}& \dots & \alpha_{2j}& \dots 
                                    &\alpha_{2k}\\ 
\alpha_{31}& \alpha_{32} & 0& \dots & \alpha_{3j}& \dots  
                                    &\alpha_{3k}\\ 
\vdots     & \vdots &\vdots &       & \vdots     &  
                                    & \vdots    \\ 
\alpha_{j1}&\alpha_{j2}&\alpha_{j3}& \dots & 0   & \dots 
                                    &\alpha_{jk}\\ 
\vdots     & \vdots &\vdots &       & \vdots     &        
                                    & \vdots    \\ 
\alpha_{k1}&\alpha_{k2}&\alpha_{k3} & \dots & \alpha_{kj}
                                                 & \dots 
                                    &  0        \\ 
-\alpha_1  &-\alpha_2&-\alpha_3 &\dots & -\alpha_j 
                                                 & \dots 
                                       & - \alpha_k    
\end{matrix}
\right)
\] 
has rank exactly two, hence all  minors of size $3$ are 
zero. In particular, the minor on rows $1,j,k+1$ and 
columns $1,2,j$ is zero, which yields that for each 
$j\ge 3$ we have 
\[
\alpha_{j1}\alpha_j = 
\alpha_{j1}\alpha_2 - \alpha_{j2}\alpha_1. 
\]
Therefore we get 
\[
Z = \alpha_1\left( -\alpha_2 w + 
    w_1 + \sum_{j=3}^k \alpha_{j2}w_j \right) -  
    \alpha_2\left( -\alpha_1 w + 
    w_2 + \sum_{j=3}^k \alpha_{j1}w_j \right). 
\]
The desired conclusion is now immediate since 
$
-\alpha_2 w + w_1 + 
\sum_{j=3}^k \alpha_{j2}w_j=
\phi(v_2-\alpha_2 u)=0
$ 
and similarly 
$
-\alpha_1 w + w_2 + 
\sum_{j=3}^k \alpha_{j1}w_j=
\phi(v_1-\alpha_1 u)=0
$. 
\end{proof}

Now we can give the proof of 
Theorem~\ref{T:the-T-bullets-are-complexes}, 
in this way  
completing also the proofs of 
Theorem~\ref{T:exactness-combinatorial} 
and Theorem~\ref{T:T-bullet-of-sum}.  

\begin{proof}
[Proof of Theorem~\ref{T:the-T-bullets-are-complexes}]  
Assume first that $n\ge 1$. 
We need to show that if $z\in T_I$ for some connected 
T-flat $I$ of level $n$ then 
$[\phi_{n-1}\circ\phi_n](z)=0$. 
Thus it is enough to show that the component 
of  $[\phi_{n-1}\circ\phi_n](z)$ in $T_A$ is zero for 
every $A\in\mc T_{n-2}(\bfM)$. By the definitions of 
$\phi_n$ and $\phi_{n-1}$ it is clear that this 
component is zero when $A$ is not a subset of $I$; and 
when $A$ is a subset of $I$ this component is precisely 
\[
(-1)^{|A|}\sum_{i=1}^k 
(-1)^{|J_i|} [\phi_{n-1}^{J_i,A}\circ
                     \phi_n^{I,J_i}](z) 
\]
where we use the same notation as in  
Theorem~\ref{T:composition-zero-high-degree}, 
hence is zero by that theorem.

It remains to show that $\phi_{-1}\circ\phi_0=0$. 
However, when $I$ is a circuit the map $\phi|I$ 
represents a uniform matroid and therefore the 
sequence $T_\bullet(\phi|I)^{+}$ is a complex, 
see Example~\ref{E:uniform-matroid-2}. Since by 
Theorem~\ref{T:multiplicity-spaces-for-restriction} 
the sequence $T_\bullet(\phi|I)^{+}$ is the same as 
the sequence 
\[
\begin{CD}
0 \lra T_I(\phi) 
@> \phi_0 >> 
U_I 
@> \phi_{-1} >> 
W \lra 0, 
\end{CD}
\]
we obtain that 
$(\phi_{-1}\circ\phi_0)\big|_{T_I} = 0$.   
It is now immediate that 
$\phi_{-1}\circ\phi_0 = 0$.  
\end{proof}

\section{The definition of  
$
(\pi.Y)_{\bullet}^{\phi}\: 
T_\bullet(\phi.Y)\lra T_\bullet(\phi)
$}
 
In this section $\phi\: U_S \lra W$ is a 
representation over 
$\Bbbk$ of a matroid $\bfM$ on a finite set $S$, 
the set $Y$ is a subset of $S$ 
such that $S\smsm Y$ is independent, and 
$
\pi_Y^{\phi} = \pi\: 
W \lra \ol W = W/V_{S\smsm Y}(\phi)
$ 
is the canonical projection. 

Our main goal is to present the definition and 
some basic properties of the morphism of 
complexes 
\[
(\pi.Y)_{\bullet}^{\phi} \: 
T_\bullet(\phi.Y) \lra T_\bullet(\phi)
\]
mentioned in the statement of 
Theorem~\ref{T:contraction-embedding}. 
We will use these in Section~15  
for the proofs of both 
Theorem~\ref{T:contraction-embedding} and 
Theorem~\ref{T:exactness-of-complex}. We begin 
by introducing some notation.

\begin{definition}\mlabel{D:extra-notation} 
Let $A$ be a T-flat of $\bfM.Y$, and let 
$B=S\smsm (Y\smsm A)^{\bsfC_{\bfM^*}}$. 

(a) By repeatedly using 
Corollary~\ref{T:T-flats-on-S-a} 
and contracting one element at a time, 
it is clear that one always has an equality 
$r_{B\smsm A}^{\bfM} + r_A^{\bfM.Y}=r_B^{\bfM}$ 
and an exact sequence 
\begin{equation}\elabel{E:contraction-sequence} 
0 \lra V_{B\smsm A}(\phi) \lra V_B(\phi) 
  \lra V_A(\phi.Y) \lra 0. 
\end{equation}  

(b) Since $S\smsm Y$ is independent,  
$r_{B\smsm A}^{\bfM}=|B\smsm A|$ 
and we have canonical isomorphisms  
\[
{\bf c}_{Y,A}^{\phi}=
{\bf c}_Y^{\phi} = {\bf c}\: \ 
\twedge V_A(\phi.Y)^*  
\otimes 
\twedge V_{B\smsm A}(\phi)^* 
\ \lra \ 
\twedge V_B(\phi)^* 
\] 
and 
\[
{\bf d}_{Y,A}^{\phi} = 
{\bf d}_Y^{\phi} = {\bf d}\: \ 
\Bbbk \ \lra \ 
\twedge U_{B\smsm A}\otimes 
\twedge V_{B\smsm A}(\phi)^*
\]
induced by the canonical projection 
$\pi\: W \lra \ol{W}$ 
and by the map $\phi$, respectively. 
\end{definition} 

We are now ready to proceed with 
the definition of the canonical 
morphism of complexes 
$
(\pi.Y)^{\phi}_\bullet \: 
T_\bullet(\phi.Y)\lra T_\bullet(\phi)
$. 

\begin{definition}\mlabel{D:contraction-morphism}
(a) Let $A$ be a T-flat of $\bfM.Y$, and let 
$B=S\smsm (Y\smsm A)^{\bsfC_{\bfM^*}}$. 
We define a canonical inclusion homomorphism 
\[
(\pi.Y)^{\phi}_A \: T_A(\phi.Y) \lra T_B(\phi) 
\]
as the composition 
\begin{equation}\elabel{E:contraction-composition} 
\begin{CD}
T_A(\phi.Y) 
@> 1\otimes{\bf d} >> 
T_{AB}^{(1)} 
@> \tau >>  
T_{AB}^{(2)} 
@> \pi^*_{Y,A}\otimes\wedge\otimes{\bf c} >> 
T_B(\phi),   
\end{CD}
\end{equation}
where 
\begin{equation}\elabel{E:contraction-embedding}  
\begin{aligned} 
T_{AB}^{(1)} &= 
T_A(\phi.Y)    \otimes \twedge U_{B\smsm A}
               \otimes \twedge V_{B\smsm A}(\phi)^*, 
\\ 
&= 
S_A(\phi . Y)^*\otimes \twedge U_A 
               \otimes \twedge V_A(\phi.Y)^* 
               \otimes \twedge U_{B\smsm A}
               \otimes \twedge V_{B\smsm A}(\phi)^*, 
\\[+5pt]  
T_{AB}^{(2)} &= 
S_A(\phi . Y)^*\otimes \twedge U_A 
               \otimes \twedge U_{B\smsm A}
               \otimes \twedge V_A(\phi.Y)^* 
               \otimes \twedge V_{B\smsm A}(\phi)^*,   
\end{aligned}
\end{equation}
and the map $\tau$ is, as usual, the canonical 
isomorphism that simply rearranges the components 
of the tensor product in the indicated order.  

(b) For each $n\ge 0$ 
we define a canonical injective homomorphism 
\[
(\pi.Y)^{\phi}_n \: T_n(\phi.Y) \lra T_n(\phi) 
\]
by the requirement that it restricts to the 
component $T_A(\phi.Y)$ 
of $T_n(\phi.Y)$ as 
\[
(\pi.Y)^{\phi}_n\big|_{T_A(\phi.Y)} = 
(-1)^{n|B\smsm A|}(\pi.Y)^{\phi}_A 
\] 
for every T-flat $A$ in $\bfM.Y$ of level $n$. 

(c) We write $(\pi.Y)^{\phi}_\bullet$ for the 
sequence of homomorphisms 
$\{(\pi.Y)^{\phi}_n\}_{n\ge 0}$. When 
$\phi$ is clear from context we omit 
superscripts and 
write $(\pi.Y)_A$, $(\pi.Y)_n$, and 
$(\pi.Y)_\bullet$.   
\end{definition}

The main results in this section are the 
following two theorems. 

\begin{theorem}
\mlabel{T:contraction-embedding-composes} 
Let $Y$ be such that $S\smsm Y$ is independent 
in $\bfM$, and let $Y\subseteq Y_1$. Let $A$ be 
a T-flat in $\bfM.Y$, let 
$B=S\smsm (Y\smsm A)^{\bsfC_{\bfM^*}}$, and let 
$A_1=B\cap Y_1$. Then 
\[
(\pi.Y)^\phi_A = 
(\pi.Y_1)^{\phi}_{A_1} \circ (\pi.Y)^{\phi.Y_1}_A.  
\]
Furthermore, 
$
(\pi.Y)^{\phi}_\bullet = 
(\pi.Y_1)^{\phi}_\bullet \circ 
(\pi.Y)^{\phi.Y_1}_\bullet.  
$
\end{theorem}

\begin{theorem}
\mlabel{T:contraction-embedding-is-morphism} 
Let \ $Y$ be a subset of \ $S$ such that \ 
$S\smsm Y$ is independent in $\bfM$. 

Then $(\pi.Y)^{\phi}_\bullet$ is an injective 
morphism of complexes from 
$T_\bullet(\phi.Y)$ into $T_\bullet(\phi)$. 
\end{theorem}

The proof of 
Theorem~\ref{T:contraction-embedding-is-morphism} 
is involved, and we present it in the next section. 
We conclude this section with the proof of 
Theorem~\ref{T:contraction-embedding-composes}.

\begin{proof}
[Proof of Theorem~\ref{T:contraction-embedding-composes}] 
The second assertion of the theorem is an immediate 
consequence of the first, and we proceed with the 
proof of that first assertion.
  
Since $S\smsm Y$ 
is independent, we have that $B$ and $A_1$ are 
T-flats of level $n$ in $\bfM$ and $\bfM.Y_1$, 
respectively. For the same 
reason the sets $A_1\smsm A$ and $B\smsm A$ are 
independent in $\bfM.Y_1$ and $\bfM$ respectively. 
Therefore $B\smsm A_1$ is also independent in $\bfM$ 
and we have a canonical exact 
sequence of $\Bbbk$-vector spaces
\[
\begin{CD}
0\lra V_{B\smsm A_1}(\phi) 
@> \subseteq >> 
V_{B\smsm A}(\phi) 
@> \pi_{Y_1}^{\phi} >> 
V_{A_1\smsm A}(\phi.Y_1) \lra 0 
\end{CD}
\]
which induces an isomorphism 
\[
{\bf e}\: 
\twedge V_{A_1\smsm A}(\phi.Y_1)^* \otimes 
\twedge V_{B\smsm A_1}(\phi)^*   
\lra 
\twedge V_{B\smsm A}(\phi)^*. 
\] 
To make our commutative diagrams more compact, 
for the rest of this proof we use the following 
notation:  given  subsets 
$X, X_1, \dots, X_t$ we set 
$\ol{\ol V}_X=V_X(\phi.Y)$ and   
$\ol V_X=V_X(\phi.Y_1)$, as well as 
$
\twedge U_{X_1,\dots,X_t}= 
\twedge U_{X_1}\otimes\dots\otimes 
\twedge U_{X_t}
$. 
Using the associativity of the wedge product and 
the definitions of the maps ${\bf c}$, ${\bf d}$, 
and ${\bf e}$, it is now straightforward to 
verify that the diagrams 
\[
\begin{CD} 
\twedge \ol{\ol V}_A^* \otimes
\twedge \ol{V}_{A_1\smsm A}^* \otimes 
\twedge V_{B\smsm A_1}^* 
@> 1\otimes {\bf e} >>   
\twedge \ol{\ol V}_A^*\otimes
\twedge V_{B\smsm A}^*
\\ 
@V {\bf c}_Y^{\phi.Y_1}\otimes 1 VV 
@VV {\bf c}_Y^{\phi} V 
\\ 
\twedge \ol{V}_{A_1}^* \otimes 
\twedge V_{B\smsm A_1}^* 
@> {\bf c}_{Y_1}^{\phi} >> 
\twedge V_B^* 
\end{CD} 
\]
and 
\[ 
\begin{CD} 
\twedge U_{A,A_1\smsm A,B\smsm A_1}  
@> 1\otimes\wedge >>  
\twedge U_{A,B\smsm A} 
\\ 
@V \wedge\otimes 1 VV 
@VV \wedge V 
\\ 
\twedge U_{A_1,B\smsm A_1} 
@> \wedge >> 
\twedge U_B 
\end{CD} 
\]
together with 
\[ 
\begin{CD} 
\Bbbk 
@> \tau\ \circ\ 
   \left( {\bf d}_{Y,A}^{\phi.Y_1}\otimes
          {\bf d}_{Y_1,A_1}^{\phi}\right)  >>  
\twedge U_{A_1\smsm A,B\smsm A_1} \otimes 
\twedge \ol{V}_{A_1\smsm A}^* \otimes 
\twedge V_{B\smsm A_1}^*
\\ 
@V = VV 
@VV \wedge\otimes {\bf e} V 
\\ 
\Bbbk
@> {\bf d}_Y^\phi >> 
\twedge U_{B\smsm A} \otimes 
\twedge V_{B\smsm A}^* 
\end{CD} 
\]
are commutative. It follows that the diagram 
\[
\begin{CD} 
\twedge U_A\otimes \ol{\ol V}_A^* 
@> = >> 
\twedge U_A\otimes \ol{\ol V}_A^* 
@> = >> 
\twedge U_A\otimes \ol{\ol V}_A^* 
\\ 
@V \tau\circ(1\otimes{\bf d}) VV 
@V \tau\circ(1\otimes{\bf d}\otimes{\bf d}) VV 
@V \tau\circ(1\otimes{\bf d}) VV 
\\ 
\twedge U_{A,A_1\smsm A} \otimes 
\widetilde V_3 
@> \tau\circ(1\otimes{\bf d}) >> 
\twedge U_{A,A_1\smsm A,B\smsm A_1} \otimes 
\widetilde V 
@> 1\otimes\wedge\otimes 1\otimes {\bf e} >> 
\twedge U_{A,B\smsm A}\otimes 
\widetilde V_1
\\ 
@V \wedge\otimes{\bf c} VV 
@V \wedge\otimes 1\otimes{\bf c}\otimes 1 VV 
@V \wedge \otimes {\bf c} VV 
\\
\twedge U_{A_1} \otimes 
\twedge \ol{V}_{A_1}^* 
@> \tau\circ(1\otimes{\bf d}) >> 
\twedge U_{A_1,B\smsm A_1} \otimes 
\widetilde V_2
@> \wedge\otimes{\bf c} >> 
\twedge U_B \otimes 
\twedge V_B^* 
\end{CD} 
\]
is also commutative, where 
\[
\begin{aligned} 
\widetilde V  & = 
\twedge \ol{\ol V}_A^* \otimes  
\twedge \ol{V}_{A_1\smsm A}^* \otimes 
\twedge V_{B\smsm A_1}^* 
\\ 
\widetilde V_1 & = 
\twedge \ol{\ol V}_A^* \otimes 
\twedge V_{B\smsm A}^* 
\\ 
\widetilde V_2 & = 
\twedge \ol V_{A_1}^* \otimes 
\twedge V_{B\smsm A_1}^* 
\\ 
\widetilde V_3 & = 
\twedge \ol{\ol V}_A^*\otimes 
\twedge \ol V_{A_1\smsm A}^*.  
\end{aligned} 
\]
The desired conclusion now follows from  
Remark~\ref{T:contraction-composition} 
and Definition~\ref{D:contraction-morphism}. 
\end{proof}

\section{$(\pi.Y)^{\phi}_\bullet$ is a 
morphism of complexes}

In this section $\phi\: U_S \lra W$ is a 
representation over 
$\Bbbk$ of a matroid $\bfM$ on a finite set $S$. 
Our goal is to present 
the proof of 
Theorem~\ref{T:contraction-embedding-is-morphism}.  
The following result is a key ingredient.

\begin{theorem}
\mlabel{T:contraction-morphism-for-S-a} 
Let $a\in S$ be such that $\{a\}$ is independent 
in $\bfM$. Let $A$ be a T-flat in $\bfM.S_a$ of 
level $n\ge 1$, let $A'$ be a maximal proper 
T-flat of $\bfM.S_a$, let 
$B=S\smsm (S_a\smsm A)^{\bsfC_{{\bf M}^*}}$, 
and let 
$B'=S\smsm (S_a\smsm A')^{\bsfC_{{\bf M}^*}}$. 
Then  
\begin{equation}
\elabel{E:commutativity-for-pi-sub-A}
\begin{CD} 
T_A(\phi . S_a) 
@> (-1)^{|A'|}(\phi . S_a)_n^{AA'} >> 
T_{A'}(\phi . S_a)                                        
\\ 
@V (-1)^{n|B\smsm A|}(\pi.S_a)_A VV 
@VV (-1)^{(n-1)|B'\smsm A'|}(\pi.S_a)_{A'} V 
\\ 
T_B(\phi)     
@> (-1)^{|B'|}\phi_n^{BB'} >> 
T_{B'}(\phi)   
\end{CD}
\end{equation} 
is a commutative diagram.    
\end{theorem}

\begin{proof}  
Throughout this proof we use the 
following notation: for $X\subseteq S_a$ we set 
$\ol{r}_X=r^{\bfM.S_a}_X$ and 
$\ol{V}_X=V_X(\phi.S_a)$, 
as well as $\ol T_X = T_X(\phi.S_a)$.  
By Theorem~\ref{T:T-partition-of-contraction} 
the set $B'$ is a maximal 
proper T-flat of $\bfM$ inside $B$ and as in 
\altref{E:contraction-sequence} we have an 
equality $r_{B'\smsm A'} + \ol{r}_{A'}= r_{B'}$ 
and an exact sequence. 
\[
0
\lra V_{B'\smsm A'}\lra V_{B'}
\lra \ol{V}_{A'}
\lra 0. 
\]
Since $B'\smsm A'\subseteq B\smsm A$, 
if $B=A$ then $B'=A'$ 
and it is straightforward 
to verify from the definitions that the 
diagram \altref{E:commutativity-for-pi-sub-A}
is commutative. Our next goal is to show that 
this diagram commutes also in the remaining 
case $B=A\cup\{a\}$. This is obvious if 
$\ol T_A=0$, therefore until stated otherwise, 
in the sequel we will always assume 
that $B=A\cup\{a\}$ and that the T-flat $A$ 
is connected in $\bfM.S_a$.

Now we look at the following commutative 
diagram with exact rows: 
\[
\begin{CD}
0 
@>>>  
V_{B'\smsm A'} 
@>>> 
V_{B'}
@>\pi >> 
\ol{V}_{A'}
@>>> 0 
\\ 
@.       
@VVV               
@VVV       
@VVV           
@. 
\\  
0 
@>>>  
V_{B\smsm A}  
@>>> 
V_B    
@>\pi >>  
\ol{V}_A  
@>>> 0. 
\end{CD}
\] 
In it all the vertical maps are canonical 
inclusions.  Furthermore, the cokernel of the 
middle one is precisely 
$K_{BB'}$, and we write 
$\ol{K}_{AA'}= \ol{V}_A/\ol{V}_{A'}$ 
for the cokernel of the inclusion to the right. 
Thus $\pi$ induces a canonical map 
$\ol\pi\: K_{BB'}\lra \ol{K}_{AA'}$. 
Also, since we are assuming $B=A\cup\{a\}$, we 
have $V_{B\smsm A}=V_a$. Since 
either $B'=A'$ or $B'=A'\cup\{a\}$, we treat 
these two cases separately. 

Assume first that $B'=A'\cup\{a\}$. Then  
we have the canonical isomorphisms $\ol\pi$ and 
\[
\bold{c'}\: 
\twedge \ol{V}^*_{A'}\otimes V^*_a \lra 
\twedge V^*_{B'},
\] 
and an equality $B\smsm B'=A\smsm A'$. 
Furthermore, since $A$ is connected the T-part 
$A\smsm A'$ is independent in $\bfM.S_a$, 
therefore the map $\pi$ induces an isomorphism 
from  $V_{B\smsm B'}$ to $\ol V_{A\smsm A'}$. 
We write $\pi'$ for the canonically induced 
inverse isomorphism  
\[
\pi'\: 
\twedge\ol V_{A\smsm A'} \lra 
\twedge V_{B\smsm B'}.
\] 
It is now routine to  
check that the canonically induced diagrams 
\begin{equation}\elabel{E:canonical-one}
\begin{CD} 
\twedge\ol{V}_A^*\otimes V_a^* 
@> \bold c  >> 
\twedge V_B^*  
\\ 
@V \bold{\ol a} \ \otimes 1 VV 
@VV \bold a V 
\\ 
\twedge\ol{K}^*_{AA'} \otimes 
\twedge\ol{V}_{A'}^*\otimes V_a^*  
@> \wedge\ol\pi^*\otimes \bold{c'} >>  
\twedge K^*_{BB'} \otimes 
\twedge\ol{V}_{B'}^* 
\end{CD}
\end{equation}
and   
\begin{equation}\elabel{E:canonical-two} 
\begin{CD}
\twedge U_A \otimes U_a 
@> \wedge >>    
\twedge U_B 
\\ 
@V \delta\otimes 1 VV 
@VV \delta V  
\\ 
\twedge U_{A\smsm A'} \otimes 
\twedge U_{A'}\otimes U_a  
@> 1\otimes\wedge >>  
\twedge U_{B\smsm B'} \otimes 
\twedge U_{B'} 
\end{CD}
\end{equation}
together with 
\begin{equation}\elabel{E:canonical-three} 
\begin{CD}
(\ol V_{A'}\cap\ol V_{A\smsm A'})^*\otimes
\twedge\ol K_{AA'}^* 
@> \pi^*\otimes\wedge\ol\pi^* >> 
(V_{B'}\cap V_{B\smsm B'})^*\otimes 
\twedge K_{BB'}^*      
\\ 
@V  \ol{\bf b} VV   
@VV {\bf b} V 
\\ 
\twedge\ol V_{A\smsm A'}^* 
@> \wedge\pi^* >> 
\twedge V_{B\smsm B'}^* 
\end{CD}
\end{equation} 
and 
\begin{equation}\elabel{E:canonical-four}
\begin{CD} 
\twedge\ol V_{A\smsm A'}^*\otimes 
\twedge\ol V_{A\smsm A'} 
@> \wedge\pi^*\otimes\wedge\pi' >> 
\twedge V_{B\smsm B'}^*\otimes
\twedge V_{B\smsm B'} 
\\ 
@V  \mu\otimes 1 VV 
@VV \mu\otimes 1 V 
\\ 
\Bbbk  @> = >> \Bbbk 
\end{CD}
\end{equation}
are commutative. 
Write $\ol{S}_A$ and $\ol{S}_{A'}$ for 
for $S_A(\phi . S_a)$ and $S_{A'}(\phi . S_a)$, 
respectively. Since the maps 
$\pi_n\: S_nW \lra S_n\ol W$ commute with 
the multiplication in the symmetric algebras, 
by taking duals we obtain a commutative diagram 
\begin{equation}\elabel{E:canonical-five}
\begin{CD} 
\ol S_A^* 
@> \pi_{S_a,A}^* >> 
S_B^* 
\\
@V \ol\Delta VV 
@VV \Delta V 
\\
(\ol V_{A'}\cap\ol V_{A\smsm A'})^*\otimes
\ol S_{A'}^* 
@> \pi^*\otimes \pi_{S_a,A'}^* >> 
(V_{B'}\cap V_{B\smsm B'})^*\otimes S_{B'}^*.  
\end{CD}
\end{equation}
Now we consider the following diagram:  
\begin{equation}\elabel{E:big-diagram}
\begin{CD} 
\ol{T}_A 
@> 1\otimes{\bf d} >> T_{AB}^{(1)} 
@> \tau >> T_{AB}^{(2)} 
@> \pi_{S_a,A}^*\otimes\wedge\otimes{\bf c} 
 > \phantom{
   \scriptsize{
   \pi^*\otimes\wedge\ol\pi^*
        \otimes 1\otimes \pi_{S_a,A'}^* 
        \otimes\wedge 
        \otimes{\bf c'} }} >  
T_B 
\\ 
@V \ol\Delta\otimes\bfd\otimes\ol{\bf a} VV          
@V \ol\Delta\otimes\bfd
            \otimes\ol{\bf a}\otimes 1 VV              
@VV \ol\Delta\otimes\bfd\otimes 1
             \otimes\ol{\bf a}\otimes 1 V              
@VV \Delta\otimes\bfd\otimes{\bf a} V   
\\ 
\ol{\mc Q}_{AA'}
@>> 1\otimes{\bf d} >         
\mc Q^{(1)} 
@>> \tau >                    
\mc Q^{(2)} 
@> \phantom{
   \scriptsize{
   \pi^*\otimes\wedge\ol\pi^*
        \otimes 1\otimes \pi_{S_a,A'}^* 
        \otimes\wedge 
        \otimes{\bf c'} }} 
 > \pi^*\otimes \pi_{S_a,A'}^*
        \otimes 1\otimes\wedge 
        \otimes\wedge\ol\pi^*
        \otimes{\bf c'}  >    
\mc Q_{BB'} 
\\ 
@V \tau VV          
@V \tau VV              
@VV \tau V              
@VV \tau V 
\\ 
\ol{\mc R}_{AA'}
@> 1\otimes{\bf d} >>         
\mc R^{(1)} 
@> \tau >>                    
\mc R^{(2)} 
@> \pi^*\otimes\wedge\ol\pi^*
        \otimes 1\otimes \pi_{S_a,A'}^* 
        \otimes\wedge\otimes{\bf c'} >>    
\mc R_{BB'} 
\\ 
@V \ol{\bf b}\otimes\ol\phi\otimes 1 VV          
@V \ol{\bf b}\otimes\ol\phi\otimes 1 VV              
@VV \ol{\bf b}\otimes\ol\phi\otimes 1 V              
@VV {\bf b}\otimes\phi\otimes 1 V 
\\ 
\ol{\mc S}_{AA'}
@> 1\otimes{\bf d} >>          
\mc S^{(1)} 
@> \tau  >>                    
\mc S^{(2)} 
@> \wedge\pi^*\otimes \pi' 
              \otimes \pi_{S_a,A'}^* 
              \otimes\wedge\otimes{\bf c'} 
 > \phantom{
   \scriptsize{
   \pi^*\otimes\wedge\ol\pi^*
        \otimes 1\otimes \pi_{S_a,A'}^* 
        \otimes\wedge 
        \otimes{\bf c'} }} >   
\mc S_{BB'}
\\ 
@V \mu\otimes 1 VV          
@V \mu\otimes 1 VV              
@VV \mu\otimes 1 V              
@VV \mu\otimes 1 V 
\\ 
\ol{T}_{A'}
@> 1\otimes{\bf d} >>              
T_{A'B'}^{(1)} 
@> \tau >>                         
T_{A'B'}^{(2)} 
@> \pi_{S_a,A'}^*\otimes\wedge
                 \otimes{\bf c'} 
 > \phantom{
   \scriptsize{
   \pi^*\otimes\wedge\ol\pi^*
        \otimes 1\otimes \pi_{S_a,A'}^* 
        \otimes\wedge\otimes{\bf c'} }} >  
T_{B'} 
\end{CD}
\end{equation}
where we use $\ol{\phantom{X}}$ \ to denote 
objects for $\bfM.S_a$, we use $\twedge X$ 
for the $d$-th exterior power of a vector space 
$X$ of dimension $d$, and   
\[
\begin{aligned}
\mc Q^{(1)} &=  
\ol{\mc Q}_{AA'}\otimes U_a\otimes V_a^* 
\\ 
&=  
(\ol V_{A'}\cap\ol V_{A\smsm A'})^*\otimes
\ol S_{A'}^*\otimes 
\twedge U_{A\smsm A'}\otimes
\twedge U_{A'}\otimes 
\twedge\ol K_{AA'}^*\otimes
\twedge\ol V_{A'}^* \otimes 
U_a\otimes V_a^*; 
\\ 
\mc Q^{(2)} &= 
(\ol V_{A'}\cap\ol V_{A\smsm A'})^*\otimes
\ol S_{A'}^*\otimes 
\twedge U_{A\smsm A'}\otimes
\twedge U_{A'}\otimes U_a\otimes 
\twedge\ol K_{AA'}^*\otimes
\twedge\ol V_{A'}^* \otimes V_a^*;  
\\[+8pt] 
\mc R^{(1)} &= 
\ol{\mc R}_{AA'}\otimes U_a\otimes V_a^* 
\\
&= 
(\ol V_{A'}\cap\ol V_{A\smsm A'})^*\otimes
\twedge\ol K_{AA'}^*\otimes
\twedge U_{A\smsm A'}\otimes
\ol S_{A'}^*\otimes\twedge U_{A'}\otimes 
\twedge\ol V_{A'}^* \otimes U_a\otimes V_a^*;  
\\ 
\mc R^{(2)} &= 
(\ol V_{A'}\cap\ol V_{A\smsm A'})^*\otimes
\twedge\ol K_{AA'}^*\otimes
\twedge U_{A\smsm A'}\otimes
\ol S_{A'}^*\otimes\twedge U_{A'}\otimes 
U_a\otimes\twedge\ol V_{A'}^* \otimes V_a^*;   
\\[+8pt] 
\mc S^{(1)} &= 
\ol{\mc S}_{AA'}\otimes U_a\otimes V_a^* 
\\
&= 
\twedge\ol V_{A\smsm A'}^*\otimes
\twedge\ol V_{A\smsm A'}\otimes 
\ol S_{A'}^*\otimes\twedge U_{A'}\otimes 
\twedge\ol V_{A'}^*\otimes U_a\otimes V_a^*;  
\\ 
\mc S^{(2)} &= 
\twedge\ol V_{A\smsm A'}^*\otimes
\twedge\ol V_{A\smsm A'}\otimes 
\ol S_{A'}^*\otimes\twedge U_{A'}\otimes 
U_a\otimes\twedge\ol V_{A'}^*\otimes V_a^*. 
\end{aligned}
\]
Since the second column (from the left) is 
obtained from the first column by tensoring 
with $U_a\otimes V_a^*$, the squares that 
involve only these two columns commute. 
The squares that involve only the second 
and third columns commute for trivial 
reasons. Next, look at the squares that 
involve only the third and the last column. 
The top one commutes because the diagrams 
\altref{E:canonical-five}, 
\altref{E:canonical-one}, and 
\altref{E:canonical-two} commute. 
The second one (from the top) 
commutes for trivial reasons. 
The third and the last squares  
commute because the canonical diagrams 
\altref{E:canonical-three} and 
\altref{E:canonical-four}  
commute. Thus the diagram 
\altref{E:big-diagram} is 
commutative. Note that the composition 
of the maps in its top and bottom rows give 
$(\pi.S_a)_A$ and $(\pi.S_a)_{A'}$, respectively. 
Also, the composition of the maps in the first 
and the last columns give $(\phi.Y)_n^{AA'}$ 
and $\phi_n^{BB'}$ respectively. 
Since $|A'|+1=|B'|$, this yields the desired 
commutativity of 
\altref{E:commutativity-for-pi-sub-A} in the 
case when $B'=A'\cup\{a\}$. 

Now we turn to the only other 
possible case: when $B'=A'$. Then  
the map $\pi\: V_{B'}\lra \ol V_{A'}$ is an 
isomorphism, and the map $\ol\pi$ induces 
as usual a 
canonical isomorphism  
\[
\bold{c''}\: 
\twedge\ol{K}^*_{AA'}\otimes V^*_a \lra 
\twedge K^*_{BB'}.
\]  
Furthermore, $\pi$ induces a canonical isomorphism 
\[
\bold{b'}\: 
\twedge\ol V_{A\smsm A'}^*\otimes V_a^* \lra 
\twedge V_{B\smsm B'}^*,
\]
and, in conjunction with the exact sequence 
\[
0\lra U_a \overset{\phi}\lra V_{B\smsm B'}
           \overset{\pi}\lra\ol V_{A\smsm A'}\lra 0,
\]   
also the canonical isomorphism   
\[
\pi''\: \twedge\ol V_{A\smsm A'}\otimes U_a \lra 
        \twedge V_{B\smsm B'}.
\]
In particular, the diagram 
\begin{equation}\elabel{E:canonical-six} 
\begin{CD}
\twedge U_{A\smsm A'}\otimes U_a 
@> \wedge >> 
\twedge U_{B\smsm B'} 
\\
@V \wedge\ol\phi\otimes 1 VV 
@VV \wedge\phi V 
\\
\twedge\ol V_{A\smsm A'}\otimes U_a 
@> \pi'' >> 
\twedge V_{B\smsm B'}
\end{CD}
\end{equation} 
is commutative. Furthermore, it is  routine to  
check that the canonically induced diagrams 
\begin{equation}\elabel{E:canonical-one-bis}
\begin{CD} 
\twedge\ol{V}_A^*\otimes V_a^* 
@> = >> 
\twedge\ol{V}_A^*\otimes V_a^* 
@> \bold c >>     
\twedge V_B^* 
\\ 
@V \bold{\ol a} \ \otimes 1  VV  
@.  
@VV (-1)^{(\ol r_{A'})}\bold a V 
\\ 
\twedge\ol{K}^*_{AA'} \otimes 
\twedge\ol{V}_{A'}^*\otimes V_a^*  
@> \tau >> 
\twedge\ol{K}^*_{AA'} \otimes V_a^* \otimes  
\twedge\ol{V}_{A'}^*  
@> \bold{c''}\otimes\wedge\pi^* >>  
\twedge K^*_{BB'} \otimes 
\twedge\ol{V}_{B'}^* 
\end{CD}
\end{equation}
and   
\begin{equation}\elabel{E:canonical-two-bis} 
\begin{CD}
\twedge U_A \otimes U_a 
@> = >>    
\twedge U_A \otimes U_a 
@> \wedge >>    
\twedge U_B 
\\ 
@V \bfd\otimes 1 VV  
@. 
@VV (-1)^{|A'|}\bfd V  
\\ 
\twedge U_{A\smsm A'} \otimes 
\twedge U_{A'}\otimes U_a  
@> \tau >>  
\twedge U_{A\smsm A'} \otimes U_a\otimes 
\twedge U_{A'}  
@> \wedge\otimes 1 >>  
\twedge U_{B\smsm B'} \otimes 
\twedge U_{B'} 
\end{CD}
\end{equation}
together with 
\begin{equation}\elabel{E:canonical-three-bis} 
\begin{CD}
(\ol V_{A'}\cap\ol V_{A\smsm A'})^*\otimes
\twedge\ol K_{AA'}^* \otimes V_a^*  
@> \pi^*\otimes{\bf c''} >> 
(V_{B'}\cap V_{B\smsm B'})^*\otimes 
\twedge K_{BB'}^*      
\\ 
@V  \ol{\bf b}\otimes 1 VV   
@VV {\bf b} V 
\\ 
\twedge\ol V_{A\smsm A'}^* \otimes V_a^* 
@> {\bf b'} >> 
\twedge V_{B\smsm B'}^* 
\end{CD}
\end{equation} 
and 
\begin{equation}\elabel{E:canonical-four-bis}
\begin{CD} 
\twedge\ol V_{A\smsm A'}^*\otimes 
\twedge\ol V_{A\smsm A'}\otimes U_a\otimes V_a^*  
@>  \mu\otimes{\bf d}^{-1} >>  
\Bbbk 
\\ 
@V \tau VV  
@VV = V   
\\
\twedge\ol V_{A\smsm A'}^*\otimes V_a^* \otimes  
\twedge\ol V_{A\smsm A'}\otimes U_a   
@. 
\Bbbk 
\\ 
@V {\bf b'}\otimes\pi'' VV  
@VV = V 
\\ 
\twedge V_{B\smsm B'}^*\otimes\twedge V_{B\smsm B'} 
@> \mu >>  
\Bbbk 
\end{CD}
\end{equation}
are commutative. We write ${\bf q_1}$ (${\bf q_2}$, 
respectively) for the composition of the maps in the 
bottom row of \altref{E:canonical-one-bis} (of   
\altref{E:canonical-two-bis}, respectively).     
Now we consider the 
diagram 
\begin{equation}\elabel{E:big-diagram-bis}
\begin{CD} 
\ol{T}_A 
@> 1\otimes{\bf d} >> 
T_{AB}^{(1)} 
@> \tau >> 
T_{AB}^{(2)} 
@> \pi_{S_a,A}^*\otimes\wedge\otimes{\bf c} 
 > \phantom{
   \scriptsize{
   \pi^*\otimes\wedge\ol\pi^*
        \otimes 1\otimes \pi_{S_a,A'}^* 
        \otimes\wedge 
        \otimes{\bf c'} }} >  
T_B 
\\ 
@V \ol\Delta\otimes\bfd\otimes\ol{\bf a} VV          
@VV \ol\Delta\otimes\bfd\otimes\ol{\bf a}
                        \otimes 1 V              
@VV \ol\Delta\otimes\bfd\otimes 1
             \otimes\ol{\bf a}\otimes 1 V              
@V (-1)^n\Delta\otimes\bfd\otimes{\bf a} VV   
\\ 
\ol{\mc Q}_{AA'}
@>> 1\otimes{\bf d} >         
\mc Q^{(1)} 
@>> \tau >                    
\mc Q^{(2)} 
@> \phantom{
   \scriptsize{
   \pi^*\otimes\wedge\ol\pi^*
        \otimes 1\otimes \pi_{S_a,A'}^* 
        \otimes\wedge 
        \otimes{\bf c'} }} 
 > \pi^*\otimes \pi_{S_a,A'}^*
        \otimes{\bf q_2}\otimes{\bf q_1} >  
\mc Q_{BB'} 
\\ 
@V \tau VV          
@VV \tau V              
@VV \tau V              
@V \tau VV 
\\ 
\ol{\mc R}_{AA'}
@> 1\otimes{\bf d} >>         
\mc R^{(1)} 
@> \tau >>                    
\mc R^{(3)} 
@> \pi^*\otimes{\bf c''}
        \otimes\wedge 
        \otimes \pi_{S_a,A'}^* 
        \otimes 1 
        \otimes\wedge\pi^* >>    
\mc R_{BB'} 
\\ 
@V \ol{\bf b}\otimes\ol\phi\otimes 1 VV          
@VV \ol{\bf b}\otimes\ol\phi\otimes 1 V              
@VV \ol{\bf b}\otimes 1\otimes\ol\phi
                       \otimes 1 V              
@V {\bf b}\otimes\phi\otimes 1 VV 
\\ 
\ol{\mc S}_{AA'}
@> 1\otimes{\bf d} >>          
\mc S^{(1)} 
@> \tau  >>                    
\mc S^{(3)} 
@> {\bf b'}\otimes \pi'' 
           \otimes \pi_{S_a,A'}^* 
           \otimes 1  
           \otimes\wedge\pi^*  
 > \phantom{
   \scriptsize{
   \pi^*\otimes\wedge\ol\pi^*
        \otimes 1\otimes \pi_{S_a,A'}^* 
        \otimes\wedge 
        \otimes{\bf c'} }} >   
\mc S_{BB'}
\\ 
@V \mu\otimes 1 VV          
@VV \mu\otimes 1\otimes{\bf d}^{-1} V              
@. 
@V \mu\otimes 1 VV 
\\ 
\ol{T}_{A'}
@> = >>              
T_{A'B'}^{(1)} 
@> = >>                         
T_{A'B'}^{(2)} 
@> \pi_{S_a,A'}^*\otimes 1\otimes\wedge\pi^* 
 > \phantom{
   \scriptsize{
   \pi^*\otimes\wedge\ol\pi^*
        \otimes 1\otimes \pi_{S_a,A'}^* 
        \otimes\wedge 
        \otimes{\bf c'} }} >  
T_{B'},  
\end{CD}
\end{equation}
where 
\[
\begin{aligned}
\mc R^{(3)} &= 
(\ol V_{A'}\cap\ol V_{A\smsm A'})^* \otimes
\twedge\ol K_{AA'}^*\otimes V_a^*\otimes 
\twedge U_{A\smsm A'}\otimes U_a \otimes 
\ol S_{A'}^*\otimes\twedge U_{A'}\otimes 
\twedge\ol V_{A'}^*;   
\\
\mc S^{(3)} &= 
\twedge\ol V_{A\smsm A'}^*\otimes V_a^*\otimes 
\twedge\ol V_{A\smsm A'}\otimes U_a \otimes 
\ol S_{A'}^*\otimes\twedge U_{A'}\otimes 
\twedge\ol V_{A'}^*. 
\end{aligned}
\]
This diagram differs in several places from 
\altref{E:big-diagram}. We will show that it is 
commutative. The squares that involve the first two 
columns (from the left) are commutative because they 
are the same as in \altref{E:big-diagram}, except 
for the bottom one, which commutes for trivial reasons. 
Among the three squares that involve the second and third 
column, we only need to consider the one involving 
$\mc S^{(3)}$, but that one is commutative for trivial 
reasons as well. Next comes the parallelogram with 
corners in $\mc S^{(1)}$, $\mc S_{BB'}$, $T_{B'}$, and 
$T^{(1)}_{A'B'}$; it is commutative because of the 
commutativity of \altref{E:canonical-four-bis}.     
Finally, we deal with the three squares that involve 
the third and the fourth columns. The top one commutes 
because \altref{E:canonical-five}, 
\altref{E:canonical-one-bis}, and 
\altref{E:canonical-two-bis} commute, 
plus the fact that 
\[
(-1)^{\ol r_{A'}}(-1)^{|A'|}=(-1)^{|A'|-\ol r_{A'}}=
(-1)^n.
\]
The second 
square from the top commutes for trivial reasons, 
and the third square from the top 
commutes because  \altref{E:canonical-six} and 
\altref{E:canonical-three-bis}
are commutative. Therefore the diagram 
\altref{E:big-diagram-bis} is commutative. The desired 
commutativity of \altref{E:commutativity-for-pi-sub-A} 
is now immediate by looking at the outer rows and 
columns of \altref{E:big-diagram-bis}. The proof of 
the proposition is now complete. 
\end{proof}

\begin{proof}
[Proof of Theorem~\ref{T:contraction-embedding-is-morphism}]
Since each map $\pi_{Y,A}^*$ is injective 
(as the dual of a surjective map) the map 
$(\pi.Y)_n$ is injective for each 
$n\ge 0$. Thus it remains to show 
that $(\pi.Y)_\bullet$ is 
a morphism of complexes. Furthermore, by 
Theorem~\ref{T:contraction-embedding-composes} 
it suffices to prove our theorem under the assumption  
that $Y=S_a$ with $\{a\}$ an independent set in $\bfM$. 
So let $A$ be a T-flat in $\bfM.S_a$ of level $n\ge 1$. 
It suffices to show that 
\[
(\pi.S_a)_{n-1}\circ 
(\phi.S_a)_n\big|_{T_A(\phi.S_a)} = 
\phi_n \circ 
(\pi.S_a)_n\big|_{T_A(\phi.S_a)}. 
\]
Let $B=S\smsm (S_a\smsm A)^{\bsfC_{\bfM^*}}$, and 
observe that if $B''$ is a T-flat inside $B$ of 
level $n-1$ which is not of the form 
$S\smsm (S_a\smsm A')^{\bsfC_{\bfM^*}}$ for some 
T-flat $A'\subset A$ in $\bfM.S_a$ of level $n-1$, 
then by Theorem~\ref{T:T-partition-of-contraction} 
we must have that $(B\smsm B'')\cap S_a=\varnothing$. 
Therefore (since the T-part $B\smsm B''$ is nonempty) 
we must have that $B\smsm B''=\{a\}$, hence $\{a\}$ 
is a T-part of $B$ and $B''=A$. In particular, 
by Theorem~\ref{T:multiplicity-kernels} we obtain 
that the composition 
$\phi_n^{BB''}\circ\pi_{S_a,A}^*=0$. 
The desired conclusion is now immediate in 
view of the commutativity of 
\altref{E:commutativity-for-pi-sub-A}. 
\end{proof}

\section{Acyclicity}

In this section $\phi\: U_S \lra W$ is a 
representation over $\Bbbk$ of a matroid $\bfM$ 
on a finite set $S$.  We tie the remaining loose 
ends by presenting the proofs of 
Theorem~\ref{T:contraction-embedding}, 
Theorem~\ref{T:exactness-of-complex}, 
Theorem~\ref{T:independence-from-phi}, and 
Theorem~\ref{T:multiplicity-exactness}. 
We begin with a lemma. 

\begin{lemma}\mlabel{T:image-of-phi-sub-zero} 
The image of the augmentation map 
$\phi_0\: T_0(\phi)\lra U_S$ is precisely 
$\im(\phi_0)=\Ker(\phi)$. 
\end{lemma} 

\begin{proof} 
Let $I$ be a circuit of $\bfM$. Since $T_I$ is 
a $1$-dimensional vector space and for each $a\in I$ 
the map $\phi_0^{I,a}\: T_I \lra U_a$ is an 
isomorphism, the map 
$\phi_0\big|_{T_I}= - \sum_{a\in I} \phi_0^{I,a}$ 
is injective; therefore 
$0\ne \phi_0(T_I)\subseteq \Ker(\phi|I)$ and hence 
$\phi_0(T_I)=\Ker(\phi|I)$ because $\Ker(\phi|I)$ is 
$1$-dimensional due to the fact that $I$ is a circuit. 
Thus we obtain 
\[
\phi_0(T_I) = \sum_{I\in\mc T_0(\bfM)} \Ker(\phi|I)
\]
and so it suffices to show that 
$\Ker(\phi)=\sum_{I\in\mc T_0(\bfM)} \Ker(\phi|I)$. 
For each circuit $I$ let 
\[
v_I = \sum_{a\in I} d_{I,a} e_a 
\]
be a basis vector for the $1$-dimensional space 
$\Ker(\phi|I)$. By 
Remark~\ref{T:rank-1-intersection-b}(b) we have 
$d_{I,a}\ne 0$ for each $a\in I$. 
For a vector $v=\sum_{a\in S}c_ae_a$ in $U_S$ let 
$supp(v)=\{a\in S\mid c_a\ne 0\}$. Suppose that 
$v$ is a vector in $\Ker(\phi)$ but not in 
$K'=\sum_{I\in\mc T_0(\bfM)} \Ker(\phi|I)$ and 
of smallest support. Then $A=supp(v)$ is a 
dependent set in $\bfM$, hence contains a 
circuit $I$. Let $a$ be an element of $I$. 
But then the vector $w=d_{I,a}v - c_av_I$ 
has smaller support than $v$. 
\end{proof}

\begin{lemma}\mlabel{T:isomorphism-of-kernels} 
Let $Y$ be a subset of $S$ such that $S\smsm Y$ is 
independent in $\bfM$. Then the diagram 
\[
\begin{CD} 
T_0(\phi.Y)   
@> (\phi.Y)_0 >>  
U_Y 
\\ 
@V (\pi.Y)_0^\phi VV 
@AA {\bf proj.} A 
\\ 
T_0(\phi) 
@> \phi_0 >> 
U_S 
\end{CD}
\]
is commutative. 
\end{lemma}

\begin{proof} 
By Theorem~\ref{T:contraction-embedding-composes} 
it is enough to consider the case when $Y=S_a$ 
with $\{a\}$ independent in $\bfM$, so we 
assume this is the situation. 
Then it suffices to show that for each circuit $A$ 
in $\bfM.S_a$ and each element $c\in A$ the diagram 
\[
\begin{CD} 
T_A(\phi.S_a)   
@> (\phi.S_a)_0^{A,c} >>  
U_c  
\\ 
@V (\pi.S_a)_A^\phi VV 
@VV = V 
\\ 
T_B(\phi) 
@> \phi_0^{B,c} >> 
U_c 
\end{CD}
\]
is commutative, where 
$B=S\smsm (S_a\smsm A)^{\bsfC_{\bfM^*}}$. 
This is trivially true when $B=A$. In the 
only other possible case 
$B=A\cup\{a\}$, we consider the diagram 
\[
\begin{CD} 
\ol T_A 
@> \tau\circ(\bfd\otimes 1) >> 
\twedge U_{A\smsm c}\otimes 
\twedge \ol V_{A\smsm c}^*\otimes U_c 
@> \mu\circ(\wedge\phi\otimes 1) >> 
U_c 
\\ 
@V \tau\circ(1\otimes{\bf d}) VV 
@V \tau\circ(1\otimes{\bf d}) VV 
@VV = V 
\\
T_{AB}^{(2)} 
@> \tau\circ(\bfd\otimes 1) >> 
\twedge U_{A\smsm c}\otimes  
\twedge \ol V_{A\smsm c}^*\otimes U_a
        \otimes V_a^*\otimes U_c 
@> (\mu\otimes\mu)\circ
   (\wedge\phi\otimes 1\otimes 
          \phi\otimes 1) >> 
U_c 
\\ 
@V \wedge\otimes{\bf c} VV 
@V (\wedge\otimes{\bf c}) \circ \tau VV 
@VV = V 
\\ 
T_B 
@> \tau\circ(\bfd\otimes 1) >> 
\twedge U_{B\smsm c}\otimes 
\twedge V_{B\smsm c}^* \otimes U_c 
@> \mu\circ(\wedge\phi\otimes 1) >> 
U_c 
\end{CD} 
\]
where we use notation as in the proof of 
Theorem~\ref{T:contraction-morphism-for-S-a}. 
The two squares that involve the left column of 
this diagram commute for trivial reasons, and 
it is a routine exercise in multilinear algebra 
to verify that the two squares involving the 
right column also commute. Thus this diagram 
is commutative. The assertion of the lemma is 
now immediate. 
\end{proof}

\begin{proof} 
[Proof of Theorems~\ref{T:contraction-embedding}, 
\ref{T:exactness-of-complex}, and 
\ref{T:multiplicity-exactness}]
Note that by 
Lemma~\ref{T:image-of-phi-sub-zero}, 
Theorem~\ref{T:contraction-embedding-composes}, and 
Theorem~\ref{T:contraction-embedding-is-morphism},  
it is enough to 
prove by induction on $m=|S|$ that 
\begin{enumerate}
\item The complex $T_\bullet(\phi)$ is a 
      resolution of $\Ker(\phi)$; 
\item Theorem~\ref{T:multiplicity-exactness} holds; 
      and 
\item If $\{a\}$ is independent in $\bfM$ then 
      the morphism 
      $
      (\pi.S_a)^\phi_\bullet\: 
      T_\bullet(\phi.S_a)\lra 
      T_\bullet(\phi)
      $ 
      is an isomorphism in homology. 
\end{enumerate}

These statements are true trivially when $m=0$. 
More generally,   
when $\phi$ is a representation of a  
uniform matroid assertion 
(1) is true by Example~\ref{E:uniform-matroid-2};  
the matroid $\bfM.S_a$ is also uniform 
and therefore in view of Example~\ref{E:uniform-matroid-1} 
assertion (2) is just a well known property of 
symmetric powers; while (3) follows from (1) in view of 
Lemmas~\ref{T:isomorphism-of-kernels} 
and~\ref{T:image-of-phi-sub-zero},   
and the fact that by the independence of $\{a\}$   
the canonical projection $U_S\lra U_{S_a}$ 
induces an isomorphism 
of $\Ker(\phi)$ onto $\Ker(\phi.S_a)$. 
In particular, assertions (1), (2), and (3) are 
true also in the case $m=1$. 

Next, we assume that $m\ge 2$, that $\bfM$ is not 
uniform, and that our assertions 
are true for smaller values of $|S|$. Furthermore,  
by Theorem~\ref{T:T-bullet-of-sum},  
Theorem~\ref{T:multiplicity-spaces-for-restriction}, 
and Theorem~\ref{T:multiplicity-zero}  
we may assume that $S$ is connected  
(hence the unique largest T-flat is $S$ itself 
and therefore $\ell_S^{\bfM}\ge 1$ by the 
non-uniformity of $\bfM$). 

In addition, suppose $c\in S$ is such that $\{c\}$ 
is not a T-part of $S$. Then $\{c\}$ is independent 
and is not a T-part of any T-flat of $\bfM$, hence by 
Theorem~\ref{T:multiplicity-kernels} the canonical 
morphism 
$
(\pi.S_c)_\bullet^\phi\: 
T_\bullet(\phi.S_c)\lra T_\bullet(\phi)
$ 
is an isomorphism and therefore (1) is 
true by the induction hypothesis. Also, if $a=c$ 
then (3) is trivial. Furthermore, when $\{a\}$ is 
independent and $a\ne c$ then both $\{a\}$ is 
independent in $\bfM.S_c$ as well as $\{c\}$ is 
independent in $\bfM.S_a$ 
(otherwise $\{a,c\}$ is a circuit in $\bfM$ 
hence either $S=\{a,c\}$ which contradicts 
$\ell_S^{\bfM}\ge 1$, or 
$\{a,c\}$ is a T-part of $S$ which contradicts 
the connectedness of $S$, or 
$\{a,c\}$ is a disjoint union of at least two 
nonempty T-parts of $S$ which is ruled out because 
$\{c\}$ is not a T-part of $S$);    
thus by 
Theorem~\ref{T:contraction-embedding-composes} 
\[
\begin{CD} 
T_\bullet(\phi.S_{ac}) 
@> (\pi.S_{ac})_\bullet^{\phi.S_a} >> 
T_\bullet(\phi.S_a) 
\\
@V  (\pi.S_{ac})_\bullet^{\phi.S_c} VV 
@VV (\pi.S_a)_\bullet^\phi V 
\\
T_\bullet(\phi.S_c) 
@>  (\pi.S_c)_\bullet^\phi >> 
T_\bullet(\phi)
\end{CD}
\] 
is a commutative diagram, where 
$S_{ac}=S\smsm \{a,c\}$. It follows 
that $(\pi.S_a)_\bullet^\phi$ is an 
isomorphism in homology 
because the other three maps in that diagram are 
isomorphisms in homology (here 
we are applying  to the morphisms 
$(\pi.S_{ac})_\bullet^{\phi.S_a}$ and 
$(\pi.S_{ac})_\bullet^{\phi.S_c}$ our 
induction hypothesis). Thus (3) is true. Finally, 
let $\{a\}$ be independent in $\bfM$, and   
let $A$ be a T-flat of $\bfM.S_a$ such 
that $\{a\}$ is a T-part of 
$B=S\smsm (S_a\smsm A)^{\bsfC_{\bfM^*}}$. If 
$B\ne S$ then (2) holds for $A$ and $B$ by passing 
to $\phi|B$ and applying our induction hypothesis. 
If $B=S$ then necessarily $A=S_a$ and $\{a\}$ 
is a T-part of $S$. It follows that $a\ne c$ 
and that $\{a\}$ is an independent T-part of $S_c$ 
in $\bfM.S_c$. Then by 
Theorem~\ref{T:multiplicity-kernels} we have 
$S_{S_a}(\phi.S_a)\cong S_{S_{ac}}(\phi.S_{ac})$ 
and $S_S(\phi)\cong S_{S_c}(\phi.S_c)$ as well as  
$S_{S_a}(\phi)\cong S_{S_{ac}}(\phi.S_c)$, 
thus by  our induction hypothesis  
\[
\dim S_S(\phi) = 
\dim S_{S_a}(\phi.S_a) + \dim S_{S_a}(\phi)
\]
and therefore (2) holds in this last remaining 
case by Theorem~\ref{T:multiplicity-kernels}(2). 

Summarizing the last paragraph, if $S$ has a 
T-part that is not a singleton, 
then (1), (2), and (3) hold. 
Therefore in the sequel we will always assume that 
every T-part of the connected T-flat $S$ is 
(independent and) a singleton. 

Let $a\in S$. By our induction hypothesis the 
complex $T_\bullet(\phi.S_a)$ is a 
resolution of $\Ker(\phi.S_a)$,  
and for every T-flat $A$ in $\bfM.S_a$ such that 
$A\ne S_a$ and such that $\{a\}$ is a T-part of 
$B=S\smsm(S_a\smsm A)^{\bsfC_{\bfM^*}}$ the 
sequence 
\[
\begin{CD} 
0\lra S_A(\phi)\otimes V_a(\phi) 
\overset{\nu}\lra S_B(\phi)
@> \pi_{S_a,A}^\phi >> 
S_A(\phi.S_a) \lra 0 
\end{CD}
\] 
is exact. It remains to show that $T_\bullet(\phi)$ 
is a resolution of $\Ker(\phi)$, that the sequence 
\begin{equation}\elabel{E:exactness-to-prove}
\begin{CD} 
0\lra 
S_{S_a}(\phi)\otimes V_a(\phi) 
\overset{\nu}\lra 
S_S(\phi)
@> \pi_{S_a,S_a}^\phi >> 
S_{S_a}(\phi.S_a) \lra 0 
\end{CD} 
\end{equation} 
is exact, and that $(\pi.S_a)_\bullet^\phi$ is an 
isomorphism in homology. Consider the exact 
sequence of complexes 
\[
\begin{CD}
0\lra T_\bullet(\phi.S_a)
@> \ (\pi.S_a)_\bullet^\phi \ >> 
T_\bullet(\phi)\lra C \lra 0
\end{CD}
\] 
where by definition $C=\Coker(\pi.S_a)_\bullet^\phi$. 
Let $n=\ell_S^{\bfM}=\ell_{S_a}^{\bfM.S_a}$. 
Since $\bfM$ is not uniform and 
every T-part of $S$ is an independent singleton 
we have $n\ge 2$. Thus   
$\HH_{n-1}(T_\bullet(\phi.S_a))=0$, and since 
$T_n(\phi)=T_S(\phi)$ the map  
$\phi_n\: T_n(\phi) \lra T_{n-1}(\phi)$ 
is injective by 
Remark~\ref{T:injective-differential}. Thus  
the long homology exact sequence yields  
$\HH_n(C)=0$, and therefore the differential 
$c_n\: C_n \lra C_{n-1}$ is injective. On the other 
hand, 
\[ 
C_n = 
\bigl(
S_S^*(\phi)/\pi_{S_a,S_a}^*[S_{S_a}(\phi.S_a)^*]
\bigr)
\otimes\twedge U_S\otimes\twedge V_S^* 
\]
and $C_{n-1}$ is the direct sum of 
$
T_{S_a}(\phi)=
S_{S_a}(\phi)^*\otimes\twedge U_{S_a}
               \otimes\twedge V_{S_a}^*
$ 
with all components of the form 
\[
\bigl(
S_X(\phi)^*/
\pi_{S_a,X\smsm a}^*[S_{X\smsm a}(\phi.S_a)^*]
\bigr)
\otimes\twedge U_X\otimes\twedge V_X^* 
\] 
where $X$ is a T-flat in $S$ of level $n-1$ such 
that $\{a\}$ is a T-part of $X$. In fact, by 
Theorem~\ref{T:T-flats-on-S-a} and 
Theorem~\ref{T:multiplicity-kernels} in every 
degree $k$ the space $C_k$ is the direct 
sum $C_k'$ of all components 
of the form $T_Z(\phi)$ for $Z$ a T-flat in $\bfM$ of 
level $k$ such that $a\notin Z$ and $Z\cup\{a\}$ is a 
T-flat in $\bfM$ of level $k+1$; plus the direct sum 
$C_{k-1}''$ of all components of the form 
\[
C_X''=
\bigl(
S_X(\phi)^*/
\pi_{S_a,X\smsm a}^*[S_{X\smsm a}(\phi.S_a)^*]
\bigr)
\otimes\twedge U_X\otimes\twedge V_X^* 
\] 
where $X$ is a T-flat in $S$ of level $k$ such that 
$\{a\}$ is a T-part of $X$. In particular we have 
$C_n=C_{n-1}''$ and $C_{n-1}'=T_{S_a}(\phi)$. 
Furthermore, it is clear from the definitions 
that $C'=C'_\bullet$ is a subcomplex of $C$, 
and that the shifted down quotient complex 
$C''=(C/C')[1]$ has component 
$C_k''$ in homological degree $k$. 

Now, for each component 
$C_X''$ of $C_k''$, the  map 
\[
(-1)^{|X\smsm\{a\}|}\phi_k^{X,X\smsm\{a\}}\ \: \  
T_X(\phi)\lra T_{X\smsm\{a\}}(\phi) 
\] 
induces  by Theorem~\ref{T:multiplicity-kernels} a 
surjective homomorphism 
$\gamma_k^X\: C_X''\lra T_{X\smsm\{a\}}(\phi)$, 
and it is straightforward to verify that this induces 
a surjective morphism of complexes $\gamma\: C''\lra C'$ 
such that $C$ is precisely the mapping cone of $\gamma$. 
Furthermore, by our induction hypothesis $\gamma_k$ is 
in fact an isomorphism except possibly for $k=n-1$. 
It follows that $\HH_k(C)=0$ for $k\le n-1$, and 
$\HH_n(C)\cong\Ker(\gamma_{n-1}^S)$. 
Since we know that $\HH_n(C)=0$ 
(because $c_n$ is injective), 
we conclude as desired that the sequence 
\altref{E:exactness-to-prove} is exact, and  
that the complex $C$ is exact. Therefore,  
the map  $(\pi.S_a)_\bullet^\phi$ is a 
quasiisomorphism and the complex 
$T_\bullet(\phi)$ is a resolution of 
$\Ker(\phi.S_a)\cong \Ker(\phi)$.      
\end{proof}

\begin{proof} 
[Proof of Theorem~\ref{T:independence-from-phi}] 
In view of 
Theorem~\ref{T:multiplicity-spaces-for-restriction} 
and Theorem~\ref{T:multiplicity-zero} 
it suffices to consider the case when $A=S$ is 
a connected T-flat of $\bfM$. We induce on the 
level $n$ of $S$. The case $n=0$ is trivial. 
When $n\ge 1$ the assertion follows from 
Theorem~\ref{T:exactness-of-complex} due to 
the fact that the nonzero component of 
$T_\bullet(\phi)$ of highest homological degree 
is $T_n(\phi)=T_S(\phi)\cong S_S(\phi)^*$ and 
its dimension is the alternating sum of the 
dimensions of the components in lower 
homological degrees. The dimensions of these 
lower components however are invariants of 
$\bfM$ by our induction hypothesis.   
\end{proof}

\section{Appendix on matroids} 

The appendix contains some general results about matroids  
that will be needed at one point or another in our 
exposition. They are elementary in nature and well known. 
Due to the lack of appropriate references, and also for 
completeness, we have included their short proofs. 

Throughout this section $\bfM$ is a matroid on a 
finite set $S$, and all sets considered are subsets of $S$.  

\begin{lemma}\mlabel{T:kernels-inclusion} 
If $J\subseteq I$ then $|J|-r_J \le |I|-r_I$. 
\end{lemma}

\begin{proof} 
Let $J'$ be a maximal independent subset of $J$, and let 
$I'$ be an extension of $J'$ to a maximal independent subset 
of $I$. Thus $J'=I'\cap J$ hence 
\[
r_I-r_J=|I'\smsm J'|\le |I\smsm J| = |I|-|J|. 
\]
The desired conclusion is now immediate. 
\end{proof}

\begin{lemma}\mlabel{T:kernels-intersection}
Let $I$ and $J$ be sets such that 
\[
|I|-r_I=|J|-r_J=|I\cup J|- r_{I\cup J}= k 
\]
for some $k$. Then also $|I\cap J| - r_{I\cap J}=k$. 
\end{lemma}

\begin{proof} 
We have $|I\cap J|=|I| + |J| - |I\cup J|$ and 
$r_{I\cap J}\le r_I + r_J - r_{I\cup J}$. So 
\[
k\ge |I\cap J| - r_{I\cap J} 
 \ge |I| + |J| - |I\cup J| - r_I - r_J + r_{I\cup J} 
 \ge k + k - k = k,   
\]
where the first inequality is by Lemma~\ref{T:kernels-inclusion}. 
The desired conclusion is now immediate. 
\end{proof}

\begin{lemma}\mlabel{T:rank-one-intersection-base}
Let $J$ be an independent subset disjoint from $I$ such 
that $r_I + r_J = r_{I\cup J} + k$. If $I'$ is a maximal 
independent subset of $I\cup J$ such that $I'\cap I$ is 
a maximal independent subset of $I$, then 
$|I'\cap J|= |J| - k$. 
\end{lemma}

\begin{proof} 
We have the equalities $r_{I\cup J}=|I'|$ and 
$r_I=|I'\cap I|$. Since $J$ is independent we get 
\[
|J| - k = r_J - k = r_{I\cup J} - r_I 
         = |I'| - |I'\cap I| = |I'\cap J|,
\]
which is the desired conclusion. 
\end{proof}

\begin{lemma}\mlabel{T:rank-one-intersection-subsets}
Let $J$ be an independent set disjoint from $I$ and such that 
$r_J + r_I = r_{I\cup J} + k$. If $J'$ is a subset of $J$ 
and $I'$ is a subset of $I$ such that 
$r_{J'} + r_{I'} \ge r_{I'\cup J'} + k$, then we have 
$r_{J'} + r_{I'} = r_{I'\cup J'} + k$. 
\end{lemma}

\begin{proof} 
Take a maximal independent subset $I_1$ of $I'$ and extend 
it to a maximal independent subset $I_2$ of $I$. Extend $I_2$ 
to a maximal independent subset $I_3$ of $I\cup J'$, and then 
$I_3$ to a maximal independent subset $I_4$ of $I\cup J$. 
Let $I_3'= I_3\cap J'$. Note that $I_4\cap (I\cup J') = I_3$, 
so $I_4\cap J' = I_3'$. On the other hand 
$|I_4\cap J|=|J| - k$ 
by Lemma~\ref{T:rank-one-intersection-base}, hence 
$|I_3'|=|I_4\cap J'| \ge |J'| - k$; thus 
\[
r_{I'\cup J'}\ge |I_4\cap I'| + |I_4\cap J'| 
             \ge r_{I'} + |J'| - k = r_{I'} + r_{J'} - k 
             \ge r_{I'\cup J'}. 
\]
The desired conclusion is now immediate. 
\end{proof}

\begin{lemma}\mlabel{T:rank-one-intersection-splitting} 
Let $J$ be an independent set disjoint from $I$ and such that 
$r_I + r_J = r_{I\cup J} + k$. If $J'$ is a subset of $J$ 
and $I'$ is a subset of $I$ such that 
$r_{J'}+ r_{I'}= r_{I'\cup J'} + k$ then we have 
$r_{J\smsm J'} + r_{I\cup J'} = r_{I\cup J}$. 
\end{lemma}

\begin{proof} 
Let $I_2$ be a maximal independent subset of $I$ such that 
$I_1= I_2\cap I'$ is a maximal independent subset of $I'$. 

Let $I_3$ be any extension of $I_2$ to a maximal independent 
subset of $I\cup J$. We claim that $|I_3\cap J'|=|J'|-k$ and 
$I_3\cap (J\smsm J')= J\smsm J'$. 
Indeed, let $I_4=I_3\cap J'$. Since $I_3\cap I=I_2$, we have 
\[
|I_3\cap J|= |I_3| - |I_2| = r_{I\cup J} - r_I = r_J - k 
           = |J| - k, 
\] 
therefore $|I_4|\ge |J'| - k$; thus we get 
\[
r_{I'}+ r_{J'} - k = r_{I'\cup J'} \ge |I_3\cap(I'\cup J')| 
                   = |I_1| + |I_4| 
                   \ge r_{I'} + r_{J'} - k    
\]
whence the claimed equalities $|I_4|=|J'| - k$ and 
$I_3\cap (J\smsm J')= J\smsm J'$. 

Next, we prove that the set $I_3'=I_3\cap (I\cup J')$ is a maximal 
independent subset of $I\cup J'$. It is certainly independent, 
hence contained in a maximal independent set $I''$ of $I\cup J'$. 
Since $I_2\subset I_3'$, we have $I''\cap I = I_2 = I_3'\cap I$. 
Extend $I''$ to a maximal independent set $I_3''$ 
of $I\cup J$. By the claim from previous paragraph  
$|I_3''\cap J'|= |J'| - k$; therefore the inclusions 
\[
I_3\cap J' = I_3'\cap J' \subseteq   I''\cap J'
                         \subseteq I_3''\cap J' 
\]
together with the equality 
$|I_3\cap J'|=|J'| - k$ yield $I''\cap J'=I_3'\cap J'$. Thus 
$I_3'=I''$, hence $I_3'$ is a maximal independent subset 
of $I\cup J'$. 

Finally, since $J\smsm J' \subseteq I_3$, we observe that 
\[
r_{I\cup J} = |I_3| = |J\smsm J'| + |I_3\cap (I\cup J')| 
                    = r_{J\smsm J'} + r_{I\cup J'}, 
\]
which is the desired conclusion.     
\end{proof}

\begin{lemma}\mlabel{T:rank-one-intersection-tower} 
Let $J\subseteq I_1 \subseteq I_2$ be sets with $J$ independent 
and such that 
\[
r_{I_2}- r_{I_2\smsm J}=r_{I_1} - r_{I_1\smsm J} = k.
\]  
If $I$ is a set with 
$I_1\subseteq I\subseteq I_2$, then $r_I- r_{I\smsm J}=k$. 
\end{lemma} 

\begin{proof} 
Let $I_1\subseteq I \subseteq I_2$. Choose a maximal independent 
set $B$ of $I_2$ such that the sets $B\cap (I_1\smsm J)$, 
$B\cap (I\smsm J)$, and $B\cap (I_2\smsm J)$ are maximal 
independent subsets of $I_1\smsm J$, $I\smsm J$, and $I_2\smsm J$, 
respectively. Let $B'=B\cap J$. Then 
\[
|B'| = |B|-|B\cap (I_2\smsm J)| = r_{I_2} - r_{I_2\smsm J} = k.
\] 
Furthermore $|B\cap I_1|= |B\cap (I_1\smsm J)| + |B'| 
= r_{I_1\smsm J} + k = r_{I_1}$, hence $B\cap I_1$ is a maximal 
independent subset of $I_1$. Then every element of $I_1$ is in 
the closure of $B\cap I_1$,  
and since every element of $I\smsm J$ is in the closure of 
$B\cap (I\smsm J)$, 
it follows that every element of $I=(I\smsm J)\cup I_1$ is in 
the closure of $B\cap I$. 
Therefore $B\cap I$ is a maximal 
independent subset of $I$, thus 
$r_I=|B\cap I|=|B\cap (I\smsm J)| + |B'| = r_{I\smsm J} + k$. 
\end{proof}

\begin{lemma}\mlabel{T:zero-intersection-ranks}
Let $J$ be an independent subset of $I$ such that 
$r_J + r_{I\smsm J} = r_I$. If $J'$ is a subset of 
$J$, then $r_{J'} + r_{I\smsm J'} = r_I$. 
\end{lemma}

\begin{proof} 
Assume $r_{J'} + r_{I\smsm J'} \ge r_I + 1$. Then, since 
$r_{J\smsm J'} + r_{J'}= r_J$, we would get 
\[
r_I=r_J + r_{I\smsm J} = r_{J'} + r_{J\smsm J'} + r_{I\smsm J} 
\ge r_{J'} + r_{I\smsm J'} \ge r_I +1, 
\]
which is impossible. 
\end{proof}

\begin{lemma}\mlabel{T:zero-intersection-bases}
Let $I_1$ and $I_2$ be sets with 
$r_{I_1}+r_{I_2}=r_{I_1\cup I_2}$, and let 
$J\subseteq I_1\cup I_2$. 

The set $J$ is a maximal 
independent subset of $I_1\cup I_2$ if and only if 
the sets $J_1=I_1\cap J$ and $J_2=I_2\cap J$ are maximal 
independent subsets of $I_1$ and $I_2$, respectively. 
\end{lemma}

\begin{proof} 
Assume first that $J$ is a maximal independent subset of 
$I=I_1\cup I_2$. Then the sets $J_1$ and $J_2$ are independent, 
and $r_{J_1} + r_{J_2} \ge r_J = r_I$. Therefore we have 
\[
r_I = r_{I_1} + r_{I_2} \ge r_{J_1} + r_{J_2} \ge r_I, 
\]
thus necessarily $r_{I_i}=r_{J_i}=|J_i|$ for $i=1,2$. 

Conversely, if $J_1$ and $J_2$ are maximal independent 
subsets of $I_1$ and $I_2$, respectively, then clearly 
every element of $I=I_1\cup I_2$ is in the closure of 
$J=J_1\cup J_2$.  
Therefore $J$ contains a maximal 
independent subset $J'$ of $I$. By the first part of the 
lemma we must have $J'\cap I_i=J_i$ for $i=1,2$, hence 
$J'=J$. 
\end{proof}

\begin{lemma}\mlabel{T:zero-intersection-subsets} 
Let $I_1$ and $I_2$ are sets with 
$r_{I_1}+r_{I_2}= r_{I_1\cup I_2}$. If $J_1\subseteq I_1$ 
and $J_2\subseteq I_2$ then $r_{J_1}+r_{J_2}=r_{J_1\cup J_2}$. 
\end{lemma}

\begin{proof} 
Let $J= J_1\cup J_2$, and 
let $J_i'$ be a maximal independent subset of $J_i$ for 
$i=1,2$. Then $r_{J_i}=|J_i'|$. Extend $J_i'$ to a maximal 
independent subset $I_i'$ of $I_i$. Then $I'=I_i'\cup I_2'$ 
is a maximal independent subset of $I=I_1\cup I_2$ by 
Lemma~\ref{T:zero-intersection-bases}, hence the set 
$J'=J_1'\cup J_2'\subseteq I'\cap J$ is an independent subset 
of $J$. Also note that $J_1'\cap J_2'$ is an independent 
subset of $I_1\cap I_2$, and since 
$0\le r_{I_1\cap I_2}\le r_{I_1}+r_{I_2}-r_{I_1\cup I_2}=0$, 
we must have $r_{I_1\cap I_2}=0$ and therefore 
$J_1'\cap J_2'=\varnothing$. 
Furthermore for each $i$ every element of $J_i$ is 
in the closure of the set $J_i'$,  
therefore every element of $J$ is in the closure of $J'$. 
Thus 
\[
r_J = |J'| = |J_1'| + |J_2'| = r_{J_1} + r_{J_2}, 
\]
which is the desired conclusion. 
\end{proof}

\begin{proposition}\mlabel{T:direct-sum-criterion} 
Let $I_1, \dots, I_k$ be subsets of $A$. Then 
$A=I_1\oplus\dots\oplus I_k$ if and only if 
$A=I_1\sqcup\dots\sqcup I_k$ and 
$r_A = r_{I_1} + \dots + r_{I_k}$. 
\end{proposition}

\begin{proof} 
The ``only if'' direction was already done in 
Section~\ref{S:matroids}. 
To prove the ``if'' direction, we induce on $k$. 
In view of the definition 
of direct sum, the case $k=2$ is a straightforward consequence of 
Lemma~\ref{T:zero-intersection-bases}. Thus we assume $k\ge 3$ 
and that the result is true for $k-1$. Let 
$A'= I_1\sqcup \dots\sqcup I_{k-1}$. Then $A=A'\sqcup I_k$, and 
we have 
\[
r_A\le r_{A'}+ r_{I_k}\le r_{I_1}+\dots + r_{I_k}=r_A. 
\] 
It follows that $r_A=r_{A'}+r_{I_k}$, and that 
$r_{A'}=r_{I_1}+\dots +r_{I_{k-1}}$.  
Since the direct sum operation is associative, 
our induction hypothesis concludes the proof. 
\end{proof} 

\begin{lemma}\mlabel{T:zero-intersection-circuits} 
Let $I$ and $J$ be subsets of $A$ such that $A=I\oplus J$.  
If $C$ is a connected subset of $A$ 
then either $C\subseteq I$ or $C\subseteq J$. 
\end{lemma}

\begin{proof} 
Let $I'=I\cap C$ and $J'=J\cap C$.  
By Proposition~\ref{T:direct-sum-criterion} and 
Lemma~\ref{T:zero-intersection-subsets} we have 
$r_{I'} + r_{J'} = r_{I'\cup J'} = r_C$,
hence $C=I'\oplus J'$. The desired conclusion is 
now immediate. 
\end{proof}

\begin{proposition}\mlabel{T:connected-components-decomposition} 
Let $A_1,\dots, A_k$ be the connected components of $A$. Then 
we have $A=A_1\oplus\dots\oplus A_k$. 
\end{proposition}

\begin{proof} 
The proposition is trivially true if $A$ is connected, hence we  
assume that $A$ is not connected. Then $A=I\oplus J$ for some 
nonempty subsets $I$ and $J$. Since each connected component 
of $A$ is either inside $I$ or inside $J$ 
(by Lemma~\ref{T:zero-intersection-circuits}) it follows that 
the connected components of $I$ together with the connected 
components of $J$ give the connected components of $A$.  
An elementary induction on the size of $A$ now completes the 
proof. 
\end{proof}

\end{document}